\documentclass[superscriptaddress,onecolumn,notitlepage,nofootinbib]{revtex4-1}

\usepackage{etoolbox}
\usepackage{amsmath,amssymb,amsthm}
\usepackage{subcaption}
\usepackage{tikz}
\usepackage{graphicx}
\usepackage{xcolor}
\usepackage[colorlinks=true, allcolors=blue]{hyperref}
\usepackage{mathtools}
\usepackage{MnSymbol}
\usepackage{algorithm} 
\usepackage[noend]{algpseudocode}
\usepackage{cleveref}

\makeatletter
\newcommand*{\declarecommand}{%
  \@star@or@long\declare@command
}
\newcommand*{\declare@command}[1]{%
  \provide@command{#1}{}%
  \renew@command{#1}%
}
\makeatother

\newtheorem{remark}{Remark}

\declarecommand{\u}{{\mathbf{u}}}
\declarecommand{\f}{{\mathbf{f}}}
\declarecommand{\x}{{\mathbf{x}}}
\declarecommand{\y}{{\mathbf{y}}}
\declarecommand{\n}{{\mathbf{n}}}
\declarecommand{\N}{{\mathbf{N}}}
\declarecommand{\X}{{\mathbf{X}}}
\declarecommand{\s}{{s}}
\declarecommand{\k}{{\mathbf{k}}}
\declarecommand{\btau}{{\boldsymbol{tau}}}
\declarecommand{\grad}{\nabla}
\declarecommand{\Id}{\mathbb{I}}
\declarecommand{\L}{\mathcal{L}}
\declarecommand{\B}{\mathcal{B}}
\declarecommand{\aOmega}{{\A}}
\declarecommand{\A}{\mathcal{A}}
\declarecommand{\I}{{\mathcal{I}}}
\declarecommand{\fOmega}{{\tilde\Omega}}
\declarecommand{\ur}{{\tilde u}}
\declarecommand{\ua}{{u_\A}}
\declarecommand{\uh}{{u_H}}
\declarecommand{\ui}{{u_I}}

\begin{document}

\title{Spectrally accurate solutions to inhomogeneous elliptic PDE in smooth geometries using function intension}
\author{David Stein}
\email{dstein@flatironinstitute.org}
\affiliation{Center for Computational Biology, Flatiron Institute, New York, NY 10010, USA}

\begin{abstract}
  We present a spectrally accurate embedded boundary method for solving linear, inhomogeneous, elliptic partial differential equations (PDE) in general smooth geometries, focusing in this manuscript on the Poisson, modified Helmholtz, and Stokes equations. Unlike several recently proposed methods which rely on function extension, we propose a method which instead utilizes function \emph{intension}, or the smooth truncation of known function values. Similar to those methods based on extension, once the inhomogeneity is truncated we may solve the PDE using any of the many simple, fast, and robust solvers that have been developed for regular grids on simple domains. Function intension is inherently stable, as are all steps in the proposed solution method, and can be used on domains which do not readily admit extensions. We pay a price in exchange for improved stability and flexibility: in addition to solving the PDE on the regular domain, we must additionally (1) solve the PDE on a small auxiliary domain that is fitted to the boundary, and (2) ensure consistency of the solution across the interface between this auxiliary domain and the rest of the physical domain. We show how these tasks may be accomplished efficiently (in both the asymptotic and practical sense), and compare convergence to several recent high-order embedded boundary schemes.
\end{abstract}

\maketitle

\section{Introduction}
\label{section:introduction}

Let $\L$ denote a constant coefficient elliptic operator, $\Omega$ a simply connected compact subset of $\mathbb{R}^2$ with smooth boundary $\Gamma$, and $\B$ a boundary operator (e.g. $\B u=u$ for Dirichlet boundary conditions). We assume that the specified inhomogeneity $f$ is smooth in $\Omega$ and the specified boundary condition $g$ is smooth on $\Gamma$. We seek to find a solution $u$ to the partial differential equation (PDE)
\begin{subequations}
  \begin{align}
    \L u &= f &&\textnormal{in }\Omega,  \\
    \B u &= g &&\textnormal{on }\Gamma.
  \end{align}
  \label{equation:general_pde}
\end{subequations}
There are two cases where optimal methods to this problem exist. The first of these is when $\Omega$ is geometrically simple: if $\Omega$ is the doubly periodic rectangle $\mathbb{T}^2$, spectral methods based on the Fast-Fourier transform (FFT) provide an optimal method for the inversion of $\L$; similar methods exist for several other simple geometries based on different spectral expansions \cite{trefethen1996finite} or adaptive integration using quadtrees \cite{greengard1996direct}. The second case is when $\Omega$ need not be geometrically simple, but $f=0$. In this case, well-conditioned boundary integral equation (BIE) methods exist for many commonly studied operators $\L$; spectrally accurate or high-order singular quadratures for the associated N\"ystrom schemes along with kernel-dependent and kernel-independent Fast Multipole Methods (FMMs) enable accurate solutions to be computed and evaluated in optimal time \cite{LIE,HW,Moura94,yingbeale,yanplatform,quaife2021hydrodynamics,sinha2016shape,nazockdast2017cytoplasmic,nazockdast2017fast,sorgentone2021numerical,HFMM2D,pvfmm}.

When the domain is complicated, as in the domain shown in \Cref{figure:domain_decomposition,figure:schematic}, and the problem is inhomogeneous (i.e. $f\neq0$) the situation is less clear. While methods have been continuously improving over the last several decades, all methods lack optimality in some way: whether due to slow convergence; ill-conditioning; or long compute times. Recently, considerable interest has focused on methods utilizing function extension, where either the inhomogeneity or the unknown solution defined on the general domains $\Omega$ is extended beyond its known values. These works include active penalty methods \cite{shirokoff2015sharp}, variations on the Immersed Boundary method \cite{stein2016immersed,stein2017immersed}, methods utilizing radial basis functions (PUX) \cite{fryklund2018partition,fryklund2020integral,af2020fast}, Fourier continuation methods \cite{bruno2010high,lyon2010high,bruno2020two,fontana2020fourier}, and those relying entirely on BIE \cite{askham2017adaptive}. Although promising, as these methods can produce relatively high-order discretizations with reasonable compute times, function extension is an inherently ill-conditioned process, as evinced by the great pains that some of these methods have taken to provide stability, and we worry that such methods will have issues both when steep boundary layers arise (as has been our own experience in the simulation of complex fluids \cite{stein2019convergent} and dissolution problems \cite{mac2021stable}), and when domains curve back on themselves.

We provide an alternative embedded boundary scheme that shares many of the benefits of those methods that utilize function extension, while eliminating some of the drawbacks. Rather than attempting to extend the inhomogeneity, we smoothly roll it off to $0$ \emph{inside of the domain}. Because this does not require extrapolation, it is both inherently \emph{stable} and relatively simple, requiring only the distance to the boundary and a regularized cutoff function. We then solve the PDE with this modified right-hand side using a regular grid method. Unfortunately, the story doesn't end there, as this candidate solution does not satisfy the PDE in the entire domain $\Omega$. To correct errors near to the boundary, we solve an \emph{annular} problem in a thin boundary-fitted annulus. The discrepancy between the annular solution and the solution on the regular grid is corrected, and boundary conditions are finally enforced by solving a homogeneous PDE utilizing a well-conditioned BIE method.

In this manuscript we present both the abstract method and a relatively simple implementation which makes use of global discretizations for both the boundary and regular grid. For problems that are not significantly multiscale in nature, this turns out to be reasonably performant, with virtually all steps having asymptotic scalings that are \emph{less} than the FFTs used to solve the regular grid problem. For simplicity, we discretize certain steps using methods that depend on dense linear algebra which have slightly worse scaling in the \emph{setup} stage of the problem, but not in the \emph{solution} stage of the problem; meaning that the overall scaling of solving repeated problems on the same geometry is \emph{the same as solving the regular PDE on a periodic grid} (albeit with worse constants). We make comments where poorly scaling methods are utilized, and how they could be improved upon (using methods already available in the literature); but in practice, these stages are rarely limiting for moderately sized problems as they make use of efficient BLAS and LAPACK routines.

This paper is organized as follows.  In \Cref{section:methods_overview}, we introduce the basic methodology of \emph{function intension} for solving PDE of the form given in \Cref{equation:general_pde}. This presentation will be simple and stripped down, both for pedagogical purposes and to emphasize the modular nature of the abstract method, where many of the substages have wide freedom of implementation, with little detailed interdependence. Then in \Cref{section:methods_preliminaries}, we introduce some basic preliminaries which will make discussion of our specific implementation easier. In \Cref{section:methods_specifics}, we revisit our abstract presentation of the method, now providing details for the specific implementational choices made throughout this paper. In \Cref{section:simple_poisson}, we explore how the method can be used to generate a solver of fixed algebraic order for any $M>0$, and show how that order can naturally be varied with the discretization to produce a spectrally accurate scheme. In \Cref{section:parameter_selection}, we discuss how to set the various parameters required by the method, and with these fixed, demonstrate large $N$ stability for the simple problem studied in the previous section. For this same problem, we show practical wall-clock timings across a wide range of problem sizes. To wrap up the presentation of the methodology, in \Cref{section:recombobulation} we provide algorithms for both the setup and solve components of our solver, with asymptotic scalings given for all compute-heavy steps.

We then turn to results for more complex problems set on more difficult domains. In \Cref{section:PUX_comparison}, we extend our method to multiply-connected domains, and compare our Poisson solver to the Partition of Unity Extension (PUX) method \cite{fryklund2018partition}. In \Cref{section:fourier_continuation}, we compare our method to the two-dimensional Fourier continuation method \cite{bruno2020two}, using this comparison to emphasize some differences between function intension and function extension.
In \Cref{section:modified_helmholtz}, we discretize and solve a modified Helmholtz problem, with a large parameter $\alpha$ as arises when time-discretizing the diffusion equation. We again compare our results to those produced by the PUX method \cite{fryklund2020integral}, finding very close agreement in the rate of convergence between the two methods, and similar stability as $\alpha$ grows. Finally, we solve a Stokes problem, comparing both errors and run-times with the IBSE method \cite{stein2017immersed}. In \Cref{section:discussion}, we conclude by discussing some of the outstanding issues with the method proposed here, and discuss some possible improvements to the method, from both a methodological and practical point of view.

\section{Methods --- overview}
\label{section:methods_overview}

\begin{figure}[h!]
	\centering
	\includegraphics[width=0.4\textwidth]{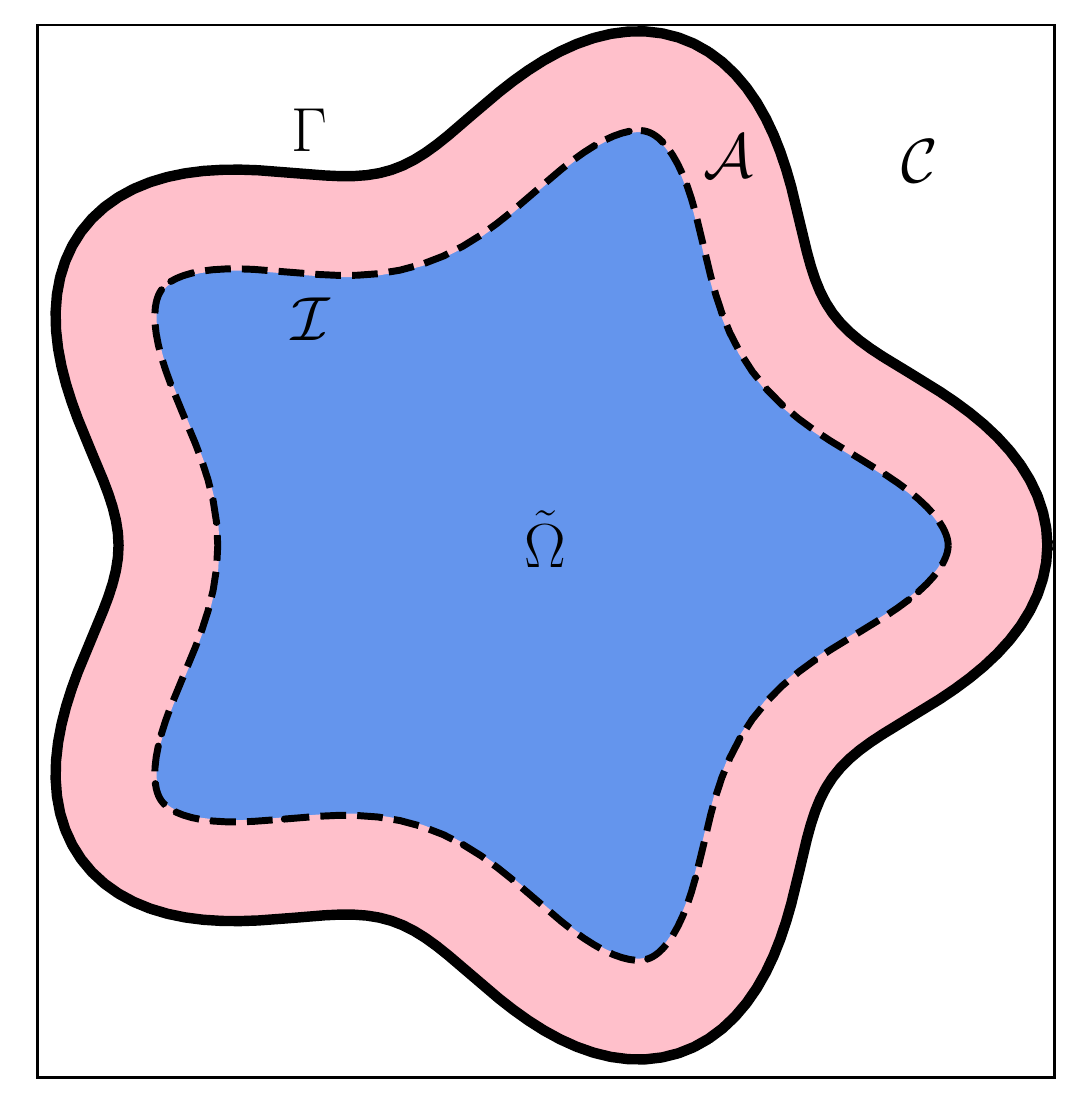}
	\caption{An example of a domain $\Omega$ (everything inside the solid black curve $\Gamma$, and the union of the pink and blue regions), embedded into a geometrically simple computational domain $\mathcal{C}$ (everything inside the rectangle that contains $\Omega$). The domain $\Omega$ is decomposed into the faithful domain $\fOmega$, shown in blue, and the annular domain $\A$, shown in pink.  The faithful domain $\fOmega$ and the annular domain $\A$ are separated by an interface $\I$, denoted by the dashed black line.}
	\label{figure:domain_decomposition}
\end{figure}

\begin{figure}[h!]
  \centering
  \begin{subfigure}[c]{0.3\textwidth}
    \centering
    \includegraphics[width=\textwidth]{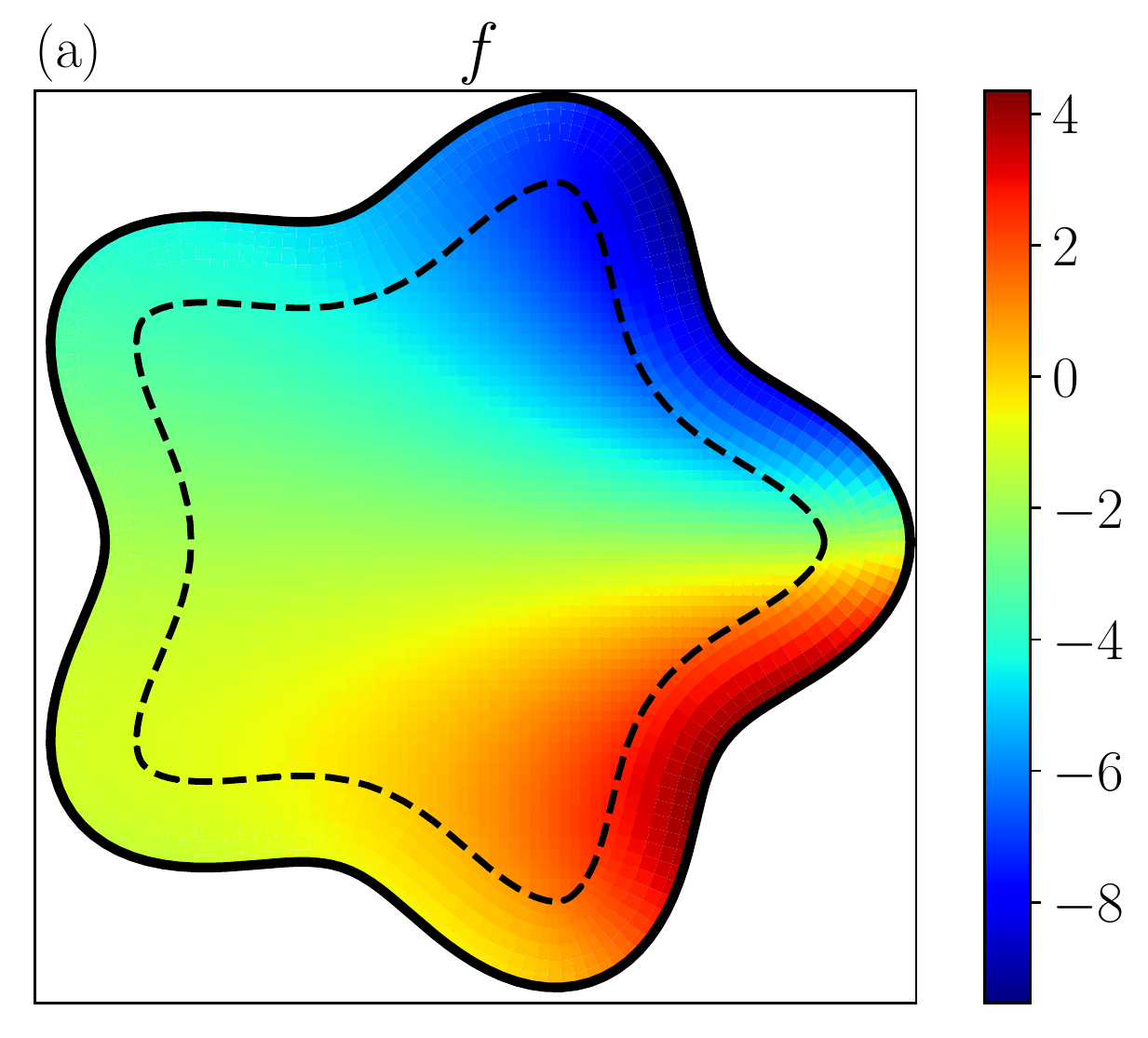}
  \end{subfigure}
  \hfill
  \begin{subfigure}[c]{0.3\textwidth}
    \centering
    \includegraphics[width=\textwidth]{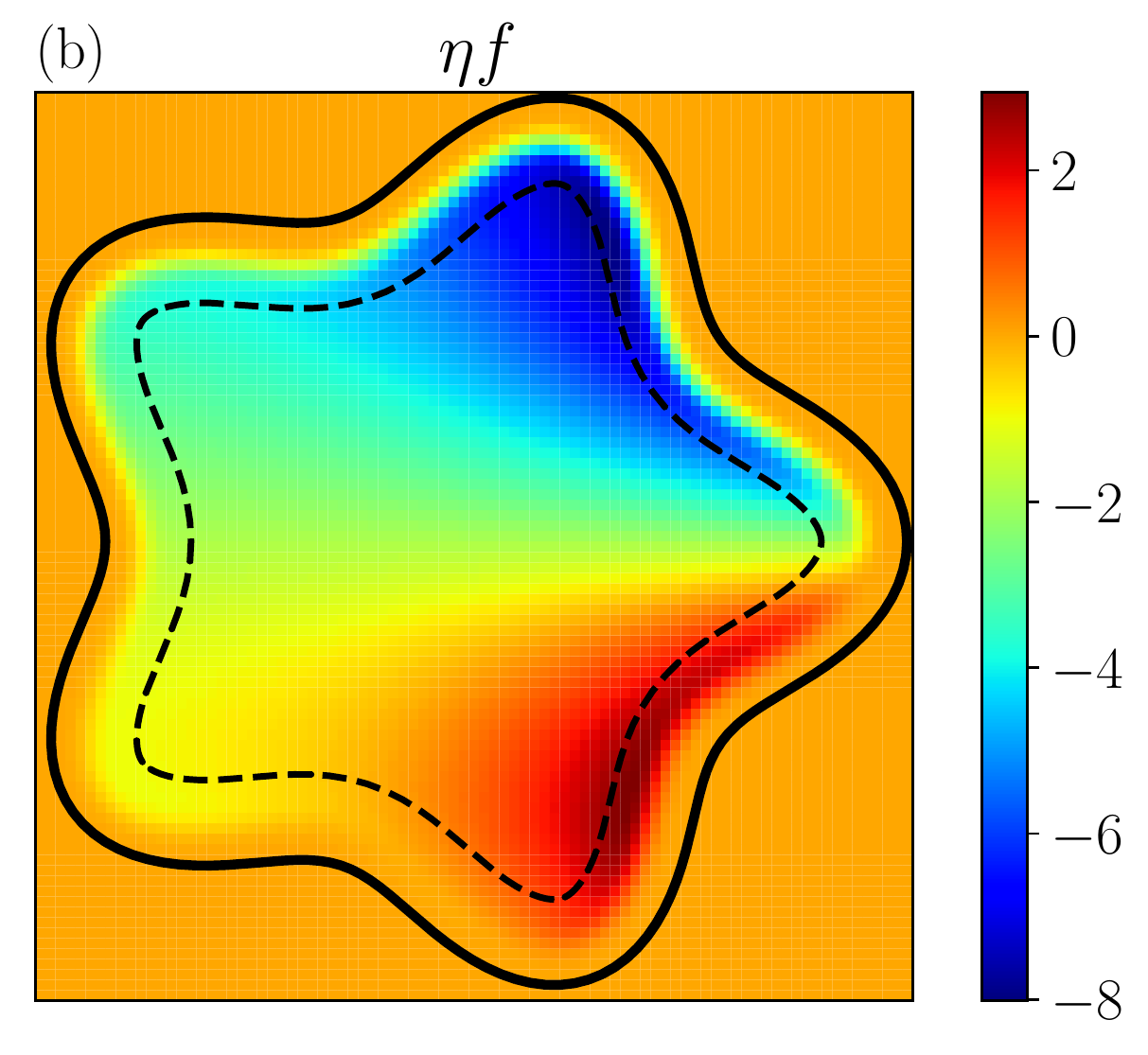}
  \end{subfigure}
  \hfill
  \begin{subfigure}[c]{0.3\textwidth}
    \includegraphics[width=\textwidth]{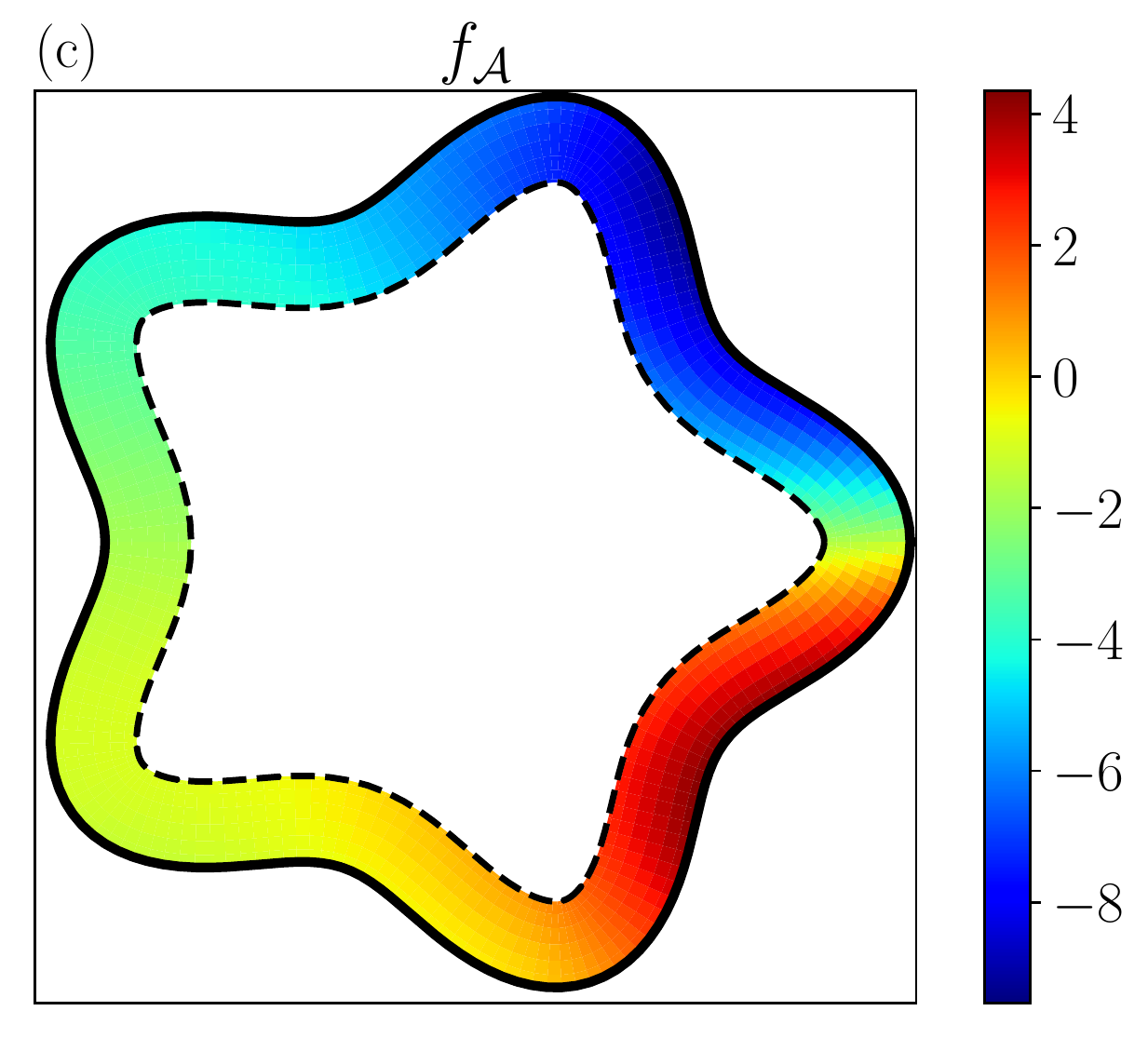}
  \end{subfigure}

  \begin{subfigure}[c]{0.3\textwidth}
    \centering
    \includegraphics[width=\textwidth]{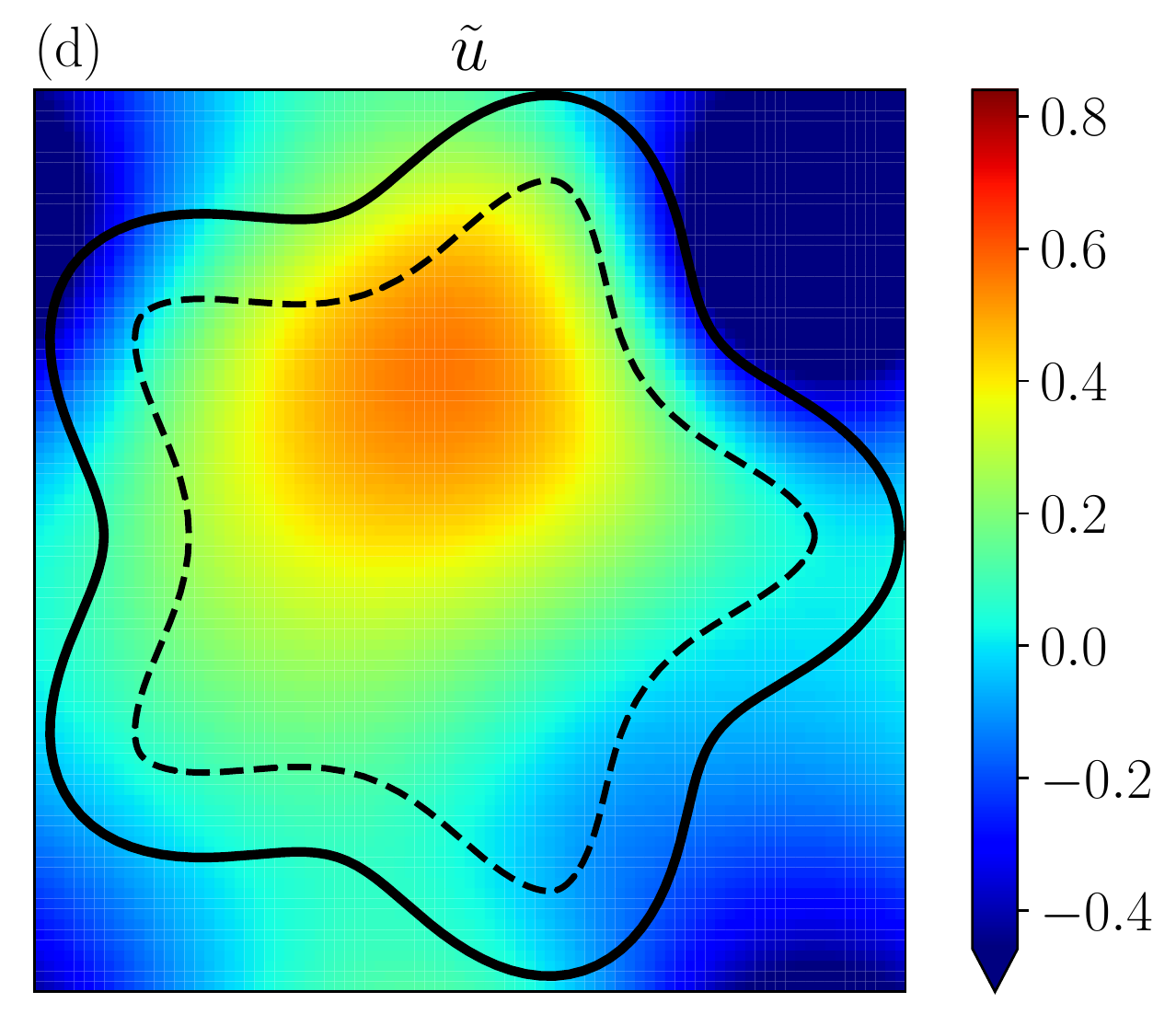}
  \end{subfigure}
  \hfill
  \begin{subfigure}[c]{0.3\textwidth}
    \centering
    \includegraphics[width=\textwidth]{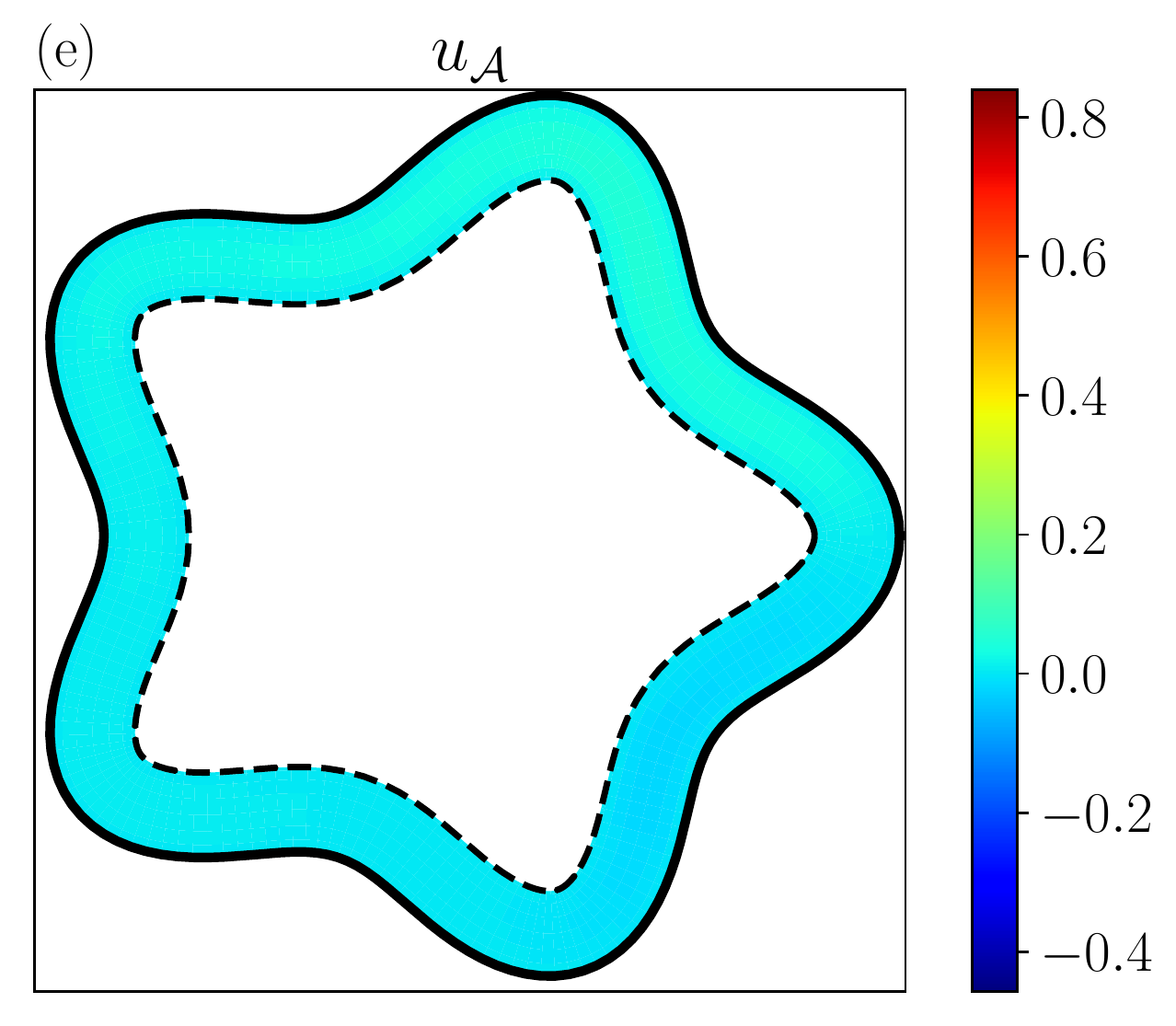}
  \end{subfigure}
  \hfill
  \begin{subfigure}[c]{0.3\textwidth}
    \centering
    \includegraphics[width=\textwidth]{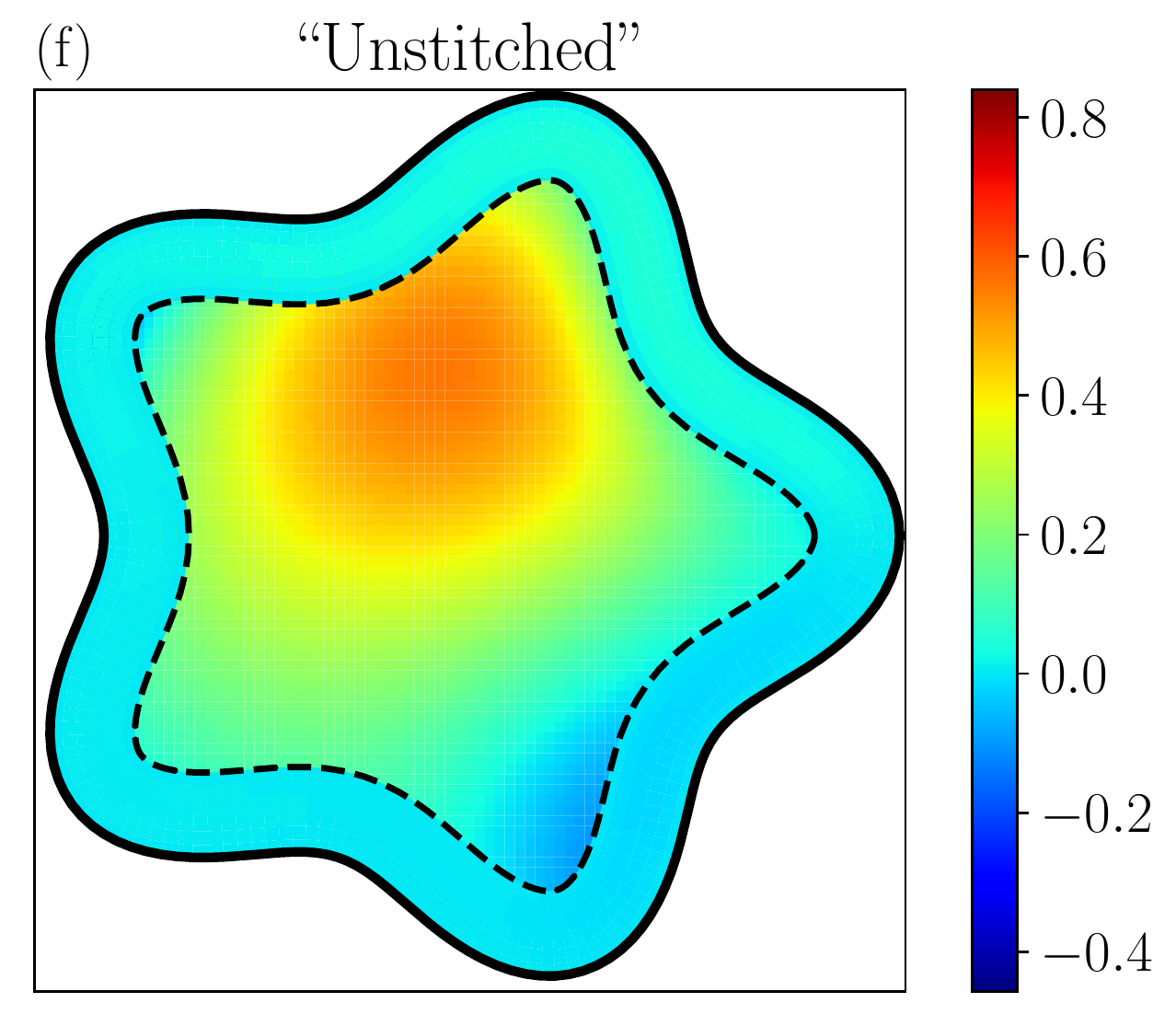}
  \end{subfigure}
  
  \begin{subfigure}[c]{0.3\textwidth}
    \includegraphics[width=\textwidth]{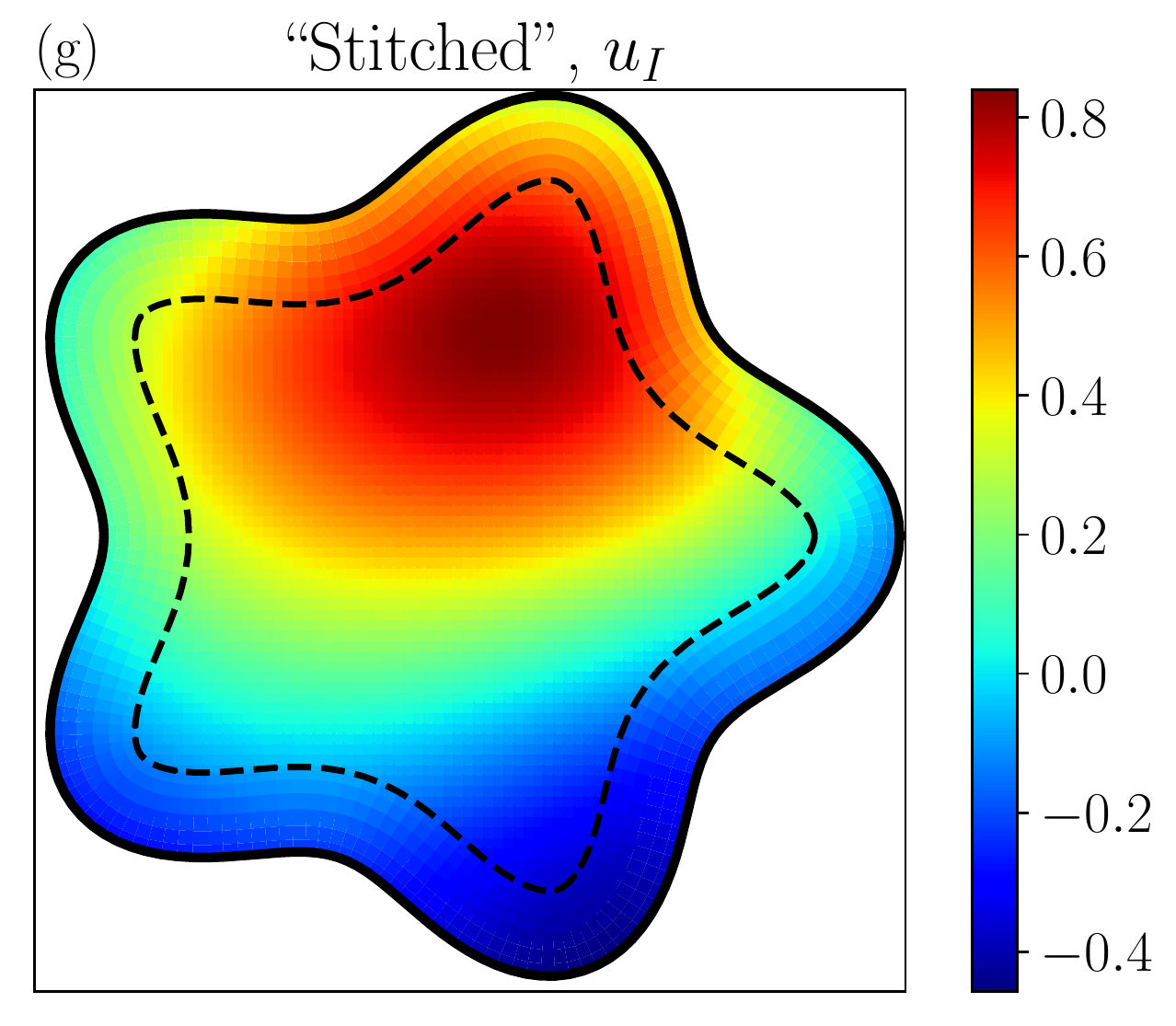}
  \end{subfigure}
  \hfill
  \begin{subfigure}[c]{0.3\textwidth}
    \includegraphics[width=\textwidth]{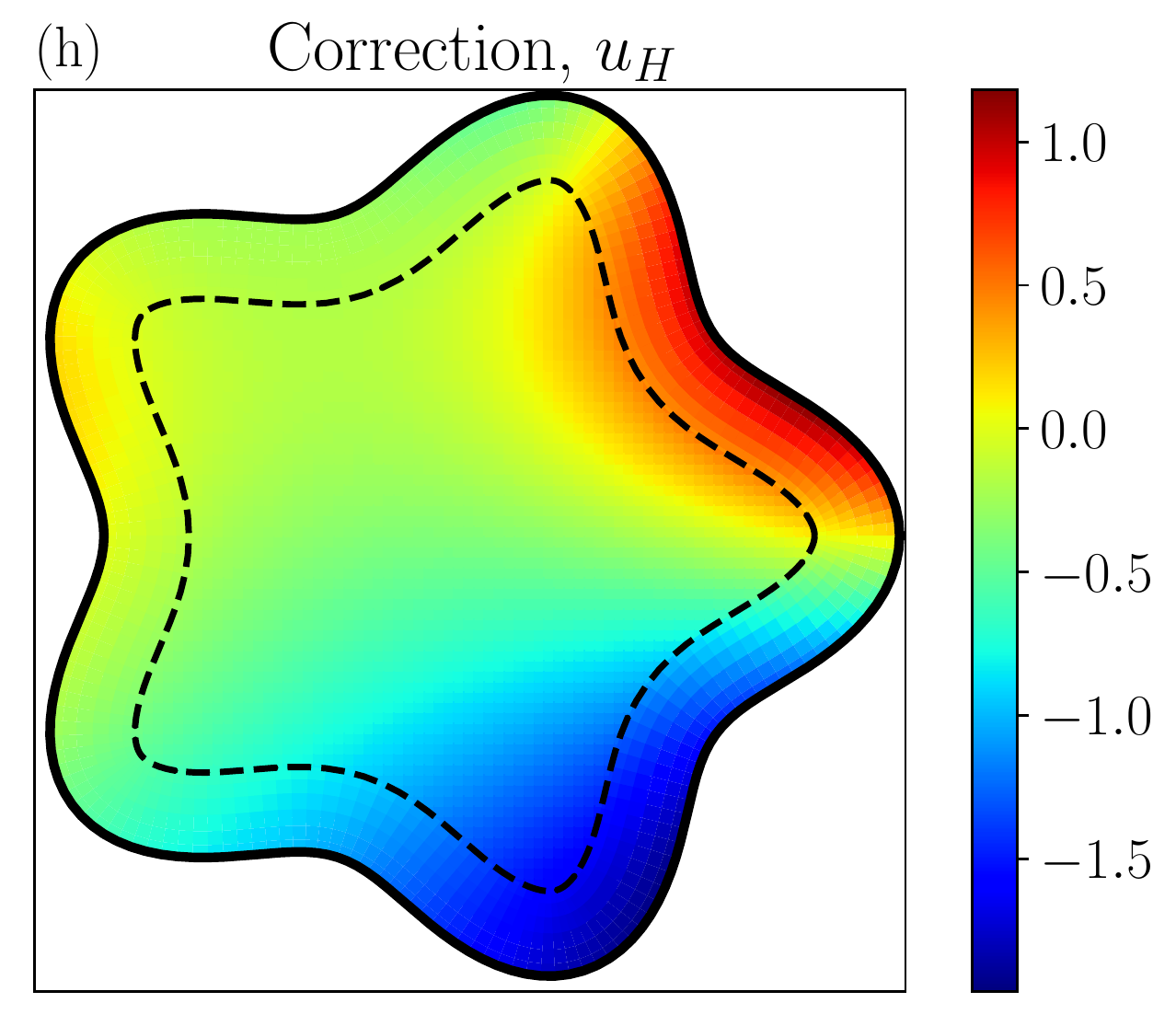}
  \end{subfigure}
  \hfill
  \begin{subfigure}[c]{0.3\textwidth}
    \includegraphics[width=\textwidth]{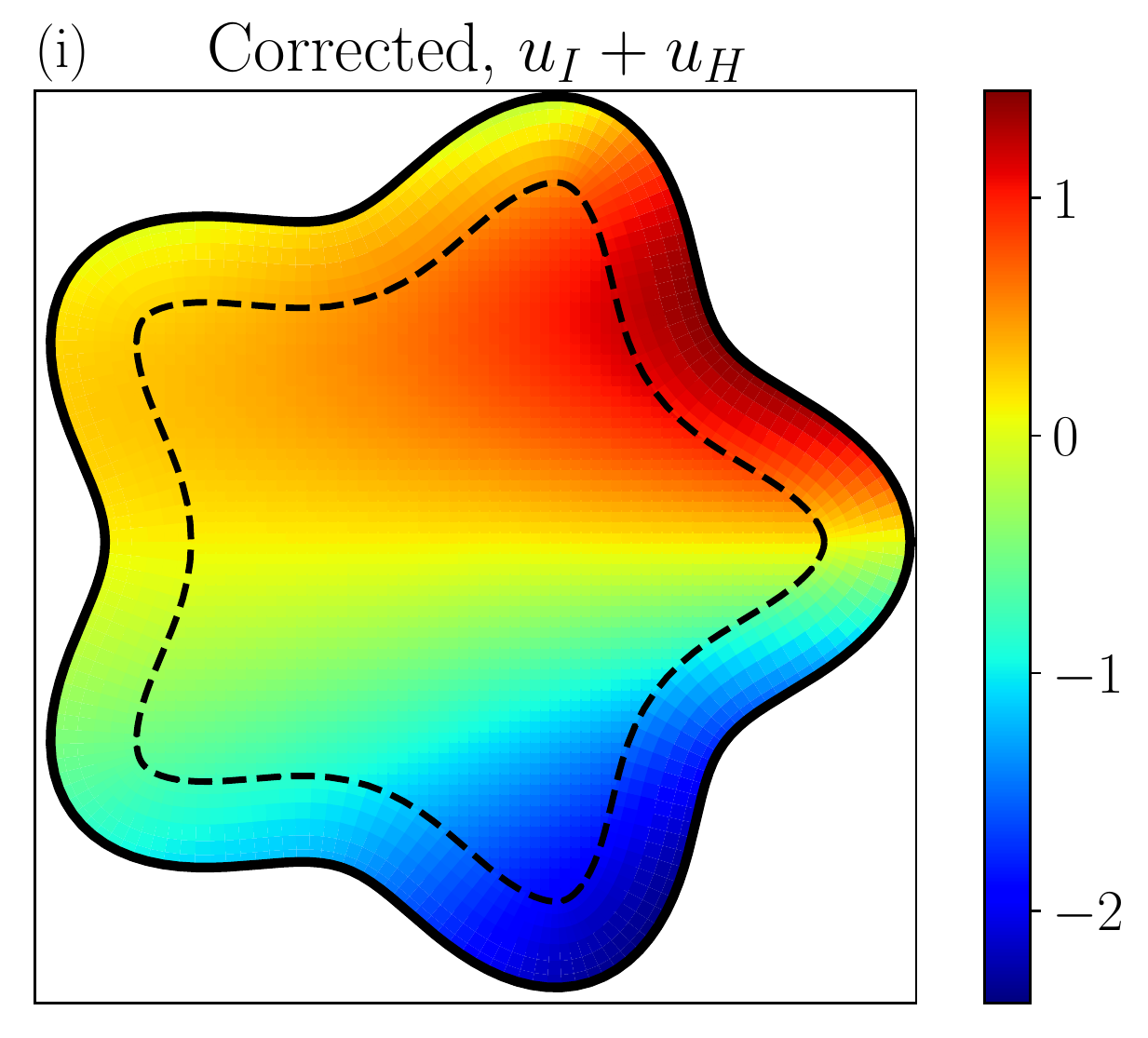}
  \end{subfigure}

  \caption{Schematic of solution process. The inhomogeneous forcing $f$ (panel a) is stably decomposed into a truncated function $\eta f$ (panel b) and an annular function $f_\A$ (panel c), via the process of \emph{function intension} (see \Cref{section:methods_overview:intension_regular}). Note that the function $\eta f=f$ everywhere inside the faithful domain $\fOmega$. The regular solution $\ur$ is obtained via standard methods on the simple domain $\mathcal{C}$ (panel d), while a thin-strip problem (see \Cref{section:methods_overview:annular_problem}) is solved to obtain $\ua$ (panel e). The ``unstitched solution'' (panel f) defined by truncating the regular solution $\ur$ to the faithful domain $\fOmega$, and by the annular solution $\ua$ inside $\A$, is clearly not smooth. Adding layer potentials, as described in \Cref{section:methods_overview:stiching}, fixes these discontinuities and produces the ``stitched'', or inhomogeneous solution $\ui$ (panel g). Correcting the boundary condition by solving a homogeneous problem (shown in panel h, see \Cref{section:methods_overview:homogeneous}) gives the estimated solution to the problem (panel i).}
  \label{figure:schematic}
\end{figure}

To avoid the main ideas being lost in a thicket of details, we begin with a highly simplified discussion of the overall method, deferring most questions of implementation to \Cref{section:methods_preliminaries,section:methods_specifics}. For further simplicity, let us choose as a specific case the interior Poisson problem with Dirichlet boundary conditions:
\begin{subequations}
	\begin{align}
		\Delta u &= f &&\textnormal{in }\Omega,	\\
		u				 &= g &&\textnormal{on }\Gamma,
	\end{align}
\end{subequations}
with $\mathcal{L}$ and $\mathcal B$ from \Cref{equation:general_pde} given by the Laplace operator and the interior trace operator for the boundary $\Gamma$, respectively. We begin by specifying a domain decomposition and defining several regions of space, boundaries, and interfaces, all shown in \Cref{figure:domain_decomposition}.
The domain $\Omega$ on which the PDE is defined will be referred to as the \emph{physical domain}, with a boundary $\Gamma=\partial\Omega$. A thin-strip region $\A$, referred to as the \emph{annular domain}, is defined along the boundary $\Gamma$, within $\Omega$. The physical region $\fOmega=\Omega\setminus\A$ that is far from the boundary will be referred to as the \emph{faithful domain}. The curve $\I$ that separates $\aOmega$ and $\fOmega$ will be referred to as the \emph{interface}. Finally, we assume that $\Omega$ is contained within a geometrically simple domain $\mathcal{C}$, referred to as the \emph{computational domain}. Our method makes use of these spaces to perform a simple and straightforward domain decomposition strategy, which is illustrated in \Cref{figure:schematic}, and described in the steps below:
\begin{enumerate}
	\item The function $f$, known only inside $\Omega$ (shown in \Cref{figure:schematic}a), is smoothly truncated, so that it is $0$ at $\Gamma$, unchanged within the faithful domain $\fOmega$ and altered only within the thin boundary-adjacent strip $\aOmega$. We refer to this process as \emph{function intension}. The truncated (or \emph{intended}) function is shown in \Cref{figure:schematic}b, and the function in the annular region $\aOmega$ is shown in \Cref{figure:schematic}c.
	\item The \emph{regular problem} is solved in a geometrically simple region $\mathcal{C}\supset\Omega$ using standard methods, generating a solution $\ur$ (shown in \Cref{figure:schematic}d) that satisfies $\Delta\tilde u=f$ in $\fOmega$.
	\item The \emph{annular problem} is solved in $\aOmega$, generating a solution $\ua$ (shown in \Cref{figure:schematic}e) that satisfies $\Delta\ua=f$ in $\aOmega$.
	\item A \emph{stitching problem} is solved to correct any mismatch at the interface $\I$ between $\ur$ and $\ua$, generating a single solution $\ui$ that satisfies $\Delta\ui=f$ everywhere in $\Omega$. The ``unstitched'' function given by $\ur$ for $\x\in\fOmega$ and $\ua$ for $\x\in\aOmega$ is shown in \Cref{figure:schematic}f, and the ``stitched'' function $\ui$, where the interface mismatch has been corrected, is shown in \Cref{figure:schematic}g.
	\item A homogeneous boundary-value problem is solved, giving a correction function $\uh$ which satisfies $\Delta\uh=0$ everywhere in $\Omega$, with boundary conditions chosen so that $u=\ui+\uh$ satisfies the physical boundary conditions $u = g$. The homogeneous correction $\uh$ is shown in \Cref{figure:schematic}h, and the solution to the full PDE $u=\ui+\uh$ is shown in \Cref{figure:schematic}i.
\end{enumerate}

\subsection{Function intension, and solving the regular problem}
\label{section:methods_overview:intension_regular}

The primary goal of this paper is to replicate the main benefit of function extension methods --- converting a problem on a complex domain to a problem on a simple domain --- without having to extend either the inhomogeneous forcing $f$ or the unknown solution $u$. Instead, the function is simply cutoff in the interior of the domain $\Omega$, a process we refer to as \emph{function intension}, which we lay out here. The basic definition is simple:
\begin{equation}
	f_I(\x) =
		\begin{cases}
      f(\x),          &\x\in\tilde\Omega,\\
			\eta(\x) f(\x), &\x\in\A,\\
			0,              &\x\in\Omega^C,
		\end{cases}
\end{equation}
with $\eta$ a smooth cutoff function that is $1$ at $\I$ and $0$ at $\Gamma$, although it is simpler to think of $f_I$ being defined simply as $f_I(\x)=\eta(\x)f(\x)$, with $\eta$ understood to be $1$ for $\x\in\tilde\Omega$ and $0$ outside of $\Omega$, so that it is irrelevant that $f$ is unknown in $\Omega^C$. The convergence of our scheme will depend asymptotically on the smoothness of $\eta$ and practically (in terms of constants) on the specific choice of the function, which we make precise in \Cref{section:methods_preliminaries:cutoff}. Because $\eta$ is smooth, $f_I=\eta f$ is smooth in $\mathcal{C}$, and so $\Delta$ may be inverted using any appropriate method on the geometrically simple domain $\mathcal{C}$. Let us thus define $\ur$ by the solution to:
\begin{equation}
	\Delta \ur = \eta f \qquad\textnormal{in }\mathcal{C},
	\label{equation:decomposition:interior}
\end{equation}
along with any appropriate far-field boundary conditions that are required. We refer to $\ur$ as the \emph{regular problem}. We note that the solution $\ur$ is not unique and depends on the far-field boundary conditions chosen for $\mathcal{C}$ and any other modifications that must be made (see \Cref{section:methods_preliminaries:peridoic_compatibility}).

How good of a guess is $\ur$ to the actual solution to \Cref{equation:general_pde}? Well, it is correct, up to numerical errors, within the faithful domain $\fOmega$, i.e. $\Delta\ur(\x)=f(\x)$ for all $\x\in\fOmega$. However, $\ur$ fails to satisfy both the inhomogeneity in $\aOmega$ and fails to satisfy the boundary conditions.

\subsection{The annular problem}
\label{section:methods_overview:annular_problem}

The region of $\Omega$ in which $\tilde u$ fails to satisfy $\Delta\tilde u=f$ is, by construction, the annular region $\mathcal{A}$. We define a second problem, referred to as the \emph{annular problem}, stated here:
\begin{subequations}
	\begin{align}
		\Delta \ua &= f	&&\text{in }\aOmega,	\\
		\ua &= 0				&&\text{on }\Gamma\text{ and on }\I.
	\end{align}
	\label{equation:decomposition:annular}
\end{subequations}
Solving this problem is nontrivial, but it is also tractable as $\A$ is reasonably geometrically simple. To do so we will make use of a body-fitted coordinate system defined only within $\A$, but defer discussion of the details to \Cref{section:methods_specifics:annular_problem}. As with $\ur$, the solution $\ua$ is not unique, as any boundary conditions which make the PDE well-posed can be chosen. It is clear that $\ua$ satisfies the inhomogeneity in $\aOmega$ --- precisely where the regular solution $\ur$ fails to. It is tempting to define a solution candidate $u$ by $\ur$ for $\x\in\fOmega$ and $\ua$ for $\x\in\aOmega$. There are two issues with this: $u$ will still fail to satisfy the boundary condition given in \Cref{equation:general_pde}, but more importantly, $u$ may have jumps in both its value and its normal derivative at the interface $\I$. We consider this more pressing problem first.

\subsection{The stitching problem}
\label{section:methods_overview:stiching}

At this point, we assume that we have access to the regular solution $\ur$ for all $\x\in\fOmega$, and the annular solution $\ua$ for all  $\x\in\aOmega$, and so are free to interpolate and differentiate. We may thus evaluate the jump in both the solution and its normal derivative at the interface:
\begin{subequations}
	\begin{align}
		\gamma &= \lim_{\x\to{\I}^{\fOmega}}\ur - \lim_{\x\to\I^\aOmega}\ua,	\\
		\sigma &= \lim_{\x\to{\I}^{\fOmega}}\partial_\n\ur - \lim_{\x\to\I^\aOmega}\partial_\n\ua,
	\end{align}
\end{subequations}
with $\partial_\n f$ denoting the normal derivative of $f$. For a given elliptic PDE with known jump conditions \cite{LIE,HW}, these discontinuities in the value and the normal derivative of the function may be corrected by adding appropriate layer potentials. This leads to the \emph{inhomogeneous solution}, which for the Poisson equation, takes the form:
\begin{equation}
	\ui(\x) = 
	\begin{cases}
		\ur(\x) + (\mathcal{S}_\I\sigma)(\x) - (\mathcal{D}_\I\gamma)(\x) &\quad\text{for }\x\in\fOmega, \\
		\ua(\x) + (\mathcal{S}_\I\sigma)(\x) - (\mathcal{D}_\I\gamma)(\x) &\quad\text{for }\x\in\aOmega,
	\end{cases}
\end{equation}
where $\mathcal{S}$ and $\mathcal{D}$ denote the single and double layer potential operators associated with the Laplace operator \cite{HW}. The inhomogeneous solution $\ui$ now both satisfies the PDE $\Delta u=f$ everywhere in $\Omega$, is continuous, and has a continuous first-derivative, and thus is as smooth as supported by $f$ --- in particular, if $f\in C^k(\Omega)$, then $u\in C^{k+2}(\Omega)$ \cite{evans2010partial}. It is still the case, unfortunately, that $\ui$ fails to satisfy the boundary condition on the original PDE.

\subsection{The homogeneous problem}
\label{section:methods_overview:homogeneous}

Finally, we measure how much $\ui$ fails to satisfy the boundary condition by evaluating the discrepancy $\delta_\Gamma=g-\ui|_\Gamma$. It then remains only to solve the now \emph{homogeneous} equation:
\begin{subequations}
	\begin{align}
		\Delta\uh &= 0, &&\qquad\text{in }\Omega,	\\
		\uh &= \delta_\Gamma, &&\qquad\text{on }\Gamma.
	\end{align}
\end{subequations}
Well-conditioned methods based on boundary-integral equations allow for the solution and fast evaluation for a wide class of common PDE. Finally, having solved for $\uh$, we may define the solution $u=\ui + \uh$ to \Cref{equation:general_pde}, valid for all $\x\in\Omega$.

\section{Methods --- preliminaries}
\label{section:methods_preliminaries}

Having presented a sketch of the solution process, we now turn to the details of our specific implementation. In this initial manuscript, we focus on a simple global implementation, which nevertheless provides spectral accuracy, along with both setup and solution times with reasonable asymptotic scalings and constants. We begin first by explicitly defining our discretization of the domain $\Omega$, and its decomposition into $\fOmega$ and $\A$.

\subsection{Domain decomposition and discretization}
\label{section:methods_preliminaries:domain}

We assume that the boundary curve $\Gamma$ is given to us as a closed parametrized curve $\X(s)=(X(s), Y(s))$, for $s\in[0,2\pi)$, with $\X(0)=\X(2\pi)$, with a counter-clockwise orientation. The speed $\phi$ of this parametrization is given by $\phi(s)=\sqrt{\X_s\cdot\X_s}$. The outward pointing unit normal vector $\n$ is defined by $\phi(s)\n(s)=Y_s(s)\hat x - X_s(s)\hat y$. Given $\Gamma$, we may now define the annular domain $\A$.

\subsubsection{Definition of the annular domain $\A$}
\label{section:methods_preliminaries:annular_domain}

Near to the curve $\Gamma$, we can define a coordinate system $(s, r)$ by $\x=\X(s)+r\n(s)$, and the annular domain $\A$ is defined to be all points $\x\in\Omega$ with $|r|<R$ for some annular radius $R$. In the example given in \Cref{figure:domain_decomposition}, the physical domain $\Omega$ is interior to $\Gamma$; for such domains the annular region is defined for $-R\leq r\leq0$; for domains exterior to $\Gamma$ the annular region is defined for $0\leq r\leq R$.

It is clear that $R$ could be chosen so large that the coordinates $(s, r)$ are no longer well defined for all $\x\in\A$. The Jacobian of the coordinate map is given by $\mathcal{J}(\mathbf{s})=|\partial\x/\partial\mathbf{s}|=\phi(s)(1+r\kappa(s))$, where the coordinate $\mathbf{s}=(s,r)$ and $\kappa$ is the signed curvature given by $\kappa=-\n\cdot\X_{ss}/\phi^3$. Let $\kappa=\kappa^+ + \kappa^-$ give the decomposition into the positive ($\kappa^+$) and negative ($\kappa^-$) parts of $\kappa$. For an interior problem, $r<0$ for $\x\in\A$, and so $R_\text{max}=1/\max{\kappa_+}$ is sufficient to ensure $\mathcal{J}>0$. For an exterior problem, $r>0$, and so $R_\text{max}=-1/\min\kappa_-$. More simply, the maximum width of our annular region should be the minimum radius of curvature of $\Gamma$, being careful with signs. Because this criteria does not depend on $\phi$, the choice of $R$ is purely geometric, and is independent of the parametrization of $\Gamma$.

\begin{remark}[Coordinate splay]
  For our method to work, coordinate collapse must be prevented by choosing $R<R_\textnormal{max}$. However, coordinate splay can also occur when $R$ is large. Although this doesn't lead to the breakdown of the method, large coordinate splay leads to the amplification of effective discrete units, so that the discrete annulus provides an inefficient discretization of space.
\end{remark}

Once $R$ is chosen, the interface $\I$ is defined by the parametrized curve $\X(s) \pm R\n(s)$ for $s\in[0,2\pi)$, with positive sign for exterior problems and negative for interior problems. So long as $R<R_\textnormal{max}$, $\I$ has the same regularity as $\Gamma$.

\subsubsection{Discretization of the boundary and annular domain}
\label{section:methods_preliminaries:annular_domain_discretization}

We take an $N$ point discretization to $\Gamma$, with the discrete nodes of the boundary given by $\X_j=\X(s_j)$, where $s_j=j\Delta s$, with $\Delta s=2\pi/N$. Because $\X$ is periodic, boundary quantities such as $\n$ and $\phi$ can be computed using FFTs. For an interior problem, the $(s, r)$ rectangle $[0, 2\pi]\times[-R, 0]$ is discretized using an $N\times M$ tensor Fourier/Chebyshev mesh, with discrete values of $r$ at the first-kind Chebyhsev points $r_k=-R(\cos(\pi(2k+1)/(2M))+1)/2$, for $k=0,1,\ldots,M-1$, and discrete values of $s$ at the same nodes $s_j$ as used to discretize the boundary. The physical nodes for the discretized annular domain are then $\x_{jk}=\X_j + r_k\n_j$. Modification to the exterior case is straightforward. We define the smallest and largest discrete grid-spacing associated with the discretization of $\Gamma$ to be $h_\textnormal{min}=\min_j\phi_j\Delta s$ and $h_\textnormal{max}=\max_j\phi_j\Delta s$, respectively.

\subsubsection{Definition and discretization of the computational domain $\mathcal{C}$}
\label{section:methods_preliminaries:computational_domain_discretization}
We choose $\mathcal{C}$ to be a rectangle, and discretize that rectangle using a simple tensor-product Fourier representation. In principle, and in contrast to methods utilizing function extension, the computational domain $\mathcal{C}$ can be chosen to be tight around the boundary $\Gamma$. It is, however, simpler to compute $\min_j(X_j)$ than $\inf_s(X(s))$, and so we define $X_\text{min}=\min_j X_j - h_\text{max}$, $X_\text{max}=\max_j X_j + h_\text{max}$, and $Y_\text{min}$ and $Y_\text{max}$, analogously. We then define $C=[X_\text{min}, X_\text{max} + w]\times[Y_\text{min}, Y_\text{max} + w]$, where $w$ is a specified amount of wiggle-room. For many PDEs, $w$ can be taken to be $0$, but for Poisson and Stokes type-problems, satisfaction of compatibility conditions for the regular solution may require, in the worst case, $w=2Mh_\text{max}$, see \Cref{section:methods_preliminaries:peridoic_compatibility}.

Finally, we choose $N_x$ and $N_y$ discrete modes, and adjust the definition of $\mathcal{C}$, to ensure that we have an isotropic discretization with an even number of modes in each direction. To be precise, we choose $N_x=2\lceil(X_\text{max}+w-X_\text{min})/(2h)\rceil$, where $\lceil\cdot\rceil$ denotes the ceiling function and $h$ gives the target resolution of the regular grid (see \Cref{section:parameter_selection}). Defining $N_y$ analogously, we finally modify $\mathcal{C}$ to be $\mathcal{C}=[X_\text{min}, X_\text{min}+hN_x]\times [Y_\text{min}, Y_\text{min}+hN_y]$.

\subsubsection{Definition of the cutoff function $\eta$}
\label{section:methods_preliminaries:cutoff}

\begin{figure}[h!]
  \centering
  \begin{subfigure}[c]{0.3\textwidth}
    \centering
    \includegraphics[width=\textwidth]{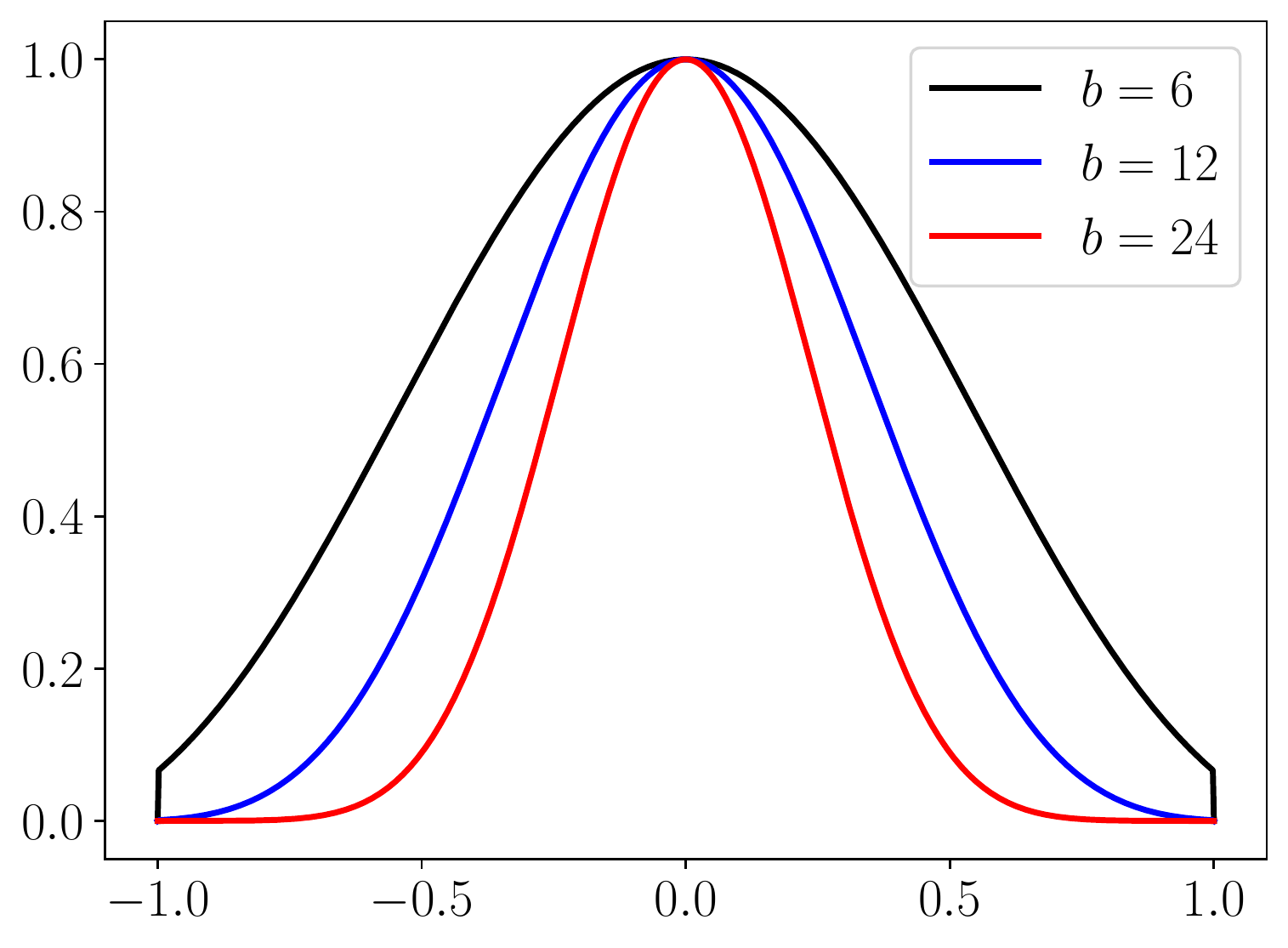}
  \end{subfigure}
  \hspace{1em}
  \begin{subfigure}[c]{0.31\textwidth}
    \centering
    \includegraphics[width=\textwidth]{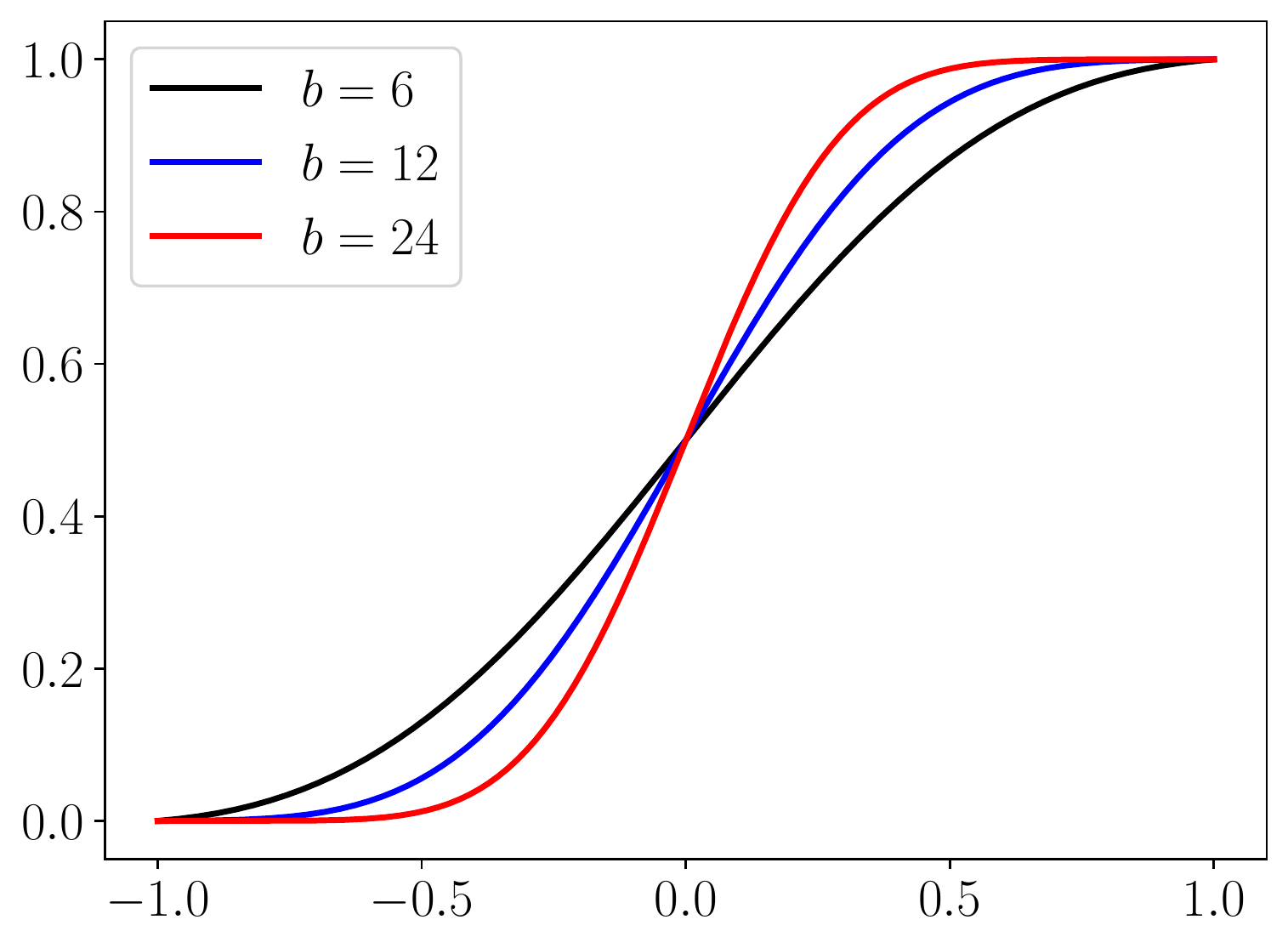}
  \end{subfigure}
  \caption{Regularized step and heaviside functions with bandwidths $b=6$ (black), $b=12$ (blue), and $b=24$ (red).}
  \label{figure:dpss}
\end{figure}

Let us assume for the moment that we have a discrete, scalar valued function $H$ of one variable that smoothly approximates a Heaviside function, with $H(x)=0$ for all $x<0$ and $H(x)=1$ for all $x>1$. For an interior problem, we can now define a function $\eta(\x)$ for any $\x\in\mathcal{C}$ by:
\begin{equation}
	\eta(\x) = 
	\begin{cases}
		0,\qquad&\text{for }\x\in\Omega^C,\\
		H(-r(\x)/R)&\text{for }\x\in\A,\\
		1,\qquad&\text{for }\x\in\fOmega,
	\end{cases}
\end{equation}
where $r(\x)$ gives the $r$-coordinate for the point $\x\in\A$ (with the obvious modification to be made for exterior problems). It is thus left to define $H$. While any smooth approximation of the Heaviside function should work, optimizing error requires that the function is both well resolved by the discretization and has continuous derivatives at $x=0$ and $x=1$. Throughout this work, we will use an integral of a prolate-spheroidal wavefunction, motivated by recent work on the non-uniform FFT \cite{barnett2019parallel}; we have found empirically that this improves performance slightly relative to simpler choices (such as a rescaled error function). To be precise, we use the standard definition of the Digital Prolate Spheroidal Sequence (DPSS), which is given by the dominant eigenvector of a matrix constructed from the sampled sinc function, with bandwidth $b/4$, as implemented in the scipy function \emph{scipy.window.dpss} \cite{2020SciPy-NMeth,barbosa1986maximum}. We sample this discrete bump function with sufficiently high frequency to ensure reconstruction to 15 digits using quintic spline interpolation; and construct $H$ as its antiderivative using adaptive quadrature (again, to 15 digits). This gives bump and step functions defined on $x\in[-1,1]$ (as shown in \Cref{figure:dpss}), which are transformed to the interval $[0,1]$ via an affine transformation.
\begin{remark}[Fast evaluation of DPSS Bump and Step Functions]
To allow fast usage, we have tabulated Chebyshev coefficients for both the bump and step functions for all $b$ such that $b=1,2,\ldots200$. The evaluation of these coefficient expansions can then be further accelerated by exploiting that (once properly centered) the bump functions are even and the step functions are odd \cite{press2007numerical}.
\end{remark}
Examples of these bump and step functions, for various values of $b$, are shown in \Cref{figure:dpss}. The choice of $b$ is considered in \Cref{section:parameter_selection:eta}, but in general higher values of $b$ will be used when $M$ is larger.  Naively, it may appear that the bump and step functions for small values of $b$, such as $b=4$ are poor choices --- they are not even continuous at $-1$ and $1$. However, these are typically used when $M$ is very small --- forcing a strong tradeoff between discrete resolution of the step function and its continuity and smoothness at the boundary.

\subsection{Common subproblems}
\label{section:methods_preliminaries:subproblems}

There are several steps in the algorithm which, while not an intimate part of the solution method, are required and not completely trivial. We discuss these here, to allow a more readable description of the specific implementation later in \Cref{section:methods_specifics}.

\subsubsection{The physical and coordinate problems}
\label{section:methods_preliminaries:coordinates}

\begin{figure}[h!]
  \centering
  \begin{subfigure}[c]{0.3\textwidth}
    \centering
    \includegraphics[width=\textwidth]{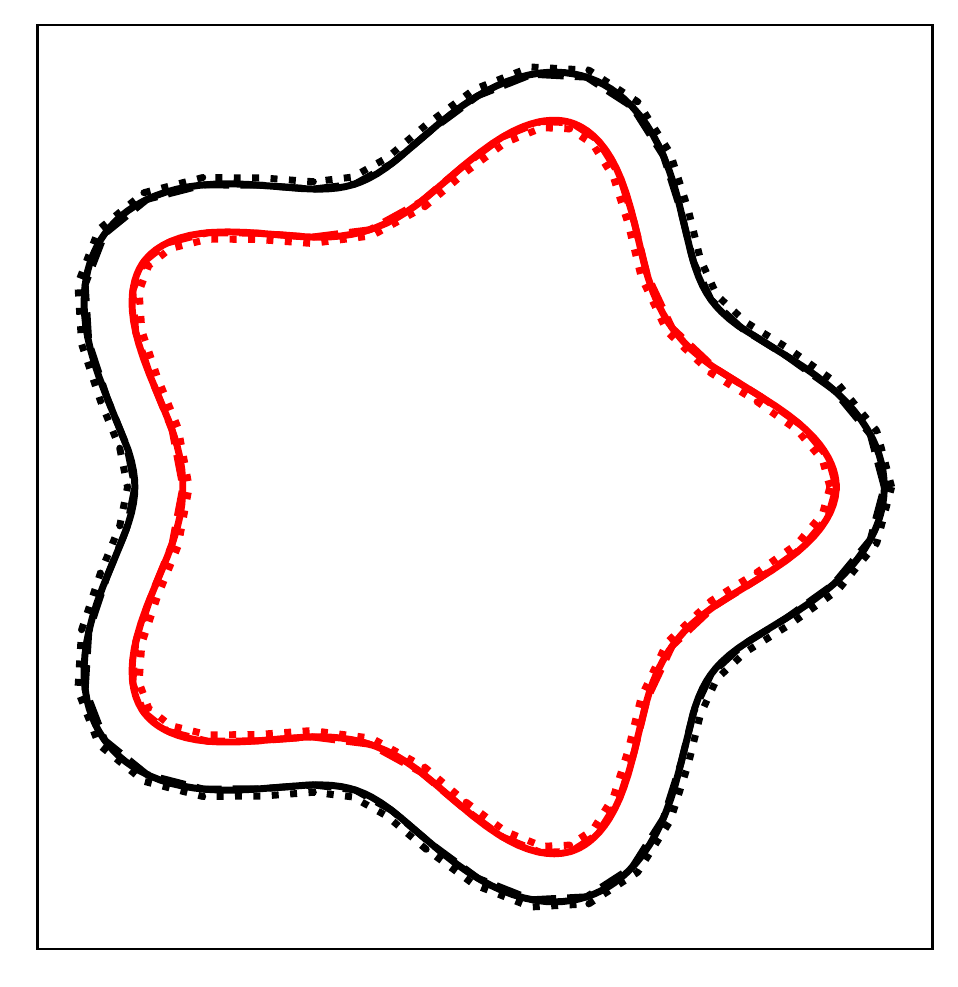}
  \end{subfigure}
  \hspace{1em}
  \begin{subfigure}[c]{0.31\textwidth}
    \centering
    \includegraphics[width=\textwidth]{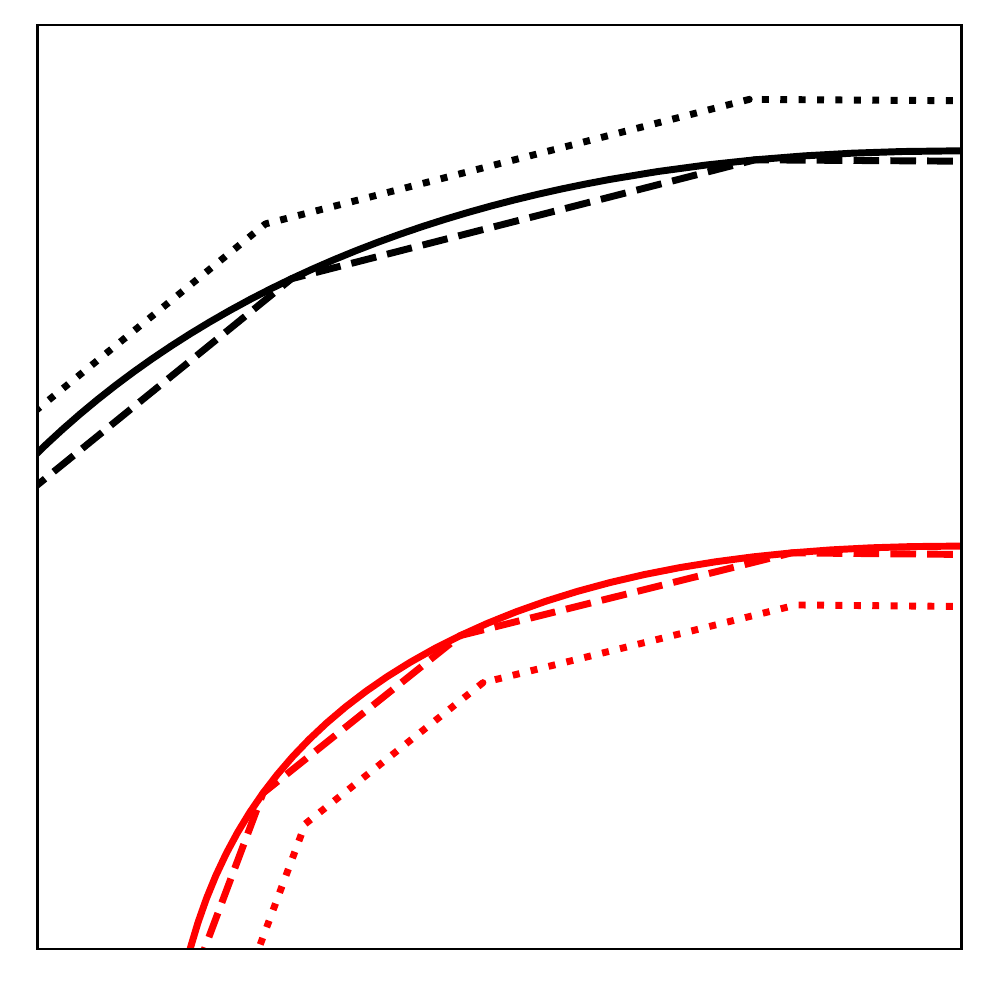}
  \end{subfigure}
  \caption{Discrete (dashed lines), continuous (solid lines), and modified discrete (dotted lines) curves for the boundary $\Gamma$ (in black) and $\I$ (in red).}
  \label{figure:coordinates}
\end{figure}

At this point, we have two separate discretizations: one for the simple computational rectangle $\mathcal{C}$ and one for the annular region $\A$. We must now connect these discretizations, which amounts to being able to accomplish the following task: given $\x\in\mathcal{C}$, determine if
\begin{enumerate}
	\item Is $\x$ interior or exterior to $\Omega$?
	\item If $\x$ is interior to $\Omega$, is it in $\A$ or $\fOmega$?
	\item If $\x\in\A$, for what values of $(r, s)$ do we have $\x=\X(s) + r\n(s)$?
\end{enumerate}
We would like to reduce questions (1) and (2) to questions about polygons, for which well-known algorithms \cite{hormann2001point} with robust implementations \cite{shapely2007} can be used. Unfortunately, it is not always the case that points inside the discrete polygon formed by connecting the points $\Gamma_j$ lie inside $\Gamma$. Instead, we define a modified discrete boundary curve $\tilde\Gamma$ by $\tilde\X = \X + 2\delta\n$, with $\delta=\max_j\delta_j$ and $\delta_j = |\kappa_j^{-1}| - \sqrt{\kappa_j^{-2}-(\Delta s\phi_j/2)^2}$. The distance $\delta_j$ comes from locally approximating the curve by a circle with radius $\kappa_j^{-1}$, with a safety buffer of 2 used. A modified discrete interface $\tilde\I$ is defined analogously. In \Cref{figure:coordinates}, we show an example with discrete and continuous curves $\Gamma$ and $\I$, along with the modified discrete curves $\tilde\Gamma$ and $\tilde\I$.
This now provides a simple algorithm for classifying $\x$, which we give for interior problems (with clear modifications for exterior problems):

\begin{algorithm}[H]
  \algrenewcommand\algorithmiccomment[2][\footnotesize]{{#1\hfill\(\triangleright\) #2}}  
  \caption{Classification of $\x$\hspace*{\fill}$\mathcal{O}(\log N)$}
  \algorithmicrequire{
  	Discrete polygons $\I$, $\tilde I$, $\Gamma$, $\tilde\Gamma$, point $\x\in\mathcal{C}$.
  }
  \begin{algorithmic}[1]
  	\If{$\x$ is exterior to $\tilde\Gamma$}
  		\Return{$\x\in\Omega^C$}
  	\ElsIf{$\x$ is interior to $\tilde\I$}
  		\Return{$\x\in\fOmega$}
  	\Else
  		\State Use Newton's method to solve for the coordinates $(s, r)$ which give $\x=\X(s)+r\n(s)$\footnotemark.
  		\If{$r>0$}
  			\Return{$\x\in\Omega^C$}
  		\ElsIf{$-R<r<0$}
  			\Return{$\x\in\A$}
  		\Else\ \Return{$\x\in\fOmega$}
  		\EndIf
  	\EndIf
  \end{algorithmic}
  \label{algorithm:coords}
\end{algorithm}
\footnotetext{For a unique solution to this Newton problem to be guaranteed, we require $R+2\delta<R_\textnormal{max}$. While we do check this condition, it fails only when both $R$ is very close to $R_\textnormal{max}$ and the boundary discretization is extremely coarse.}
Since $N^2$ points must be categorized, each step must be carefully optimized. A simple quadtree based approach, where whole blocks can be categorized as either in the exterior $\Omega^C$ or faithful domain $\fOmega$ is sufficient to reduce the lookups (steps 1 and 2) to $\mathcal{O}(\log N)$ for each point $\x$. Similarly, for a given point $\x$ whose coordinates must be determined, the Newton solver is $\mathcal{O}(1)$, so long as the routines used to interpolate $\X$ and $\n$ are. This can be accomplished, to spectral accuracy, using a type-2 nonuniform FFT, see \Cref{remark:accelerating_nufft}. Finally, explicitly solving the 2D Newton problem $\x=\X(s)+r\n(s)$ is more expensive than solving the scalar minimization problem $\inf_s|\X(s)-\x|$; once $s$ is known $r$ is easily recovered.

\subsubsection{Interpolation between domains, boundaries, and interfaces}
\label{section:methods_preliminaries:interpolation}

In multiple places in our algorithm, we will be required to interpolate between domains, or between domains and boundaries/interfaces. We discuss a few cases here: (1) Interpolation of a smooth function known everywhere in $\mathcal{C}$ to anywhere else in $\mathcal{C}$; (2) interpolation of a function known on the discrete nodes of the annular domain to arbitrary $\x\in\A$; and (3) interpolation of a function known at the discrete nodes of the annular domain to the discrete boundary and interface nodes.

Interpolation from $\mathcal{C}$ to general $\x$ (1) is the simplest. Given a function $f$ on the discrete nodes of $\mathcal{C}$, we compute its two-dimensional FFT $\hat f$. We may then, after suitably rescaling $\x$ to the unit-periodic rectangle $[0,2\pi]^2$, compute $f(\x)$ by evaluating the type-2 nonuniform FFT (NUFFT). The error in this interpolation converges spectrally fast in the number of modes $N_x$ and $N_y$ used to discretize $\mathcal{C}$, so long as $f$ is smooth and periodic.
\begin{remark}[Accelerating interpolation via NUFFT]
	In our algorithm, we will typically have to compute the NUFFT to only $\mathcal{O}(MN)$ targets, with $M\ll N$. Since there are a small number of targets, the interpolation cost is small relative to the oversampled FFT, and we can exploit the option provided in FINUFFT \cite{barnett2019parallel} to use an oversampling factor of 1.25 rather than 2. This reduces the cost from approximately 4 times a grid-sized FFT time to $\approx1.5$ times. An additional optimization comes from realizing that only the final interpolation step of the NUFFT algorithm depends on $\x$; we may thus precompute all other portions of the algorithm and interpolate to multiple different collections of target values $\x$ at a fixed cost per target $\x$.
	\label{remark:accelerating_nufft}
\end{remark}
Interpolating a function $f$ from the discrete nodes of $\A$ to any $\x\in\A$ (2) is only slightly harder. Given an arbitrary $\x\in\A$, we first compute its coordinates $(s, r)$ (see \Cref{section:methods_preliminaries:coordinates}). The function $f$ is periodic in $s$ but not in $r$; but it is defined on a Chebyshev mesh in $r$. We may thus take its even reflection and modify it to obtain a periodic function \cite{boyd2001chebyshev}, and then, upon suitably rescaling $r$ to lie within the unit periodic interval $[0,2\pi]$, again use a type-2 NUFFT to compute $f(\x)$. Errors in this interpolation converge spectrally fast in the number of boundary modes $N$ and the number of Chebyshev modes $M$.

Finally, we consider case (3), interpolating a function $f$ from the discrete nodes of $\A$ to either the discrete boundary nodes $\Gamma_j$ or the discrete interface nodes $\I_j$. Here we may utilize simple univariate Chebyshev interpolation at fixed values of $s$, see \cite{boyd2001chebyshev}. Errors in this interpolation converge spectrally fast in the number of Chebyshev modes $M$.

\subsubsection{Periodic compatibility conditions}
\label{section:methods_preliminaries:peridoic_compatibility}

For both the Poisson and Stokes problems (and potentially other PDE), it will be necessary to slightly modify the simple algorithm presented in \Cref{section:methods_overview} due to the constant null-space in the periodic operators used to solve the regular problem. We discuss the Poisson case here; the adjustments required for the Stokes problem are no different. In particular, we have the solvability condition that $\Delta u=\eta f$ is solvable on a periodic domain only if $\eta f$ has mean $0$; that is, if $\int_\mathcal{C} \eta f=0$. This will typically not be the case. There are several solutions to this problem; here, we use perhaps the simplest --- modifying the function $\eta f$ \emph{only in the exterior of $\Omega$}. This is accomplished by computing the integral of $\eta f$ over $\Omega$, and adding a regularized \emph{bump} function in the domain with a countervailing mass. To be precise, we let $\xi(\x)$ be a bump function with $\int_\mathcal{C}\xi=1$ and $\text{supp}\,\xi$ disjoint from $\Omega$. We may then solve:
\begin{equation}
	\Delta u=\eta f - \xi\int_\mathcal{C}\eta f,
\end{equation}
where the right-hand side clearly satisfies the solvability condition and is unchanged within $\fOmega$ (actually, $\Omega$). The numerical approximation may be affected, if, however, $\xi$ is poorly resolved by the discretization. One solution is to use the bump function from which $\eta$ is constructed, see \Cref{section:methods_preliminaries:cutoff}. However, because $\eta$ is constructed as an integral of the bump function, it is better resolved than the underlying bump. Instead, we use that bump function but spread out over twice the radius, and have found empirically that this does not affect the overall error in the method.

A final question is how to place the center of the bump function $\xi$. In many cases, this can be placed by the user; for example, in the problem shown in \Cref{figure:schematic}, the bump could be manually centered in the upper right hand corner, and the computational domain $\mathcal{C}$ could be taken tight to $\Omega$. A solution that always works, without manual intervention from the user, is to simply increase the size of the domain by $w=2Mh$ (see \Cref{section:methods_preliminaries:computational_domain_discretization}), in which case the center can be placed at $(X_\text{max}+Mh, Y_\text{max}+Mh)$ and guaranteed not to overlap $\Omega$. This simplified approach, which comes at the cost of slightly increasing the number of degrees of freedom in $\mathcal{C}$, is what we will use throughout this manuscript.

\subsubsection{Evaluation of layer potentials}
\label{section:methods_preliminaries:close_evaluation}

In two places --- both in the stitching problem and in the homogeneous correction, we will be required to evaluate layer potentials to a large number of points, many of which may be close to the curve from which these layer potentials emanate. This presents two problems:
\begin{enumerate}
	\item Naive evaluation near boundaries typically have $\mathcal{O}(1)$ errors that do not disappear on grid refinement \cite{barnett2015spectrally,helsing2008evaluation};
	\item Direct evaluation is $O(N_b N_t)$, with $N_b$ and $N_t$ giving the number of boundary points and target points, respectively. Here, $N_b=N$ and $N_t\propto N^2$, and so direct evaluation is $O(N^3)$.
\end{enumerate}
 For all problems in this paper, we will use a recently developed method called \emph{Quadrature by Fundamental Solutions} \cite{stein2021quadrature}, combined with a standard Fast-Multipole method library \cite{HFMM2D}, which provides a reasonable solution to both problems. This method, once setup is done, allows for highly-accurate evaluation of layer-potentials up to and on the boundary of the domain, in a kernel independent manner, in $\mathcal{O}(N_b^2+N_t)=\mathcal{O}(N^2)$ time. Unfortunately, it does come with a $\mathcal{O}(N_b^3)$ setup cost. When $N_b$ is relatively small --- up to $10,000$ or so, modern implementations of LAPACK are surprisingly fast and this is rarely the dominant cost. For larger problems, other, more scalable methods, such as panel discretizations with kernel specific close-evaluation quadratures \cite{helsing2008evaluation}, would be more efficient and maintain an asymptotic scaling of $\mathcal{O}(N^2)$, in both setup and solution stages.

\section{Methods --- specifics}
\label{section:methods_specifics}

We are now ready to return to solving a PDE, and provide specific algorithms for all stages of the computation. This section goes deeper into the method as developed in \Cref{section:methods_overview}, utilizing the tools developed in \Cref{section:methods_preliminaries}. We restate our model Poisson problem here:
\begin{subequations}
	\begin{align}
		\Delta u &= f &&\textnormal{in }\Omega,	\\
		u				 &= g &&\textnormal{on }\Gamma.
	\end{align}
\end{subequations}
For simplicity, we assume that $f$ and $g$ are given to the user as evaluatable functions on their respective domains.

Our first step is to discretize the problem: in this section we assume that the number of boundary nodes $N$, the number of Chebyshev modes $M$, and a grid-spacing $h$ are given. The boundary and annulus are discretized as described in \Cref{section:methods_preliminaries:annular_domain_discretization}, with $R=Mh$. If $R>R_\textnormal{max}$, no solution is attempted. A more principled way to set these parameters is given in \Cref{section:parameter_selection}, but we proceed in this way for now to demonstrate certain basic features that would be otherwise hidden. All discrete nodes in $\mathcal{C}$ are categorized into physical points ($\in\Omega$) and exterior points ($\in\Omega^C)$, and those that are physical points are further categorized into points within the annular region ($\in\A$) and those within the faithful region ($\in\fOmega$), via the technique described in \Cref{section:methods_preliminaries:coordinates}. The cutoff function $\eta$ is then evaluated for all discrete nodes of $\mathcal{C}$. Since the Poisson problem has a nullspace, a regularized bump function $\xi$ is also computed at all discrete nodes of $\mathcal{C}$ (see \Cref{section:methods_preliminaries:peridoic_compatibility}). At this point, our domain is discretized and we have all relevant information required to solve PDEs on the domain. We evaluate $f$ at the discrete nodes of $\mathcal{C}$ and $\A$, and $g$ at the discrete nodes $\X_j$.

\subsection{Function intension, and solving the regular problem}
\label{section:methods_specifics:regular_problem}

We can now describe explicitly how to solve the \emph{regular problem} outlined in \Cref{section:methods_overview:intension_regular}. First, the intended function $\eta f$ is computed at all discrete nodes of $\mathcal{C}$ by defining it to be $0$ for $\x\in\Omega^C$ and to be $\eta f$ for $\x\in\Omega$. A modified $f$ with mean $0$ is computed as $f^m=\eta f-\xi\int_\mathcal{C}\eta f$, and its Fourier transform $\widehat{f^m}_\k$ is computed via the standard two-dimensional FFT. The Fourier modes of $\widehat{\ur}_\k$ are then computed as $\widehat{\ur}_\k=\widehat{f^m}_\k/|\k|^2$, with $\widehat{\ur}_{\boldsymbol{0}}=0$, with $\k$ the wavevectors for the domain $\mathcal{C}$. The regular solution $\ur$ is then recovered by the standard two-dimensional inverse FFT. Since $N_x$ and $N_y$ are proportional to $N$, the total cost for solving the regular problem scales as $\mathcal{O}(N^2\log N)$.

\subsection{The annular problem}
\label{section:methods_specifics:annular_problem}

We  now seek to solve the annular problem:
\begin{subequations}
	\begin{align}
		\Delta \ua &= f	&&\text{in }\aOmega,	\\
		\ua &= 0				&&\text{on }\Gamma\text{ and on }\I.
	\end{align}
	\label{equation:decomposition:annular}
\end{subequations}
In the coordinates $s$ and $r$ for the annular region $\A$, the Laplace operator is given by:
\begin{equation}
    \Delta\ua = \frac{1}{\psi}\left[\frac{\partial}{\partial r}\left(\psi\frac{\partial\ua}{\partial r}\right) + \frac{\partial}{\partial s}\left(\frac{1}{\psi}\frac{\partial\ua}{\partial s}\right)\right],
    \label{equation:annular_laplacian}
\end{equation}
where $\psi(s, r)$ is given by $\psi(s, r) = \phi(s)(1 + r\kappa(s))$, with $\kappa$ the curvature of $\Gamma$. Inverting this operator is not completely trivial. For discretizations where $NM$ is relatively small, it is probably reasonable to form $\Delta$ and directly invert it using dense linear algebra. This, however, comes at a setup cost of $\mathcal{O}((NM)^3)$ (to factor the operator), and an application cost of $\mathcal{O}((NM)^2)$. Especially for larger $M$, this would significantly dominate the entire computation. Instead, we seek an iterative solution to the problem.

We first note that if $\Gamma$ is a circle with a uniform parametrization, then both $\kappa$ and $\phi$ are independent of $s$. This means that we can take the Fourier transform to obtain:
\begin{equation}
	\widehat{\Delta\ua}_k = \frac{1}{\psi}\left[\frac{\partial}{\partial r}\left(\psi\frac{\partial\widehat{\ua}_k}{\partial r}\right) - \frac{1}{\psi}k^2\widehat{\ua}_k\right],
	\label{equation:annular_laplacian_k}
\end{equation}
This operator can, upon discretizing $\partial/\partial r$ via standard Chebyshev operators, be formed and directly inverted, independently for each $k$, at a total cost of $\mathcal{O}(NM^3)$. We may thus solve the annular problem on a simple circle via the following algorithm:
\begin{enumerate}
	\item Given $f$, compute $\hat f_k$ via the 1D FFT. Cost: $\mathcal{O}(MN\log N)$.
	\item For each mode $k$, invert \Cref{equation:annular_laplacian_k} to obtain $\widehat{\ua}_k$. Total cost: $\mathcal{O}(NM^2)$.
	\item Compute $\ua$ via the 1D inverse FFT. Cost: $\mathcal{O}(MN\log N)$.
\end{enumerate}
While this procedure does not work for general domains, it provides a surprisingly effective preconditioner, and it is straightforward to apply the annular Laplacian given in \Cref{equation:annular_laplacian} for a total cost of $\mathcal{O}(MN\log N)+\mathcal{O}(NM^2)$. We have found, and will show in the examples, that GMRES with the circular preconditioner converges robustly, with a relative residual of $10^{-14}$ typically reached in 10-20 iterations, so long as $R$ does not approach too close to $R_\text{max}$. When $R$ approaches $R_\text{max}$, the iteration count can increase dramatically.

\begin{remark}[Spectral vs. Pseudospectral]
	We have implemented this solver discretizing both the nodal values and the spectral modes. For scalar problems, both methods work about equally well, although the spectral method sometimes saturates at slightly lower errors (and for this reason we use it throughout this manuscript). For vector problems (and in particular Stokes), we have found it simpler to obtain robust convergence by discretizing spectral modes, omitting from the solution vector the Nyquist frequency in the azimuthal direction, see \Cref{section:stokes}. In both implementations, we use a rectangular method for discretizing the Chebyshev operators and imposing the boundary conditions \cite{driscoll2016rectangular}.
\end{remark}

\begin{remark}[Choice of boundary conditions]
	For simplicity, we have chosen $u=0$ as the boundary conditions for the annular solve at both $\I$ and $\Gamma$. In certain circumstances, it could be beneficial to choose other boundary conditions, subject to the constraint that those conditions satisfy any compatibility conditions for the PDE. One such circumstance is if the user has a method to perform close-evaluation of either a single-layer or double-layer potential, but not both. Consider the case where the user only has a method to apply the single-layer potential. Rather than setting the boundary condition on $\I$ to be $0$, we could instead set it to be $\ur|_\mathcal{I}$. Looking briefly ahead to \Cref{eq:jumps}, we see that then $\gamma=0$, and so the double-layer potential in \Cref{eq:laplace_stitching} can be ignored. Alternatively, if the user only has a method to apply the double-layer potential, they could fix $\partial_\n\ua=\partial_\n\ur$ at $\I$, instead.
	\label{remark:annular_bcs}
\end{remark}

\subsection{The stitching problem}
\label{section:methods_specifics:stiching}

We now know $\ur(\x)$ for all $\x\in\mathcal{C}$, and the annular solution $\ua(\x)$ for all $\x\in\aOmega$. We seek now to compute the jumps in the solution and its normal derivative:
\begin{subequations}
	\begin{align}
		\gamma &= \lim_{\x\to{\I}^{\fOmega}}\ur - \lim_{\x\to\I^\aOmega}\ua,	\\
		\sigma &= \lim_{\x\to{\I}^{\fOmega}}\partial_\n\ur - \lim_{\x\to\I^\aOmega}\partial_\n\ua,
	\end{align}
	\label{eq:jumps}
\end{subequations}
with $\partial_\n f$ denoting the normal derivative of $f$ and $\x\to\I^W$ denoting the limit as $\x$ tends to $\I$ from within the domain $W$. Computing $\ur$ and $\ua$ on $\I$ can be done directly using the interpolation scheme given in \Cref{section:methods_preliminaries:interpolation}. The normal derivatives can be computed by computing $\partial_x\ur$ and $\partial_y\ur$ via FFT based differentiation, and $\partial_r\ua$ via Chebyshev differentiation, before again applying the interpolation operators given in \Cref{section:methods_preliminaries:interpolation}. This immediately gives $\partial_\n\ua$, and $\partial_\n\ur$ can be computed then as $\partial_\n\ur = \n\cdot\grad\ur$.

For the Poisson problem, we may correct these jumps by adding the single and double layer potentials $\mathcal{S}_\I\sigma-\mathcal{D}_\I\gamma$, given explicitly by:
\begin{equation}
	(\mathcal{S}_\I\sigma - \mathcal{D}_\I\gamma)(\x) = \int_\I G(\x,\y)\sigma(\y)\,d\mathbf{s}_\y -  \int_\I \frac{\partial G(\x,\y)}{\partial\n^\y}\gamma(\y)\,d\mathbf{s}_\y,
	\label{eq:laplace_stitching}
\end{equation}
where $G(\x,\y)=\-(2\pi)^{-1}\log|\x-\y|$. As discussed in \Cref{section:methods_preliminaries:close_evaluation}, we utilize FMM accelerated QFS to evaluate these layer potentials, but briefly describe this here in the context of the stitching step for an interior problem. Points in $\fOmega$ are inside $\I$, while points in $\A$ are outside $\I$. These require different treatment. To evaluate to all $\x\in\fOmega$, we use QFS to construct an effective inward\footnote{Note that we use the term ``inward'' here, rather than ``interior''. Inward will be used for potentials evaluated from $\I$ into $\fOmega$, while ``outward'' will be used for potentials evaluated from $\I$ into $\A$. For interior problems, inward corresponds to interior and outward cooresponds to exterior; for exterior problems, inward corresponds to exterior and outward corresponds to interior.} representation; that is, we compute a source curve $\I_\textnormal{in}$ and an effective potential $\zeta_\textnormal{in}$ such that:
\begin{equation}
	\left(\mathcal{S}_\I\sigma - \mathcal{D}_\I\gamma\right)|_\fOmega(\x) = \mathcal{S}_{\I_\textnormal{in}}\zeta_\textnormal{in}(\x),
	\label{eq:qfs_constraint}
\end{equation}
for all $\x\in\I$. Note that the layer potentials on the left-hand side are singular or principal-value, and so we specify here that the interior limit (from $\fOmega$) is taken. Once $\zeta_\textnormal{in}$ has been determined, it can be evaluted from $\I_\textnormal{in}$ using a standard periodic trapezoid rule and FMM acceleration to all points $\x\in\fOmega$:
\begin{equation}
	(\mathcal{S}_{\I_\textnormal{in}}\zeta_\textnormal{in})(\x) \approx \frac{-1}{2\pi}\sum_{j=1}^{N_s} \log|\x-\X^{\I_\textnormal{in}}(s_j)|\zeta_\textnormal{in}(s_j)w_j,
	\label{eq:discrete_qfs}
\end{equation}
where $N_s$ is the number of points discretizing the source curve $\I_\textnormal{in}$ and $w_j=2\pi\phi_j/N_s$, with $\phi_j$ the speed for the parametrization of $\I_\textnormal{in}$. The number of source points is typically $N_s=N$, although this may need to be adjusted due to geometric constraints, see \cite{stein2021quadrature}; and this method for evaluating the integral \Cref{eq:laplace_stitching} converges spectrally fast in $N$. The process of evaluating the layer potential in \Cref{eq:laplace_stitching} at $\x\in\A$ is similar; an effective outward potential $\zeta_\textnormal{out}$ on a source curve $\I_\textnormal{out}$ is computed, subject to the constraint that \Cref{eq:qfs_constraint} holds but now with the exterior limit (from $\A$) taken on the right hand side. The layer potential in \Cref{eq:discrete_qfs}, with interior and exterior quantities swapped, is evaluated via FMM at all $\x\in\A$. We then define the \emph{inhomogeneous solution}:
\begin{equation}
	\ui(\x) = 
	\begin{cases}
		\ur(\x) + (\mathcal{S}_\I\sigma-\mathcal{D}_I\gamma)(\x) &\quad\text{for }\x\in\fOmega, \\
		\ua(\x) + (\mathcal{S}_\I\sigma-\mathcal{D}_I\gamma)(\x) &\quad\text{for }\x\in\aOmega.
	\end{cases}
	\label{eq:particular}
\end{equation}
We now have a smooth particular solution $u_I$ to the PDE in $\Omega$. The algorithm for the stitching step is summarized in \Cref{algorithm:stitch}, with computational scaling given for our implementation; the total effort is $\mathcal{O}(N^2\log N)+\mathcal{O}(NM^2)$.

\begin{algorithm}[H]
  \algrenewcommand\algorithmiccomment[2][\footnotesize]{{#1\hfill\(\triangleright\) #2}}  
  \caption{Stitching problem \hspace*{\fill} $\mathcal{O}(N^2\log N + NM^2)$}
  \algorithmicrequire{
  	The Fourier transform of the regular solution $\hat\ur$, and the annular solution $\ua$.
  }
  \begin{algorithmic}[1]
  	\State [$\mathcal{O}(N^2)$] compute $\widehat{\partial_x\ur}$ and $\widehat{\partial_y\ur}$
  	\State [$\mathcal{O}(N^2\log N)$] interpolate $\ur$, $\partial_x\ur$, and $\partial_y\ur$ onto $\I$
  		\Comment{see \Cref{section:methods_preliminaries:interpolation}}
  	\State [$\mathcal{O}(N)$] compute $\partial_\n\ur=\n\cdot\grad\ur$
  	\State [$\mathcal{O}(N M^2)$] compute $\ua$ and $\partial_\n\ua$ on $\I$
  		\Comment{see \Cref{section:methods_preliminaries:interpolation}}
  	\State [$\mathcal{O}(N)$] compute $\gamma$ and $\sigma$
  		\Comment{see \Cref{eq:jumps}}
  	\State [$\mathcal{O}(N^2)$] compute inward effective potential $\zeta_\textnormal{in}$ defined on effective source curve $\I_\textnormal{in}$
  	\State [$\mathcal{O}(N^2)$] Evaluate $\mathcal{S}_{\I_\textnormal{in}}\zeta_\textnormal{in}$ at all points $\x\in\fOmega$ and add to $\ur$ to obtain $u_I$ for $\x\in\fOmega$
  	\State [$\mathcal{O}(N^2)$] compute outward effective potential $\zeta_\textnormal{out}$ defined on effective source curve $\I_\textnormal{out}$
  	\State [$\mathcal{O}(N^2)$] Evaluate $\mathcal{S}_{\I_\textnormal{out}}\zeta_\textnormal{out}$ at all points $\x\in\A$ and add to $\ua$ to obtain $u_I$ for $\x\in\A$
  	\Statex\Comment{Steps 6-9, see \Cref{section:methods_preliminaries:close_evaluation}}
  \end{algorithmic}
  \label{algorithm:stitch}
\end{algorithm}

\subsection{The homogeneous problem}
\label{section:methods_specifics:homogeneous}

Finally, we must correct $u_I$ to satisfy the physical boundary conditions. First, $u_I$ is interpolated to the boundary nodes $\X_j$ from the annular nodes using the method described in \Cref{section:methods_preliminaries:interpolation}. We now know the discrepancy $g-u_I$ at all discrete boundary nodes. It thus remains to solve the homogeneous equation:
\begin{subequations}
	\begin{align}
		\Delta u_H &= 0,	&&\text{in }\Omega,\\
		u_H &= g-u_I &&\text{on }\Gamma.
	\end{align}
\end{subequations}
For this interior Dirichlet problem, the simple representation $u_H=\mathcal{D}_\Gamma\zeta$ suffices to provide the well-conditioned second-kind boundary integral equation \cite{HW}:
\begin{equation}
	\left(\mathbb{I}/2 - D_{\Gamma,\Gamma}\right)\zeta = g - u_I,
\end{equation}
with $D_{\Gamma,\Gamma}$ denoting the principal value operator obtained by evaluating $\mathcal{D}_\Gamma$ on $\Gamma$. The left-hand side operator can be discretized with spectral accuracy in $N$ utilizing singular Kress quadrature \cite{hao2014high}, to give the $N\times N$ matrix equation:
\begin{equation}
	A^{ij}\zeta_j = (g-u_I)|_{s_j},
	\label{eq:discretized_homogeneous}
\end{equation}
which requires $\mathcal{O}(N^2)$ effort to form. This matrix can then be directly factored with $\mathcal{O}(N^3)$ effort, to allow its solution in $\mathcal{O}(N^2)$ time, or solved via GMRES in $\mathcal{O}(N^2)$ time via direct application, with an iteration count independent of $N$.
\begin{remark}[FMM and Fast direct solvers]
	In cases where $N$ is large, further accelerations are possible. If the physical problem is well conditioned, using GMRES and applying $A$ with an FMM and local corrections reduces the cost to $\mathcal{O}(N)$ []. When the physical problem itself is poorly conditioned (such as high-frequency Helmholtz problems), the number of GMRES iterations may become large. Using a fast-direct solver is one option in this case \cite{martinsson2005fast}, which we have not yet explored.
\end{remark}

Finally, once $\zeta$ is known, $u_H=\mathcal{D}_\Gamma\zeta$ can be evaluated for all $\x$ in both $\fOmega$ and $\A$, again using the method described in \Cref{section:methods_preliminaries:close_evaluation}. The total numerical cost of the homogeneous correction is $\mathcal{O}(N^2)$.

\subsection{Finishing up}
\label{section:methods_specifics:finishing}

At this point, we know $u_I$ and $u_H$ for all discrete nodes $\x$ of $\mathcal{C}$ in $\fOmega$ and discrete nodes $\x$ of $\A$. We can simply add these together to obtain $u$. We do not, however, at this point know $u$ on the discrete nodes $\x$ of $\mathcal{C}$ within the region $\A$. These values can now be obtained by interpolating $u$ from $\A$ to any $\x\in\A$ using the method described in \Cref{section:methods_preliminaries:interpolation}. This gives us a full representation of $u$; we can now integrate, differentiate, or interpolate $u$ to any location in $\Omega$ with spectral accuracy (in $N$ and $M$).

\subsection{Parameters}
\label{section:methods_specifics:parameters}

Finally, we collect together the various parameters that need to be set in order to fully define the method. These are given in \Cref{table:parameter_choices}. Note that the number of modes discretizing $\A$ is taken to be $N$, and the tolerance used for \emph{all} iterative and approximate methods is set to the same value of $\epsilon=10^{-14}$ throughout. Because this tolerance is used in multiple approximate methods throughout the solver, we do not expect solutions to achieve this tolerance, but rather a small multiple of it.

\begin{table}
	\centering
	\begin{tabular}{c|l}
	\hline
	\hline
	Parameter      & Description	\\
	\hline
	$N$				     & Number of discrete boundary modes                          \\
	$M$            & Number of Chebyshev modes discretizing $\A$	              \\
	$R$				     & Annular radius                                             \\
	$h$            & gridspacing of regular discretization of $\mathcal{C}$     \\
	$b$            & bandwidth of $\eta$                                        \\
	$\epsilon$     & tolerance used for all iterative and approximate methods   \\
	\hline
	\hline
	\end{tabular}
	\caption{Required parameter choices.}
	\label{table:parameter_choices}
\end{table}

\section{A simple example and dependence on $M$}
\label{section:simple_poisson}

We now return to the example problem shown in \Cref{figure:schematic}, and analyze the convergence properties of the scheme, deferring a discussion of parameter selection and numerical performance to \Cref{section:parameter_selection}, as immediately jumping to optimal parameter choices obscures some details. The problem we will analyze is the Dirichlet Poisson problem, with $f$ and $g$ manufactured from the known solution $u=e^{\sin(x)}\sin(2y) + \log(0.1 + \cos(y)^2)$, set on a star-shaped domain defined by the function
\begin{equation}
	\X(s) = (x_c + r\omega(s)\cos(s),\ y_c + r\phi(s)\sin(s)),
\end{equation}
with $\omega(s)=1+a\cos(ds)$. For this problem we take $x_c=y_c=0$, $r=1$, $d=5$, and $a=0.15$. This is a relatively simple domain to allow more exploratory range in some of the numerical experiments done in this section; more complex domains are considered in later examples.

\subsection{Behavior for fixed $M$}
\label{section:simple_poisson:fixed_M}

\begin{figure}[h!]
  \centering
  \begin{subfigure}[c]{0.4\textwidth}
    \centering
    \includegraphics[width=\textwidth]{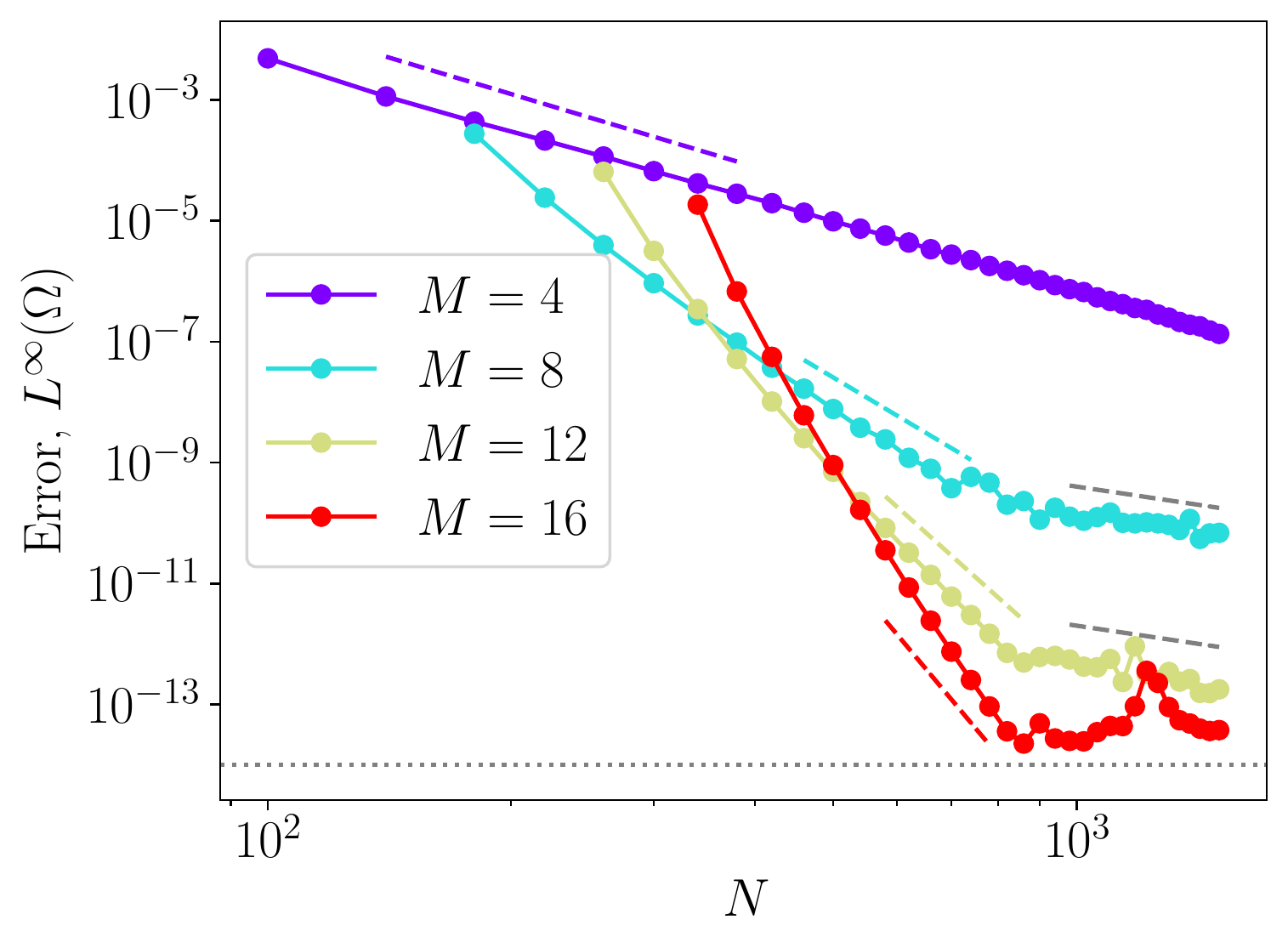}
    \caption{}
  \end{subfigure}
  \hspace{1em}
  \begin{subfigure}[c]{0.4\textwidth}
    \centering
    \includegraphics[width=\textwidth]{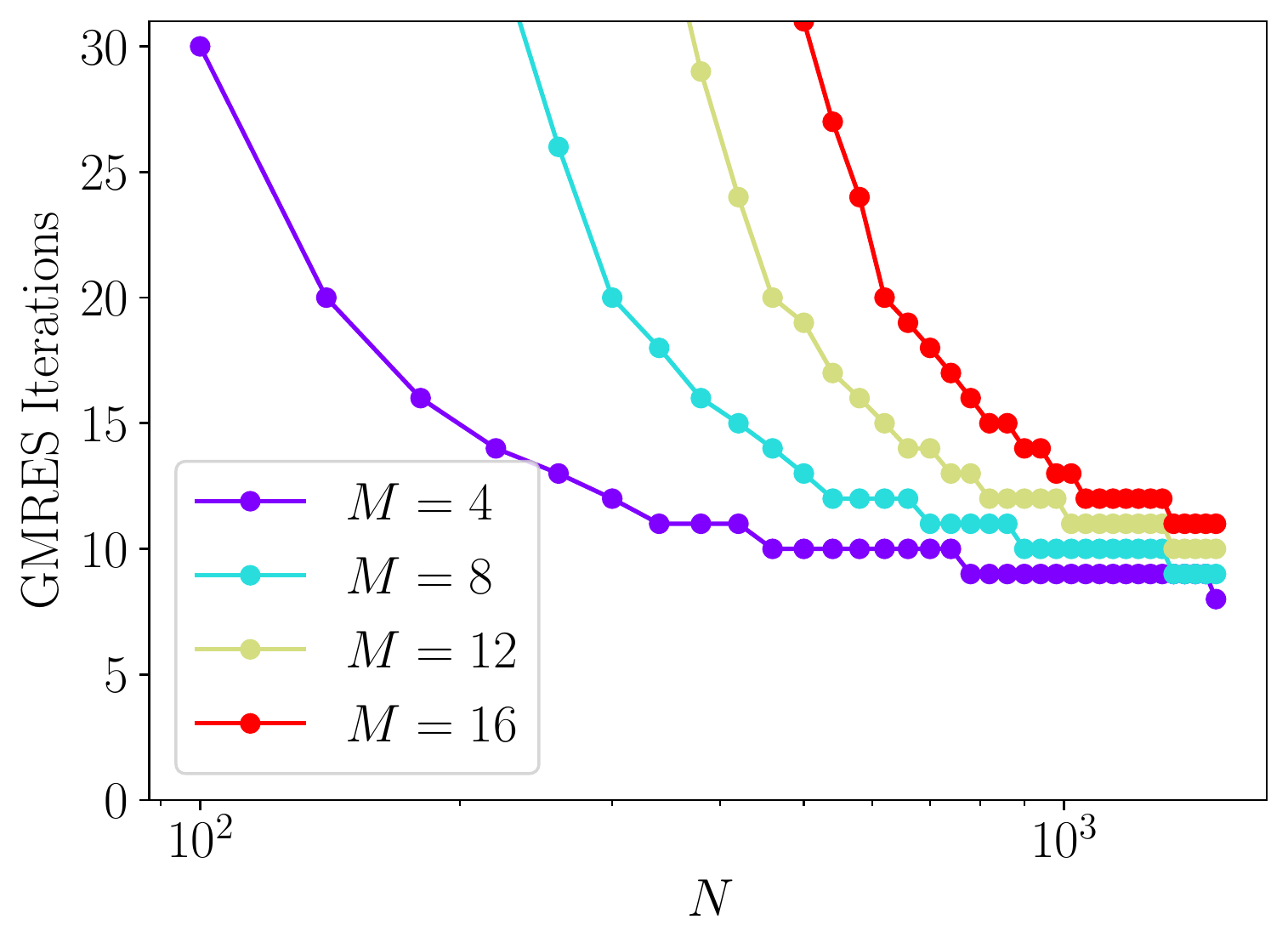}
    \caption{}
  \end{subfigure}
  \caption{Convergence study in boundary points $N$ at fixed $M$. Panel (a) shows errors in $L^\infty(\Omega)$, with matching $M$th order convergence lines (dashed). The dashed gray lines are proportional to $1/N^2$, and the dotted gray line shows the $10^{-14}$ tolerance set on all approximate or iterative portions of the solver. Panel (b) shows the number of GMRES iterations used by the Annular solver to reach a relative residual $<10^{-14}$.}
  \label{figure:fixed_m}
\end{figure}

We begin with a simple exploration of the properties of this solver for a fixed number of Chebyshev modes $M$ and a variable radius $R=Mh_\text{min}$, with $h_\text{min}$ the smallest discrete boundary gridspacing, as given in \Cref{section:methods_preliminaries:annular_domain}. We vary $N$ from $100$ to $1500$, tracking both the error, in $L^\infty(\Omega)$, and the number of GMRES iterations required to invert the annular problem, for $M=4$, $8$, $12$, and $16$. The regular gridspacing $h$ is set to be $h_\textnormal{min}/2$, $b=\lceil1.5M\rceil$, and $\epsilon=10^{-14}$. The error and number of iterations are shown in panels (a) and (b) of \Cref{figure:fixed_m}, respectively. There are several points worth making about this simple study.
\begin{enumerate}
	\item Solutions fail to exist (or have large errors) for small $N$ when $M$ is large. This is because the $R<R_\text{max}$ criteria is not obeyed by the configuration: the implied coordinates would be singular. When the criteria is just \emph{barely} obeyed, a large number of GMRES iterations are required to invert the annular problem. The iteration count decays rapidly as $N$ is increased, to a nearly $M$ independent number of $\approx10$.
	\item Higher $M$ indeed leads to faster convergence --- the accompanying dashed lines are $M$th order convergence lines.
	\item Rapid convergence \emph{stagnates} at a certain error: for $M=4$ this is off the graph, but for larger $M$ this effect is clearly apparent. This stagnation occurs when the dominant error in the problem is the resolution of the cutoff function $\eta$\footnote{Note that when $M$ is fixed, and the regular gridspacing is set proportional to the boundary gridspacing, $\eta$ varies from $0$ to $1$ over the same number of regular gridpoints regardless of $N$.}. Empirically, we observe continued second-order convergence once this floor has been reached\footnote{The continued second order convergence once occurs because the cutoff function $\eta$ multiplies $f$, but we are solving for $u$, which is two derivatives smoother.}.
	\item The horizontal gray dotted line is placed at $\epsilon=10^{-14}$. Convergence to a small multiple of $\epsilon$ is observed for the $M=16$ convergence curve.
\end{enumerate}
We do not belabor the first or the third point because they both have the same solution: letting the number of Chebyshev modes $M$ increase proportionally with $N$. In doing so, the physical width of $\A$ can be fixed geometrically so that it embeds without coordinate collapse (see \Cref{section:methods_preliminaries:annular_domain}). Convergence of the annular solution is now expected to be spectral in $N$, and because $R$ remains fixed but $N_x$ and $N_y$ increase, the regularized step function $\eta$ becomes progressively better resolved, eliminating the stagnation observed when fixed values of $M$ are used.

\subsection{Behavior for proportional $M$}
\label{section:simple_poisson:proportional_M}

\begin{figure}[h!]
  \centering
  \begin{subfigure}[c]{0.4\textwidth}
    \centering
    \includegraphics[width=\textwidth]{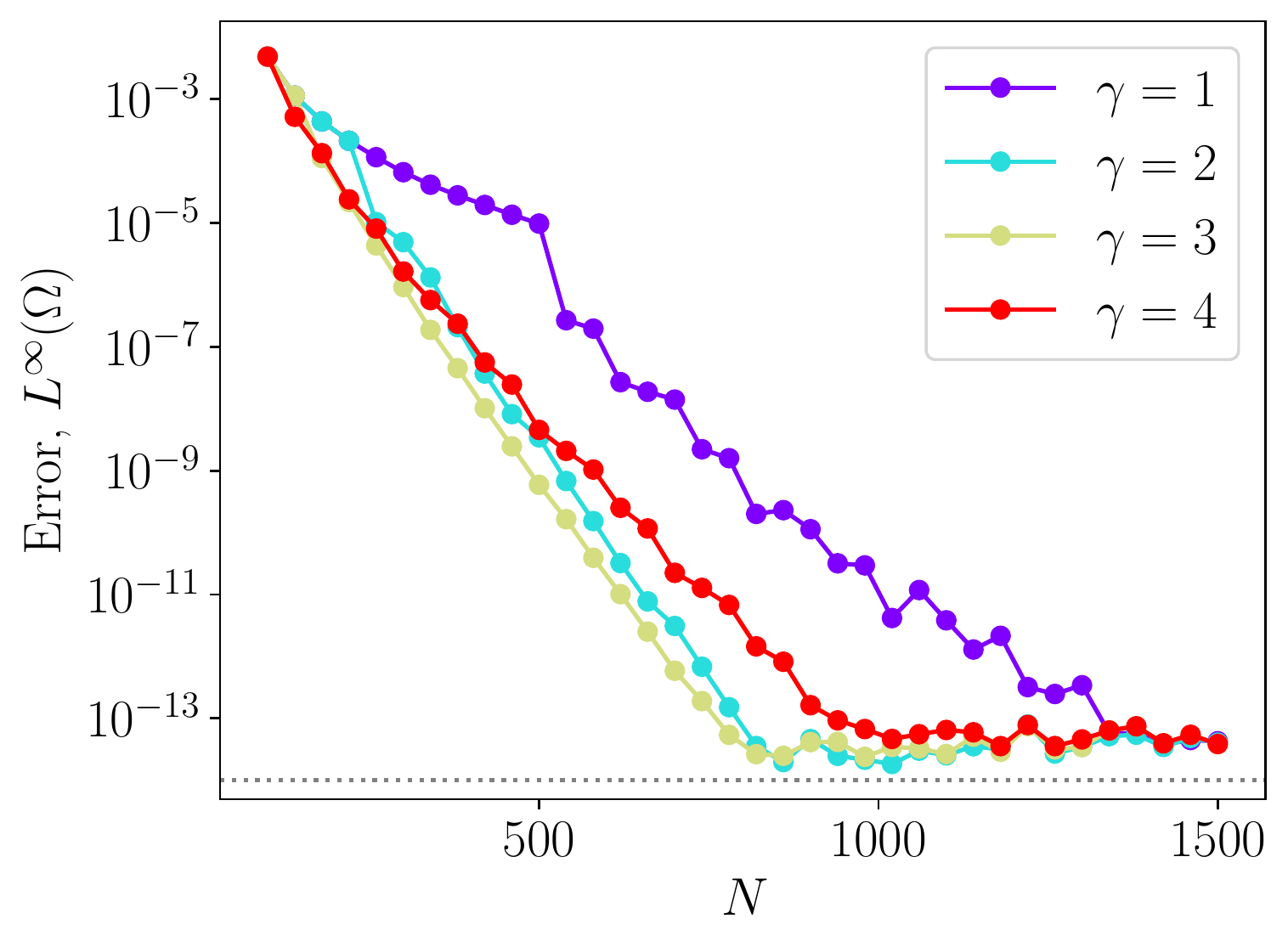}
    \caption{}
  \end{subfigure}
  \hspace{1em}
  \begin{subfigure}[c]{0.4\textwidth}
    \centering
    \includegraphics[width=\textwidth]{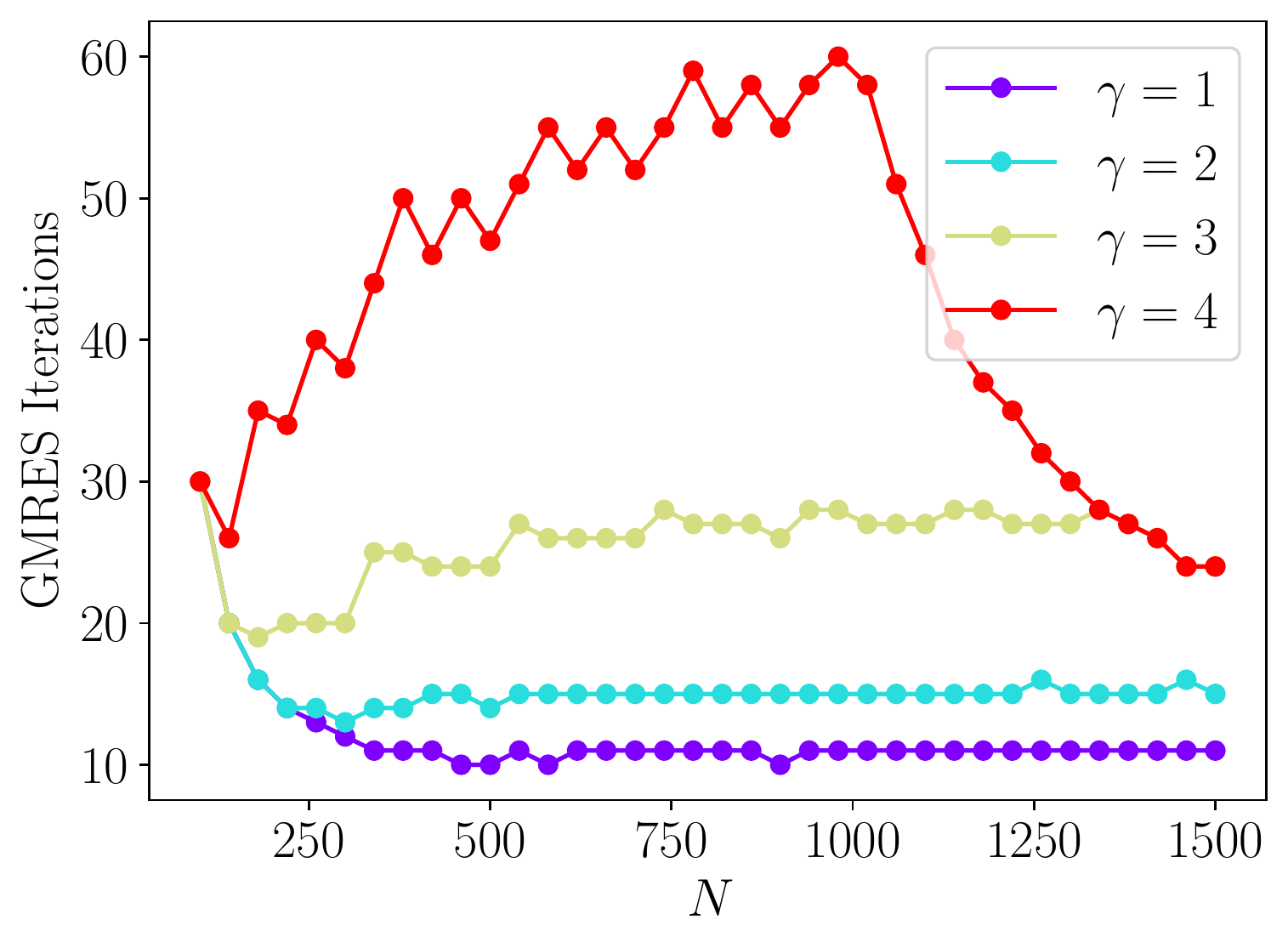}
    \caption{}
  \end{subfigure}
  \caption{Convergence study in boundary points $N$ with fixed $\gamma$ ($M$ proportional to $N$). Panel (a) shows errors in $L^\infty(\Omega)$; panel (b) shows the number of GMRES iterations used by the annular solver to reach a relative residual $<10^{-14}$.}
  \label{figure:fixed_gamma}
\end{figure}

We now redo this refinement study, with the same parameter choices as in \Cref{section:simple_poisson:fixed_M} but now scaling $M$ as a function of $N$. To be precise, we will take $M=\lfloor\gamma N/100\rfloor$, for $\gamma=1$, $2$, $3$, and $4$, restricting $M$ to a minimal value of $4$ and a maximal value of $40$. \Cref{figure:fixed_gamma} shows the $L^\infty(\Omega)$ error along with the number of GMRES iterations used in the annular solve. Notice that the $x$-axis is now linear; as expected, convergence is spectral in $N$. Using smaller values of $\gamma$ leads to lower iteration counts but a slower rate of exponential convergence; higher values typically lead to faster convergence, although beyond a certain point this advantage becomes negligible as other errors (or simply resolving $f$) begin to dominate. As in the $M=16$ case, errors saturate at a small multiple of $10^{-14}$. The jumps in the iteration counts and the staircase effect in the error (especially for $\gamma=1$) are due to integer shifts in the value of $M$ as $N$ is increased.

\section{Parameter selection}
\label{section:parameter_selection}
We are now ready to discuss how to select parameters. We suppose that the user provides a parametrized curve $\Gamma$ that bounds the domain $\Omega$, along with an inhomogeneity $f$ and boundary inhomogeneity $g$. For simplicity, we assume that $\Gamma$, $f$, and $g$ are given as evaluable functions anywhere in their respective domains. The choices we must make are given in \Cref{table:parameter_choices}. There are many ways these choices could be made, and in a highly optimized numerical solver they would likely be chosen to minimize the error given a runtime constraint. We will choose a simpler method --- enforcing that every element of the solver is resolved to a user-specified length-scale $h$.

\subsection{Selection of the annular width $R$}
\label{section:parameter_selection:annular_width}

Having observed that the maximal annular width $R_\text{max}$ is purely geometric (and not a function of the discretization or parametrization), we are free to choose $R$ before discretization is considered. The primary tradeoff here is efficiency vs. speed: the choice of $R$ is equivalent to choosing $\gamma$ in \Cref{section:simple_poisson:proportional_M}, and so the rate of exponential convergence must be traded off against an increased number of GMRES iterations. We have found that a good default is $R=R_\text{max}/2$; and this will be used throughout the remainder of the paper.

\subsection{Spatial discretization}
\label{section:parameter_selection:spatial}

We now turn to the spatial discretization. The simplest possibility is to allow the user to specify a lengthscale $h$ as the basic unit of parametrization, and to ensure that all components of the solver resolve this lengthscale. In particular, we take:
\begin{enumerate}
	\item (Selection of number of boundary points $N$): We remind the reader that the Jacobian of the coordinate transformation is $\mathcal{J}=\phi(s)(1-r\kappa(s))$, where $\phi$ and $\kappa$ are the speed and curvature, respectively. We can think of the term $(1-r\kappa(s))$ as an amplification factor for the boundary discretization. Letting $h^j_\Gamma=\phi(s_j)\Delta s$, where $\Delta s=2\pi/N$ and $N$ is the number of boundary points, we can compute that the azimuthal grid-spacing at the interface is $h_\I=\phi(s_j)(1-r\kappa(s_j))\Delta s$. We thus choose $N$ to be the smallest even integer such that $\min_j(h^j_\Gamma,\ h^j_\I)<h$; with it being clear that the grid-spacing at $\I$ and $\Gamma$ are sufficient to bound the azimuthal spacing everywhere in $\A$.
	\item (Selection of the number of Chebyshev modes $M$):  The effective radial resolution will be $h_r=\pi R/(2M)$. We select $M$ to be the smallest integer with $h_r<h$; i.e. $M=\lceil\frac{\pi R}{2h}\rceil$.
	\item Bounds for the grid are computed as discussed in \Cref{section:methods_preliminaries:computational_domain_discretization}, and $N_x$ and $N_y$ are chosen to be the minimal even integers with $h_x=h_y<h$, in accordance with the considerations discussed in \Cref{section:methods_preliminaries:computational_domain_discretization,section:methods_preliminaries:peridoic_compatibility}.
\end{enumerate}

\subsection{Choice of the bandwidth $b$ defining $\eta$}
\label{section:parameter_selection:eta}

\begin{figure}[h!]
  \centering
  \begin{subfigure}[c]{0.4\textwidth}
    \centering
    \includegraphics[width=\textwidth]{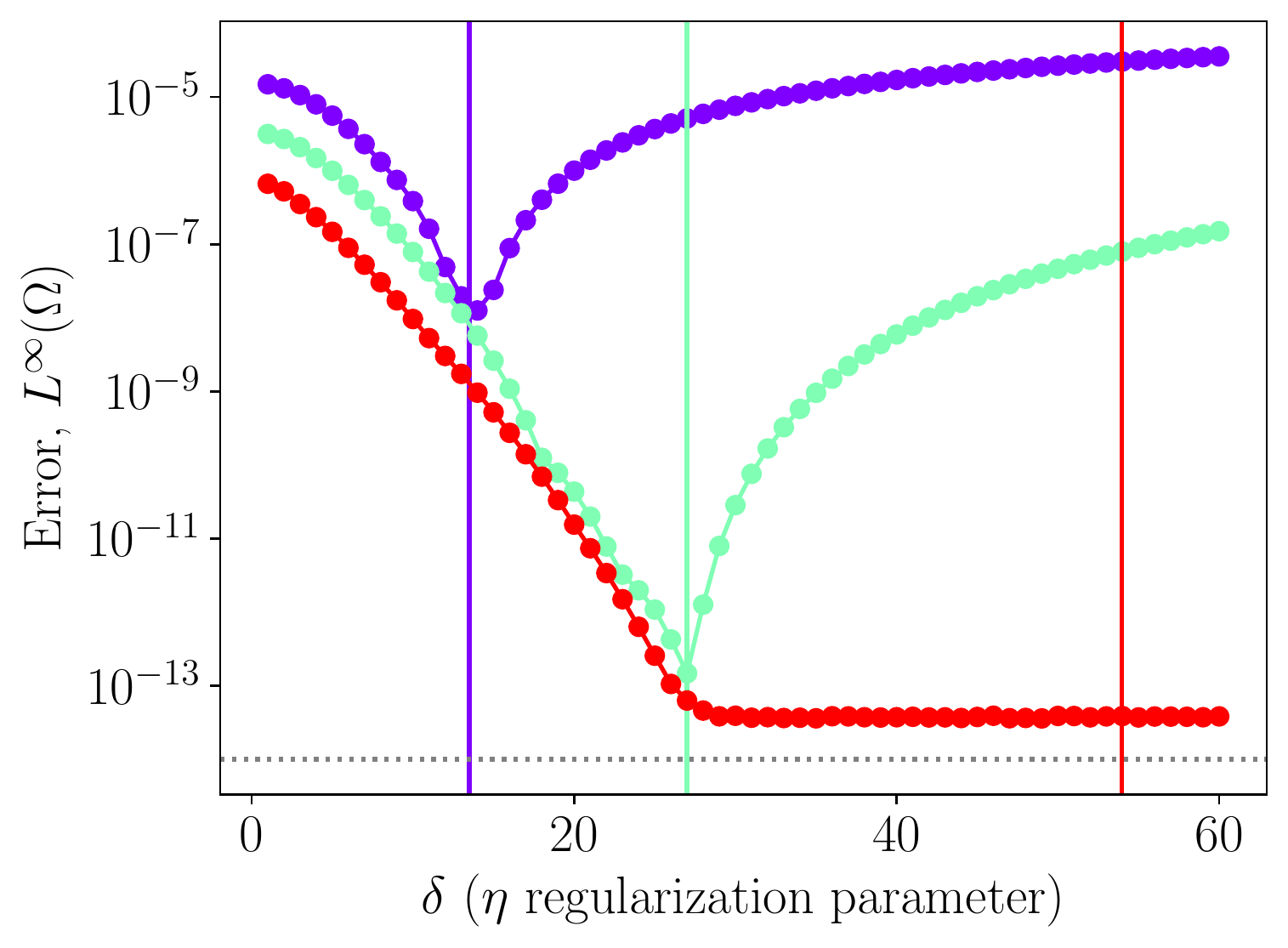}
    \caption{}
  \end{subfigure}
  \hspace{1em}
  \begin{subfigure}[c]{0.4\textwidth}
    \centering
    \includegraphics[width=\textwidth]{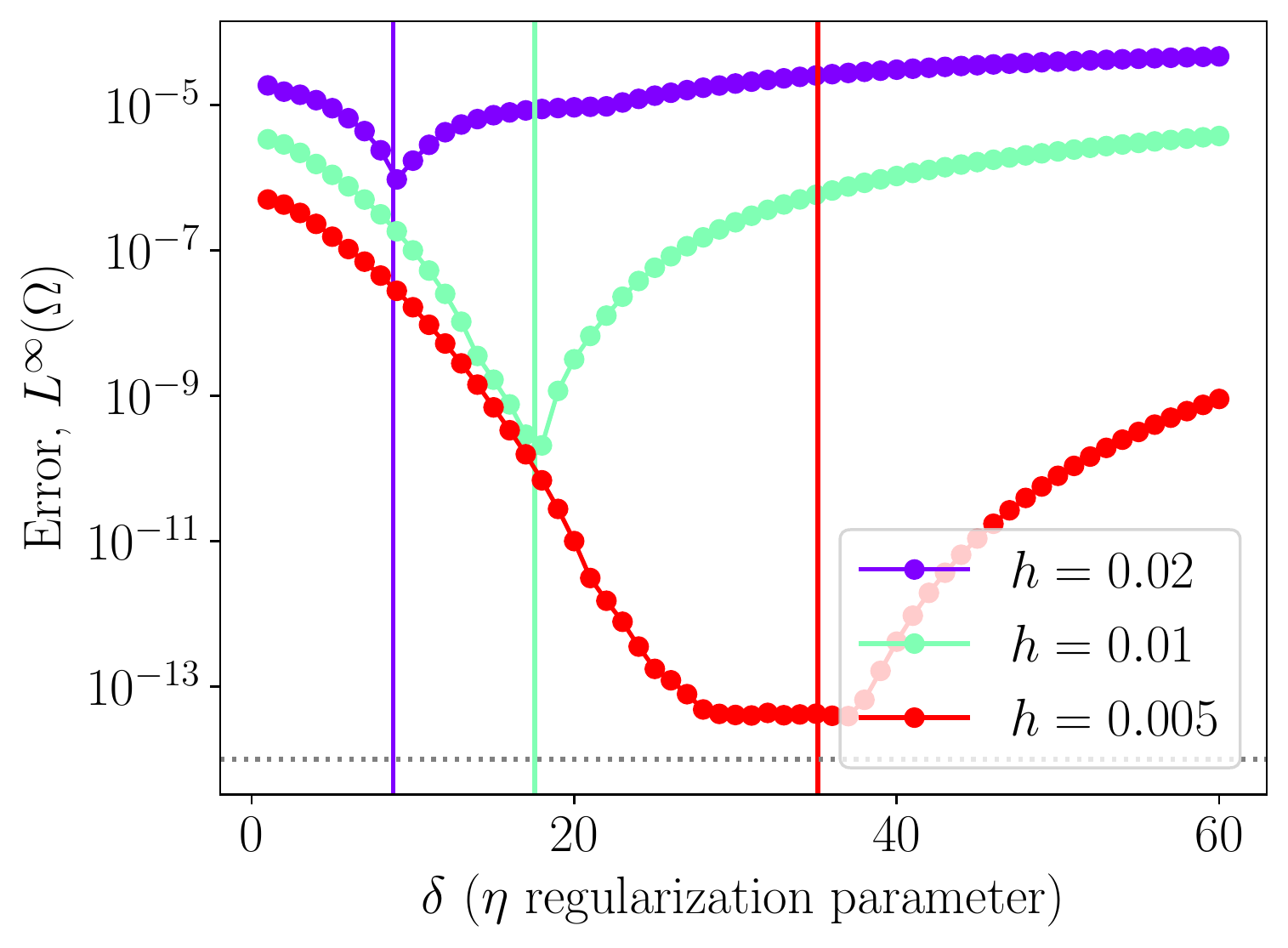}
    \caption{}
  \end{subfigure}
  \caption{Error as a function of $\delta$ (the regularization parameter for $\eta$). In both panels, different colors give different values of $h$, the underlying discretization scale, with purple a coarse discretization; green an intermediate discretization, and red a fine discretization. Solid vertical lines show the estimate $2R/h$ that is used to set $\delta$ in our simulations. Panel (a) show results computed on the relatively smooth function and simple domain used in \Cref{section:simple_poisson}; Panel (b) shows results for the more complicated function and domain used in \Cref{section:PUX_comparison}.}
  \label{figure:regularization_study}
\end{figure}

The final choice that remains to be made is the regularization parameter defining $\eta$. For the Poisson problem, the regularization parameter $\lceil2R/h\rceil$ provides near-optimal results. The best choice of this parameter will, in general, be PDE dependent, and will especially differ from our choice in near-identity problems (such as high-$k$ Helmholtz or modified Helmholtz); where the smoothing feature of the underlying elliptic operator will be less apparent. In \Cref{figure:regularization_study}, we show the error as a function of the regularization parameter $\delta$ that defines $\eta$ for two very different Poisson problems set on different domains over a range of underlying discretization parameters $h$. For both problems and all values of $h$, the estimate $2R/h$ provides a remarkably good estimate of the best observed value.

\subsection{Numerical timing and stability}
\begin{figure}[h!]
  \centering
  \begin{subfigure}[c]{0.32\textwidth}
    \centering
    \includegraphics[width=\textwidth]{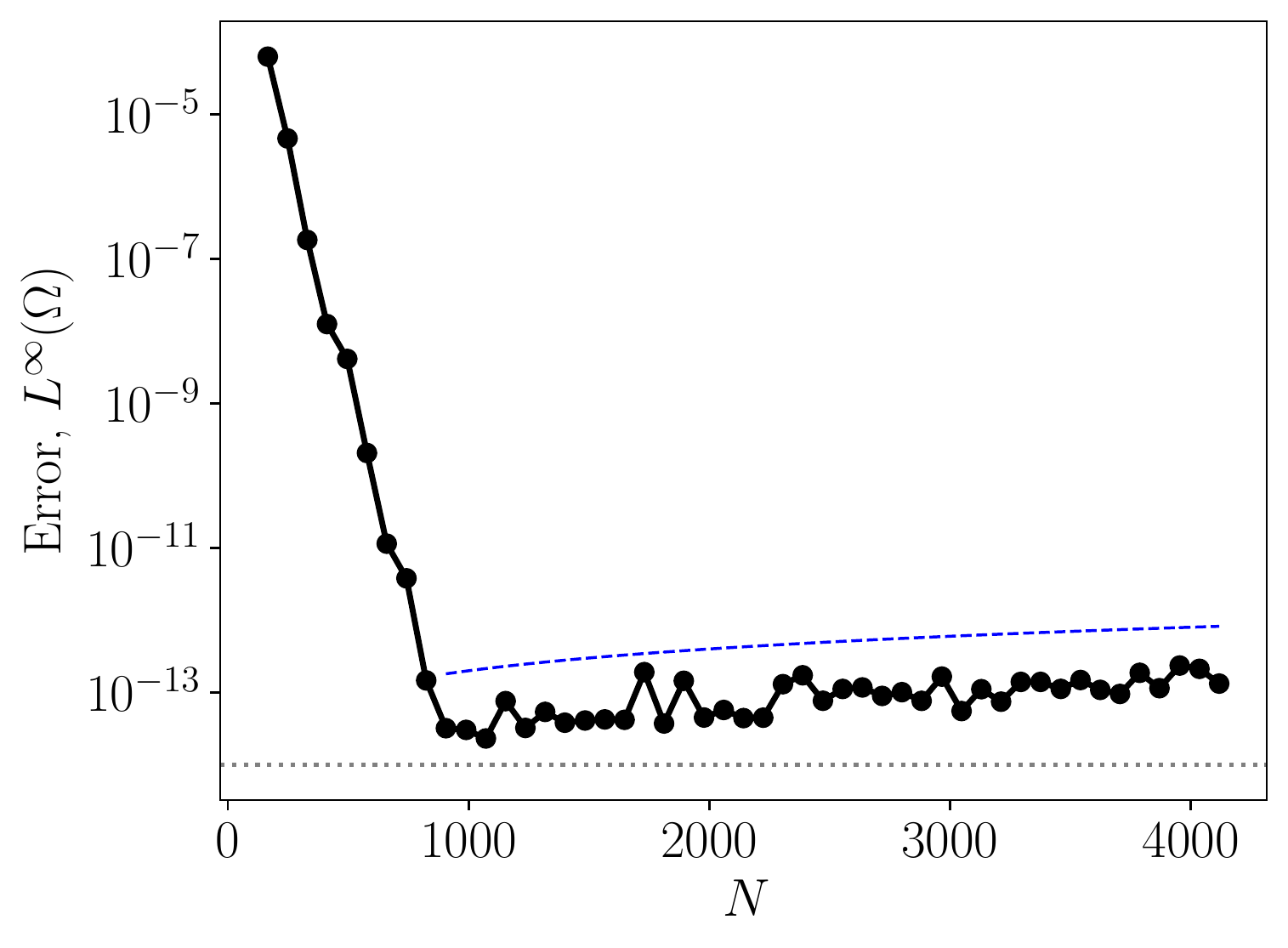}
    \caption{}
  \end{subfigure}
  \hfill
  \begin{subfigure}[c]{0.32\textwidth}
    \centering
    \includegraphics[width=\textwidth]{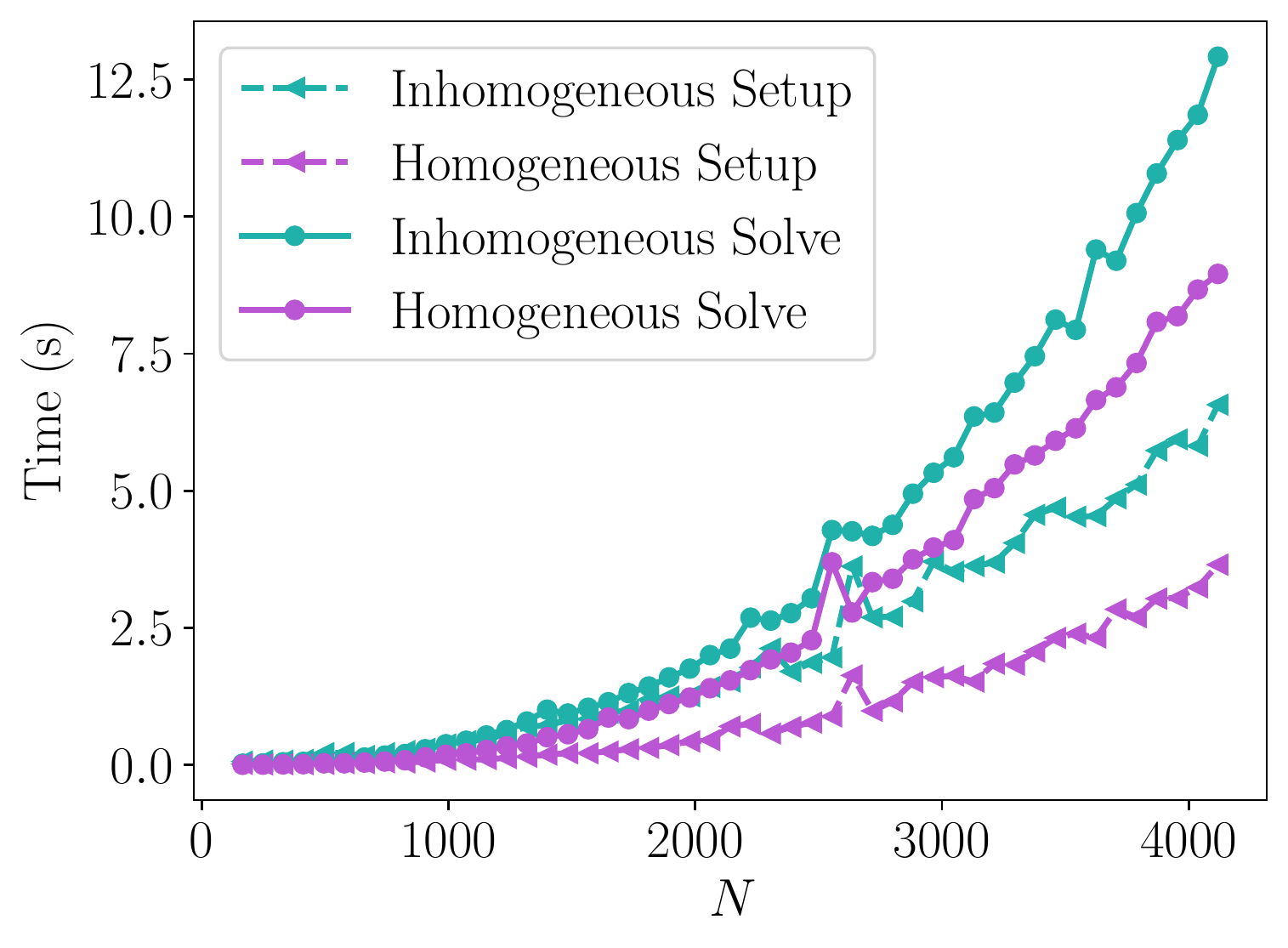}
    \caption{Macbook Pro}
  \end{subfigure}
  \hfill
  \begin{subfigure}[c]{0.32\textwidth}
    \centering
    \includegraphics[width=\textwidth]{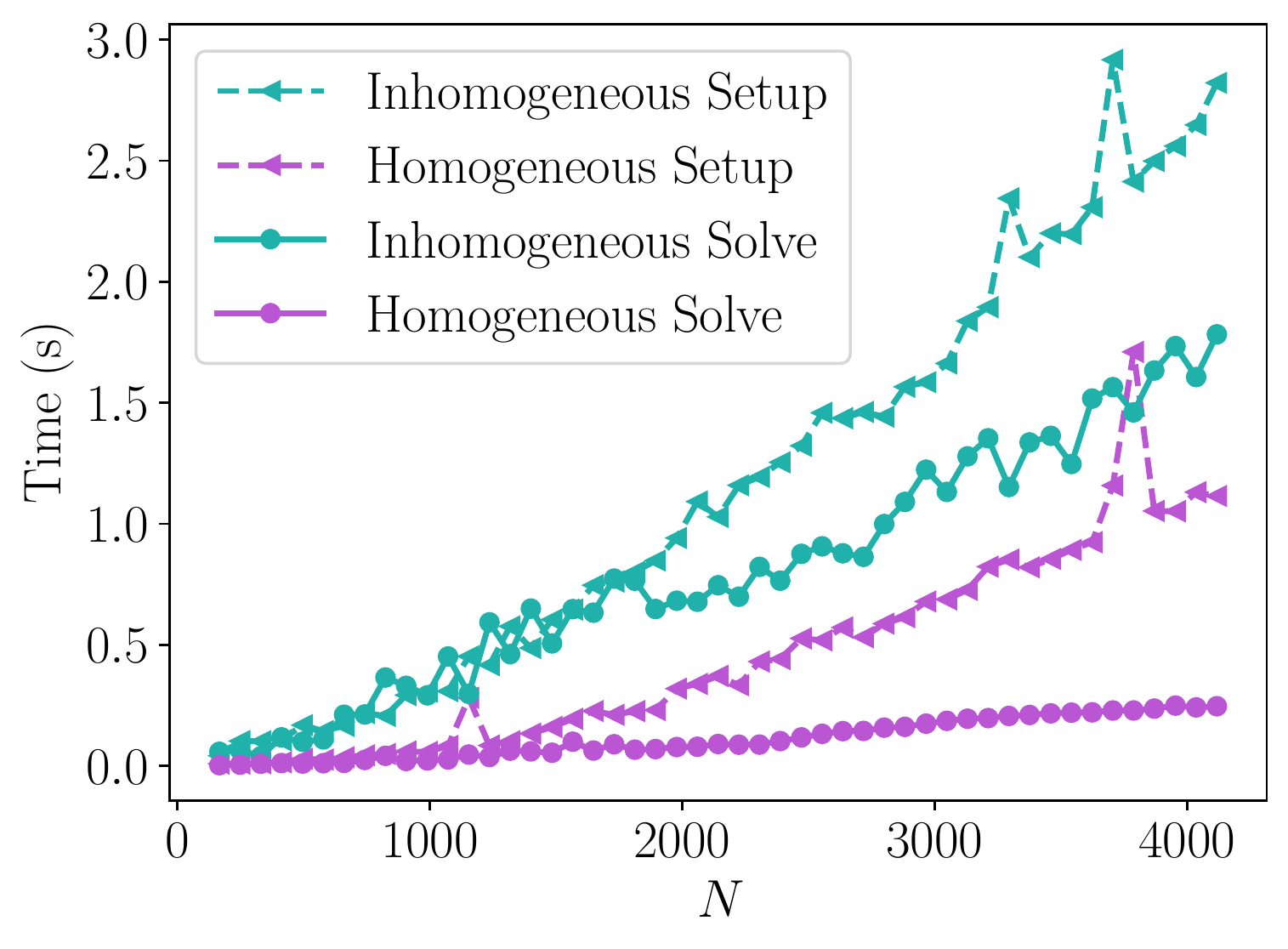}
    \caption{Intel Node}
  \end{subfigure}
  \caption{Results for the simple example problem studied in \Cref{section:simple_poisson}, using the parameter choices discussed in \Cref{section:parameter_selection}, over a large range of spatial discretization scales $h$. Panel (a) shows errors as a function of $N$ (note that $N\propto h^{-1}$).  The dotted gray line is set at $\epsilon=10^{-14}$. The dashed blue line is proportional to $N$, showing the expected (slow) loss of error as $N$ grows. Panels (b) and (c) show wall clock times on a Macbook Pro and Intel Cluster node, respectively, broken down into setup and solve times for both the inhomogeneous portion of the solve and the homogeneous correction.}
  \label{figure:error_and_timings}
\end{figure}

Finally, we rerun the problem from \Cref{section:simple_poisson} using the parameter choices outlined in this section, over a large range of $h$, examining both numerical timings and stability as the number of boundary points $N$ grows large. The error, in $L^\infty(\Omega)$, is shown in \Cref{figure:error_and_timings}(a) for $h=0.002$ to $h=0.05$ (corresponding to very fine discretizations with $N\approx 4000$ to very coarse discretizations with $N<200$). As expected, we see rapid and exponential convergence of the error to a small multiple of $\epsilon=10^{-14}$. As $N$ gets larger, the error stays relatively stable, with a very slow loss of accuracy, proportional to $N$ (shown as the dashed blue curve). This mild loss of accuracy is due to the direct computation of the gradient of $u$ used to match interface derivatives (see \Cref{section:methods_specifics:stiching}); this could be remedied by exploiting other methods to estimate the off-grid derivative.

Panels (b) and (c) of \Cref{figure:error_and_timings} show walk clock timings, broken into both setup (dashed lines) and solve (solid lines) and inhomogeneous (green) vs. homogeneous corrections (purple), on two different computers. Our implementation is not optimal but wall clock times are included to demonstrate what is practically achievable with a carefully implemented but not aggressively optimized code written in Python. The timings in panel (b) are from a quad-core Macbook Pro with a single Intel(R) Core(TM) i7-8569U CPU @ 2.80GHz and 16 GB of RAM, the timings in panel (c) are from a 40-core cluster node with two Intel(R) Xeon(R) Gold 6148 CPU @ 2.40GHz and 768GB of RAM. Timings are broken down into both setup portions (dependent on $\Gamma$ and $h$, but not on $f$ or $g$), and solve portions (dependent on $f$ and $g$). For small problem sizes (e.g. $N\approx 200$), setup and solve are both done on the 10s of milli-second timescale. When $N=2000$, solutions are produced in $<1$ second, with timings in this implementation dominated by the annular solve and calls to the FMM; setup time is $<2s$, with timing dominated by factorization of the dense homogeneous correction and QFS matrices for the close evaluation of layer-potentials. In \Cref{section:discussion}, we discuss future implementational improvements that could reduce the cost of the method.

\section{Summary of algorithm}
\label{section:recombobulation}

Finally, we present summaries of the full algorithm, with computational costs for our specific implementation given for all compute-heavy steps. We report asymptotic scalings in only $N$ and $M$, taking $N_x\propto N_y\propto N$. The algorithm for the setup portion of the solve, which depends on the domain $\Omega$, the user-specified resolution scale $h$, and tolerance $\epsilon$, is given in \Cref{algorithm:setup}; the total scaling for this stage, in our implementation, is $\mathcal{O}(N^3) + \mathcal{O}(NM^3)$. The algorithm for the solve stage, which depends on $f$, $g$, and $\epsilon$, is given in \Cref{algorithm:solve}; with a total scaling of $\mathcal{O}(N^2\log N)+\mathcal{O}(N M^2)+\mathcal{O}(NM\log N\log M)$. Recall that $N$ is the number of \emph{boundary} nodes, and so the total number of unknowns is $\mathcal{O}(N^2)$, with the complexity of the FFT used to solve the regular problem setting our benchmark scaling of $\mathcal{O}(N^2\log N)$.

\begin{algorithm}[H]
  \algrenewcommand\algorithmiccomment[2][\footnotesize]{{#1\hfill\(\triangleright\) #2}}  
  \caption{Setup procedure for domain $\Omega$ \hspace*{\fill} $\mathcal{O}(N^3 + NM^3)$}
  \algorithmicrequire{
  	Smooth parametrized coordinates $\X$ for the curve $\Gamma=\partial\Omega$, length-scale $h$, tolerance $\epsilon$.
  }
  \begin{algorithmic}[1]
  	\State compute $R_\textnormal{max}$; $R=R_\textnormal{max}/2$
  		\Comment{see \Cref{section:methods_preliminaries:annular_domain}}
    \State select number of boundary nodes $N$, Chebyshev nodes $M$, and regular grid nodes $N_x$, $N_y$, so all discretizations will resolve lengthscale $h$
    	\Comment{see \Cref{section:parameter_selection:spatial}}
    \If{compatibility conditions are required}
    	\State{augment $\mathcal{C}$ to provide space for bump function}
    	\Comment{see \Cref{section:methods_preliminaries:peridoic_compatibility}}
    \EndIf
    \State Discretize $\mathcal{C}$, $\Gamma$, $\A$
   		\Comment{see \Cref{section:methods_preliminaries:computational_domain_discretization,section:methods_preliminaries:annular_domain_discretization}}
    \State [$\mathcal{O}(N^2\log N)$] Categorize discrete nodes of $\mathcal{C}$ and compute cutoff $\eta$
    	\Comment{see \Cref{section:methods_preliminaries:cutoff,section:methods_preliminaries:coordinates}}
    \If{compatibility conditions are required}
    	\State{compute bump function $\xi$}
    	\Comment{see \Cref{section:methods_preliminaries:peridoic_compatibility}}
    \EndIf
   	\State [$\mathcal{O}(N^3)$] Compute inward QFS representation for $\Gamma$
   		\Comment{see \Cref{section:methods_preliminaries:close_evaluation}}
   	\State [$\mathcal{O}(N^3)$] Compute inward and outward QFS representations for $\I$.
   		\Comment{see \Cref{section:methods_preliminaries:close_evaluation}}
   	\State [$\mathcal{O}(N^3)$] Compute and factor singular matrix $A$
   		\Comment{see \Cref{eq:discretized_homogeneous}}
   	\State [$\mathcal{O}(NM^3)$] Construct preconditioner for annular solver
   		\Comment{see \Cref{section:methods_specifics:annular_problem}}
  \end{algorithmic}
  \label{algorithm:setup}
\end{algorithm}

\begin{algorithm}[H]
  \algrenewcommand\algorithmiccomment[2][\footnotesize]{{#1\hfill\(\triangleright\) #2}}  
  \caption{Solution procedure given $f$, $g$\hspace*{\fill}$\mathcal{O}(N^2\log N + NM^2 + NM\log N\log M)$}
  \algorithmicrequire{
  	Setup procedure as given in \Cref{algorithm:setup}; evaluatable functions $f$ and $g$ \textbf{or} discrete values of (1) $f$ known at the nodes of $\mathcal{C}$ in $\fOmega$ and nodes of $\A$ and (2) $g$ known at nodes of $\Gamma$; tolerance $\epsilon$.
  }
  \begin{algorithmic}[1]
  	\State [$\mathcal{O}(N^2)$] compute $\eta f$ at all discrete nodes of $\mathcal{C}$
  		\Comment{see \Cref{section:methods_specifics:regular_problem}}
    \If{compatibility conditions are required}
    	\State [$\mathcal{O}(N^2)$] adjust $\eta f$ to have 0 mean
    	\Comment{see \Cref{section:methods_preliminaries:peridoic_compatibility}}
    \EndIf
    \State [$\mathcal{O}(N^2)$] solve regular problem to get $\ur$ at all discrete nodes of $\mathcal{C}$
    	\Comment{see \Cref{section:methods_specifics:regular_problem}}
    \State [$\mathcal{O}(NM^2)$] solve annular problem to get $\ua$ at all discrete nodes of $\A$
    	\Comment{see \Cref{section:methods_specifics:annular_problem}}
    \State [$\mathcal{O}(N^2\log N)+\mathcal{O}(NM^2)$] solve stitching problem to get $\ui$
    	\Comment{see \Cref{section:methods_specifics:stiching,algorithm:stitch}}
    \State [$\mathcal{O}(N M^2)$] Compute boundary discrepancy at $\Gamma$
    	\Comment{see \Cref{section:methods_specifics:homogeneous}}
    \State [$\mathcal{O}(N^2)$] Solve BIE problem and add $\uh$ to $\ui$ to find $\u$
    	\Comment{see \Cref{section:methods_specifics:homogeneous}}
    \State [$\mathcal{O}(NM\log N\log M)$] Compute $\u$ at all discrete nodes of $\mathcal{C}$ in $\A$
    	\Comment{see \Cref{section:methods_specifics:finishing}}
  \end{algorithmic}
  \label{algorithm:solve}
\end{algorithm}

\section{A multi-body problem with high-frequency $f$ and comparison to PUX}
\label{section:PUX_comparison}

In this section, we extend our method to multiply connected domains, and solve a more complex problem previously solved in \cite{fryklund2018partition}, comparing errors with those generated by the high-order \emph{Partition of Unity Function Extension Method} (PUX). The definition of the domain and oscillatory RHS forcing are given in \cite{fryklund2018partition} and are not repeated here. 

\subsection{Considerations for multiply-connected domains}

Extension to multiply connected domains, with $K$ well-separated bodies where the annular domains do not overlap, is relatively straightforward and so we lay out only the relevant considerations here.

\begin{enumerate}
	\item Annular regions $\A_i$ are defined independently for each boundary, and the choice of $R_i$, $N_i$ and $M_i$ is made independently to be consistent with the single scale $h$, as described in \Cref{section:parameter_selection}.
	\item The faithful domain is defined as $\fOmega=\Omega\setminus(\cup_i\A_i)$.
	\item A single function $\eta$ is again computed, with $\eta(\x)=1$ for all $\x\in\fOmega$, $\eta(\x)=0\in\Omega^C$, and transitioning in the same way over each radial region $\A$. Note that regularization parameter $\lceil 2R_i/h\rceil$ for each region will typically be different. Errors will often be controlled by the boundary with the smallest $R$; this is somewhat inevitable when using global discretizations, as in this paper.
	\item The annular solutions $u_{\A_i}$ are computed independently.
	\item Solution of the homogeneous problem, and evaluation of the solution in the multi-boundary context is described in \cite{stein2021quadrature}.
\end{enumerate}

The stitching problem in the multi-body case is straightforward, but some simple optimizations exist when utilizing our specific implementation which make its computation more efficient, and we describe these here. Recall that the goal is to evaluate the layer potential given in \Cref{eq:laplace_stitching} for all discrete nodes of $\mathcal{C}\in\fOmega$ and all discrete nodes of $\A_i$ for each body. The first step is to compute all inward sources $\I^i_\textnormal{in}$ along with all inward effective potentials $\zeta^i_\textnormal{in}$. Using a single FMM, these can be evaluated at all $\x\in\fOmega$ and on all interfaces $\I^i$:
\begin{equation}
	v(\x) = \sum_{i=1}^{K}\frac{-1}{2\pi}\sum_{j=1}^{N_s^i} \log|\x-\X^{\I^i_\textnormal{in}}(s^i_j)|\zeta^i_\textnormal{in}(s^i_j)w^i_j,
\end{equation}
which requires only $\mathcal{O}(KN)+\mathcal{O}(N_\fOmega+KN)$ operations (with $N_\fOmega$ the number of discrete nodes of $\mathcal{C}$ within $\fOmega$, and $N_s^i$ the number of source nodes on $\I_i$, with $N_s^i\propto N$ for all $i$). It remains to evaluate these potentials at each annulus $\A_i$. We first independetly compute outward sources $\I^i_\textnormal{out}$ and effective potentials $\zeta^i_\textnormal{out}$ for each body, and on each interface we now compute an adjusted $v$:
\begin{equation}
	v^i_{adj}(\x) = v(\x) + \frac{1}{2\pi}\sum_{j=1}^{N_s^i} \log|\x-\X^{\I^i_\textnormal{in}}(s^i_j)|\zeta^i_\textnormal{in}(s^i_j)w^i_j,
\end{equation}
which is the layer potential at $\I_i$ generated by all inward sources other than the $i$th source itself. Again, using the methodology described in \cite{stein2021quadrature}, we compute a second outward potential $\xi^i$ for each body so that:
\begin{equation}
	\mathcal{S}_{\mathcal{I}^i_\textnormal{out}}\xi^i = v^i_{adj},
\end{equation}
for all $\x\in\I^i$. We may now, independently for each body, compute:
\begin{equation}
	\nu^i(\x) = \sum_{i=1}^K\left(\mathcal{S}_{\I^i}\sigma^i - \mathcal{D}_{\I^i}\gamma^i\right)|_{\A^i}(\x) = \sum_{j=1}^{N_s^i} \log|\x-\X^{\I^i_\textnormal{out}}(s^i_j)|\left[\zeta^i_\textnormal{out}(s^i_j)+\xi^i(s^i_j)\right]w^i_j,
\end{equation}
valid for any $\x\in\A^i$. Note that this last step requires $K$ independent FMM calls, for a total cost of $\mathcal{O}(KNM)$; these calls are embarassingly parallel. Finally, we define the particular solution $u_I$ in analogy to \Cref{eq:particular}, as:
\begin{equation}
	\ui(\x) = 
	\begin{cases}
		\ur(\x) + v(\x) &\quad\text{for }\x\in\fOmega, \\
		\u_{\A_i}(\x) + \nu^i(\x) &\quad\text{for }\x\in\aOmega_i.
	\end{cases}
	\label{eq:particular}
\end{equation}

\subsection{Comparison to PUX}

We solve this problem over a range of $h$, although results are reported with respect to the number of gridpoints discretizing $\mathcal{C}$ in the $x$-direction ($N_x$), for consistency with \cite{fryklund2018partition}. Results from that manuscript have been extracted \cite{Rohatgi2020} and multiplied by $\|u\|_{L^\infty}\approx 1.35$, as we report absolute errors here. We see that both methods converge rapidly to low error; although our method typically delivers slightly better accuracy for a given $N_x$. We make several comments, however, regarding the direct comparability of the results.
\begin{enumerate}
	\item The parameters and functional form of the inner boundary used in \cite{fryklund2018partition} do not seem to generate exactly the domain shown in their manuscript. We have instead changed the non-zero values of $c$ and $d$ to be $c_3=c_6=0.1$, $d_{-3}=0.05$ which generates a qualitatively similar domain (see \cite{fryklund2018partition} for details). This may shift the error curves to some degree.
	\item The results reported in the PUX method used a fixed, and very fine, discretization for the boundary (as well as the number of partitions used for the extensions), varying only the gridspacing of the underlying grid. In our results, the boundary shares a discretization scale $h$ consistent with the background grid, as defined in \Cref{section:parameter_selection}.
	\item The PUX method generates solutions when $N$ is very small (say, $100$). Our method (at least with the standard way of setting parameters given in \Cref{section:parameter_selection}), is unable to generate solutions here as the implied $M$ would be less than $4$, the minimal $M$ we allow.
\end{enumerate}
While direct comparison is not possible, it is clear that both methods produce discretizations for this problem that converge rapidly as they are refined, and deliver comparable errors.

\begin{figure}[h!]
  \centering
  \begin{subfigure}[c]{0.38\textwidth}
    \centering
    \includegraphics[width=\textwidth]{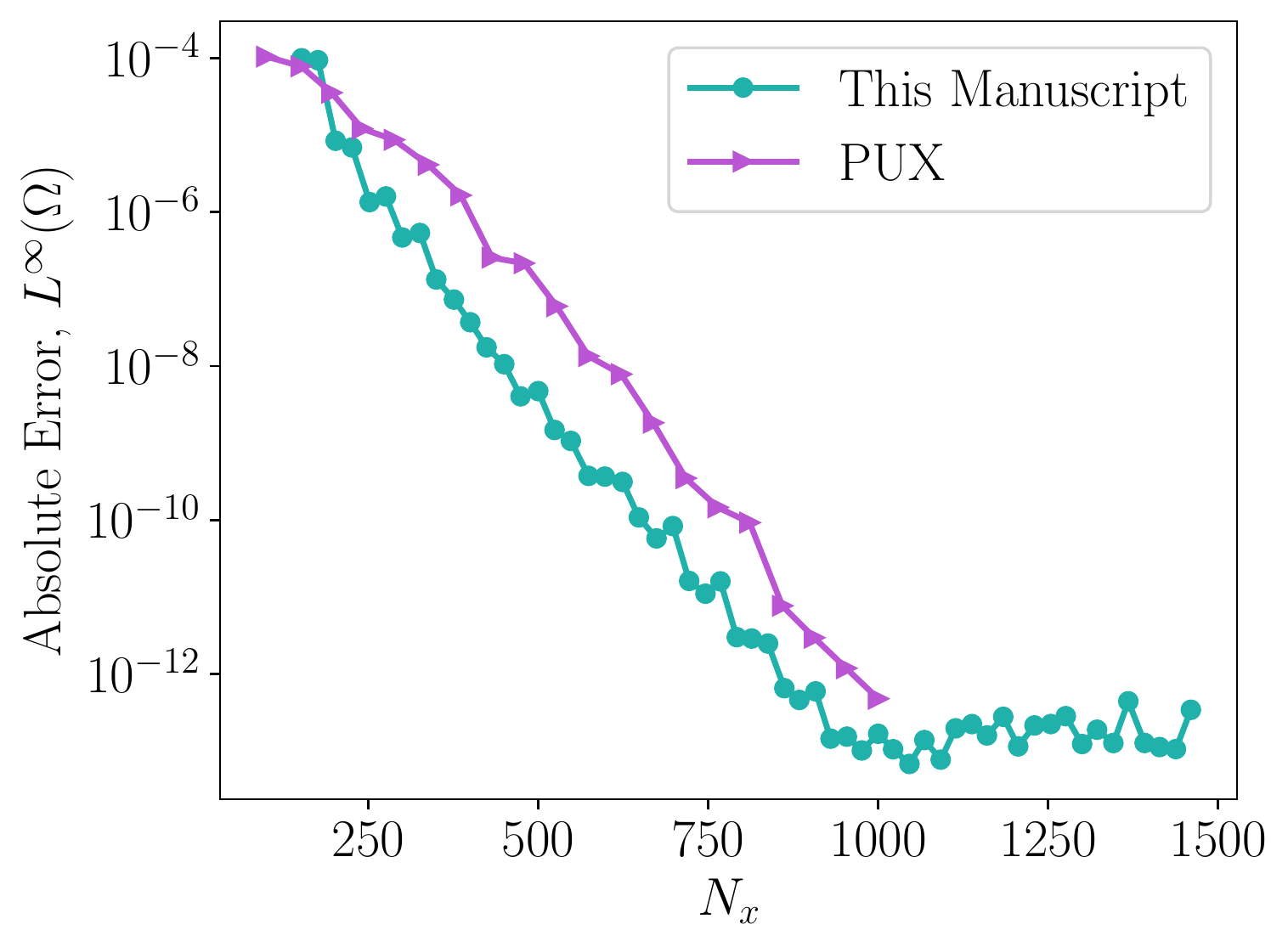}
    \vspace{-0.7em}
    \caption{Comparison with PUX}
  \end{subfigure}
  \hfill
  \begin{subfigure}[c]{0.3\textwidth}
    \centering
    \includegraphics[width=\textwidth]{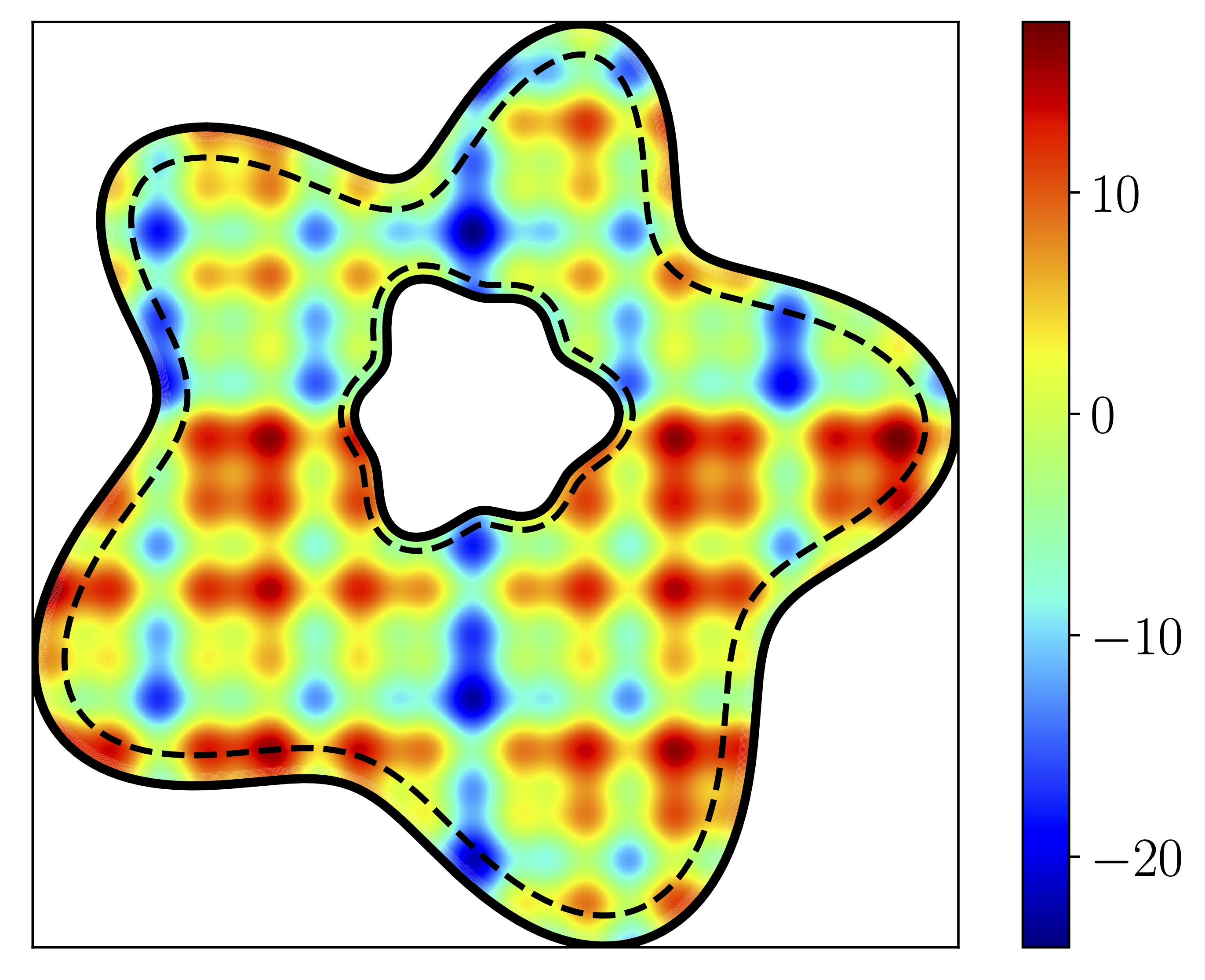}
    \vspace{0.2em}
    \caption{RHS $f$}
  \end{subfigure}
  \hfill
  \begin{subfigure}[c]{0.3\textwidth}
    \centering
    \includegraphics[width=\textwidth]{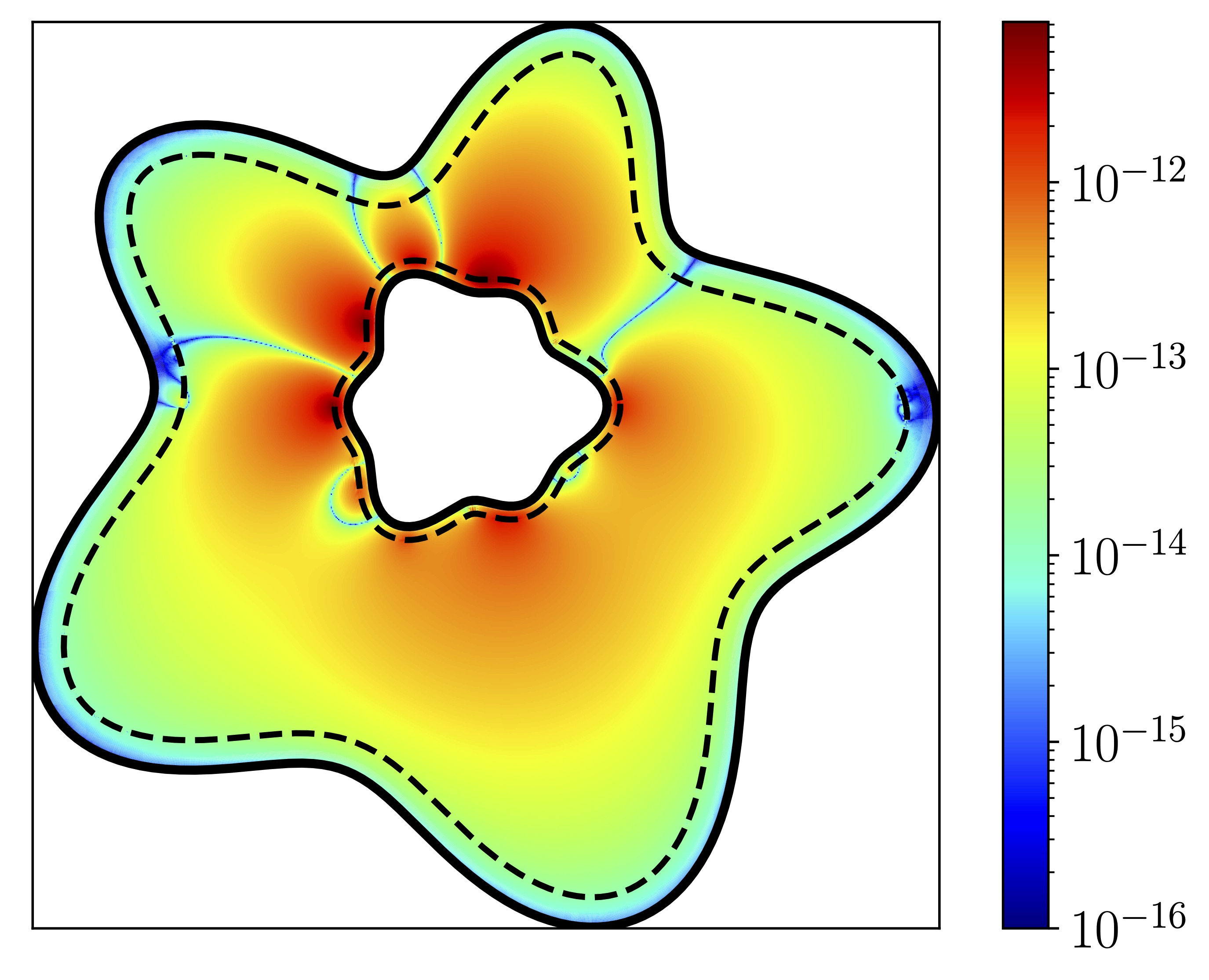}
    \vspace{0.2em}
    \caption{$L^\infty(\Omega)$ Error, $N_x=800$}
  \end{subfigure}
  \caption{Results for multi-boundary Poisson problem with high frequency RHS $f$ from \cite{fryklund2018partition}. Panel (a) shows Absolute $L^\infty(\Omega)$ errors for both our method and PUX. Panel (b) shows the right-hand side $f$ and the domain; Panel (c) shows absolute errors for $N_x=800$.}
  \label{figure:pux_comp}
\end{figure}

\section{Comparison to Fourier-Continuation method}
\label{section:fourier_continuation}

In this section we compare to the 2D Fourier-continuation method, as presented in \cite{bruno2020two}. In particular, we solve the Poisson problem with right-hand side $f=-\sin(kx)\sin(ky)$ set on the ``kite-shaped'' domain given by $\X(s)=(\cos(\theta)+0.35\cos(2\theta)-0.35)\hat{\mathbf{x}} + 0.7\sin(\theta)\hat{\mathbf{y}}$. To enable direct comparison, we solve for a smooth right-hand side with $k=2\pi$ (see Example 4.2 and Table 3 in \cite{bruno2020two}), and a highly oscillatory right-hand side, with $k=40\pi$ (see Example 4.3 and Table 4 in \cite{bruno2020two}), and report the most comparable error diagnostics possible. 

\begin{figure}[h!]
  \centering
  \begin{subfigure}[c]{0.3\textwidth}
    \centering
    \includegraphics[width=\textwidth]{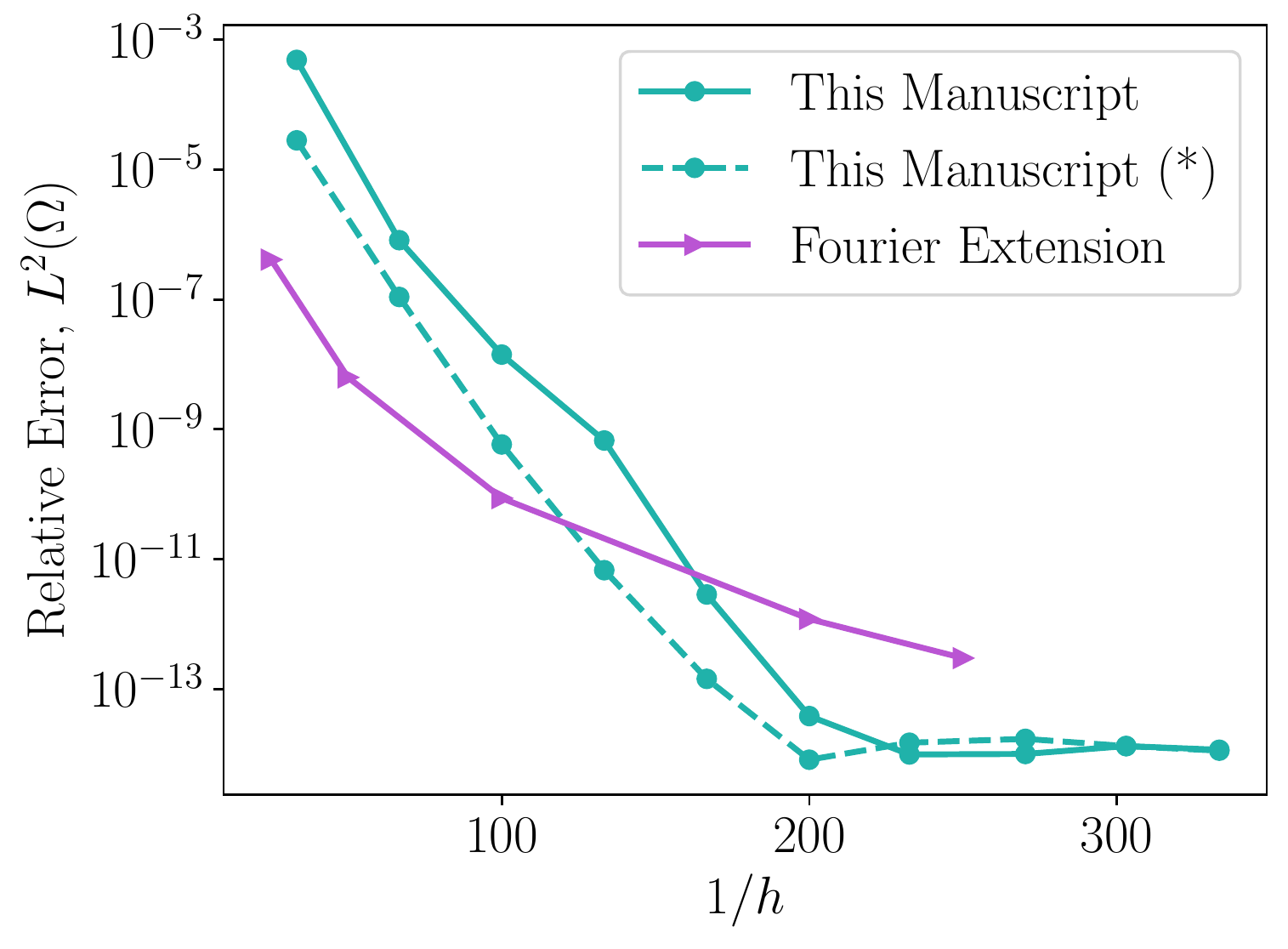}
    \caption{$k=2\pi$}
  \end{subfigure}
  \hfill
  \begin{subfigure}[c]{0.3\textwidth}
    \centering
    \includegraphics[width=\textwidth]{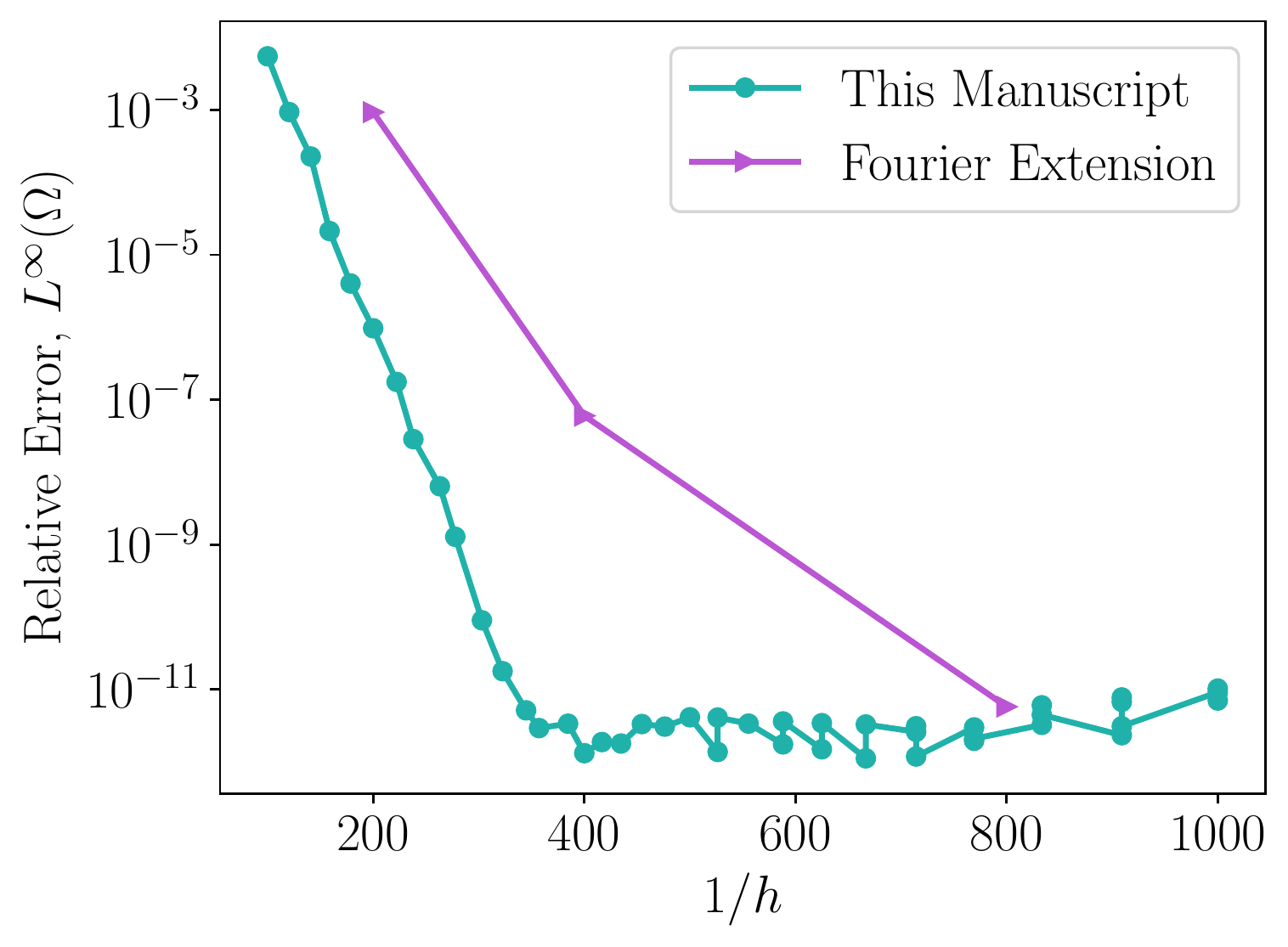}
    \caption{$k=40\pi$}
  \end{subfigure}
  \hfill
  \begin{subfigure}[c]{0.37\textwidth}
    \centering
    \includegraphics[width=\textwidth]{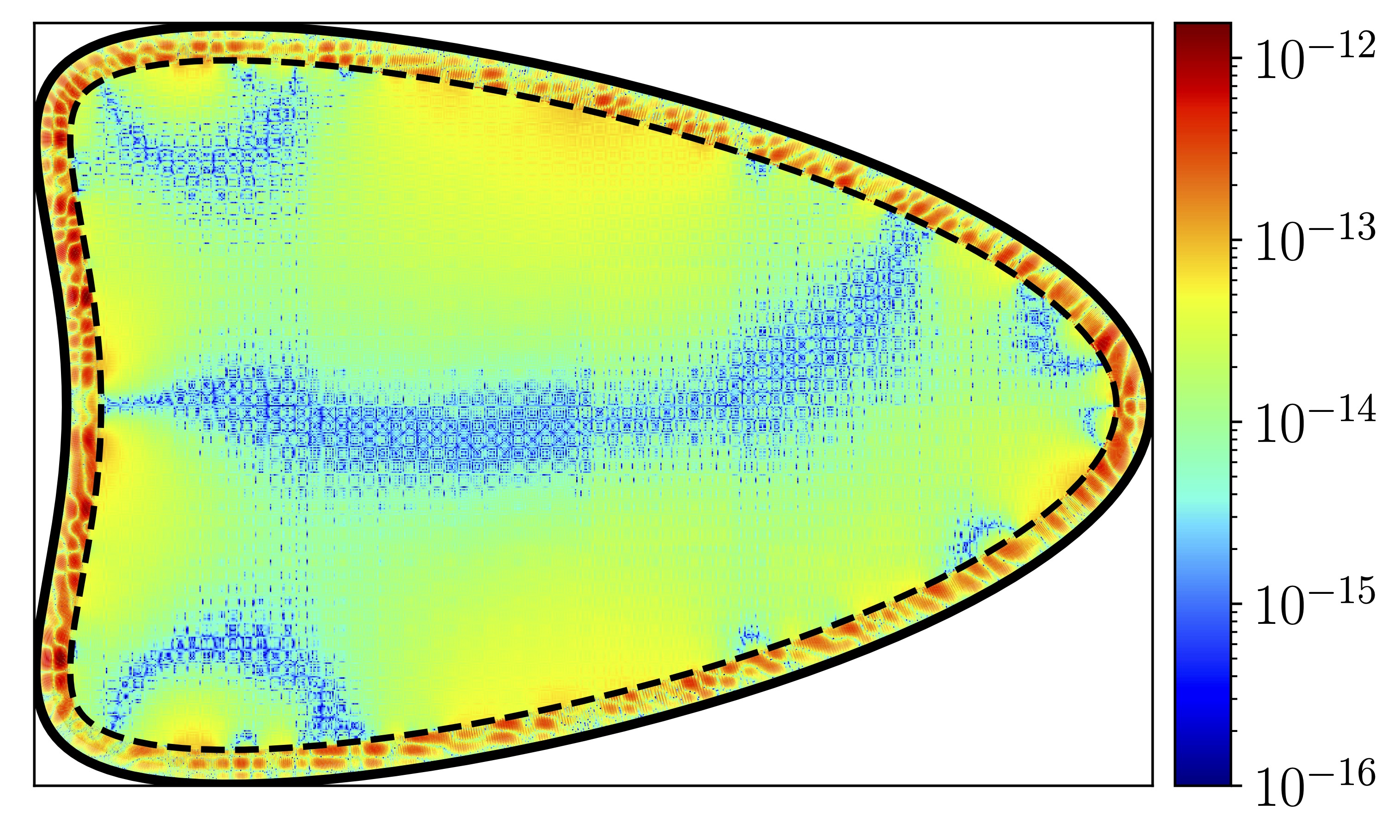}
    \caption{Error, $k=40\pi$ with $h=0.0025$.}
  \end{subfigure}
  \caption{Comparison to Fourier continuation method from \cite{bruno2020two}. In both problems $f=-\sin(kx)\sin(ky)$. Panel (a) shows relative $L^2(\Omega)$ errors for $k=2\pi$. The solid green line shows our errors from our solution with parameters set as described in \Cref{section:parameter_selection}; the dashed line (indicated with an asterisk), shows solutions with a more aggresive radial region set by $R=0.75R_\textnormal{max}$. Panel (b) shows relative $L^\infty(\Omega)$ errors for the harder problem with $k=40\pi$; note as well that our solutions are $L^\infty(\Omega)$ over the \emph{whole domain}, while those reported in \cite{bruno2020two} are reported only far from the boundary $\Gamma$. Panel (c) shows the domain for this problem, along with the error computed by our scheme when $h=0.0025$.}
  \label{figure:fe_comparison_easy_hard}
\end{figure}

In \Cref{figure:fe_comparison_easy_hard}, we show refinement studies for both values of $k$, with results from this manuscript and those reported in \cite{bruno2020two}, along with a plot of the domain and pointwise error when $k=40\pi$ and $h=0.0025$. When $k=2\pi$, the right-hand side is smooth. Errors (shown in panel a) from our method converge spectrally, achieving near-machine precision when $h\approx0.005.$ At large $h$, our method produces larger errors, with crossover achieved for $h$ slightly less than $0.01$. The reason for larger errors at large $h$ is simple: the domain here, shown in panel (c), has high-curvature regions where it is \emph{convex}, and low-curvature regions where it is \emph{concave}. For an interior problem, function \emph{intension} is hard in high-curvature convex regions: normal coordinates moving into the domain cross quickly, forcing $R_\textnormal{max}$ to be small. The reverse is true of function extension, which will have more difficulty when high-curvature regions occur in concave regions for interior problems. Panel (b) shows a similar refinement study, but now with $k=40\pi$. In this case, the right-hand side is more oscillatory, and by the time it is well resolved $M$ is large despite the fact that $R_\textnormal{max}$ is small; and so the higher-order convergence achieved by our scheme is sufficient to provide lower errors at all values of $h$. Panel (c) shows the domain, along with the associated errors when $k=40\pi$ and $h=0.0025$.

\begin{figure}[h!]
  \centering
  \begin{subfigure}[c]{0.5\textwidth}
    \centering
    \includegraphics[width=\textwidth]{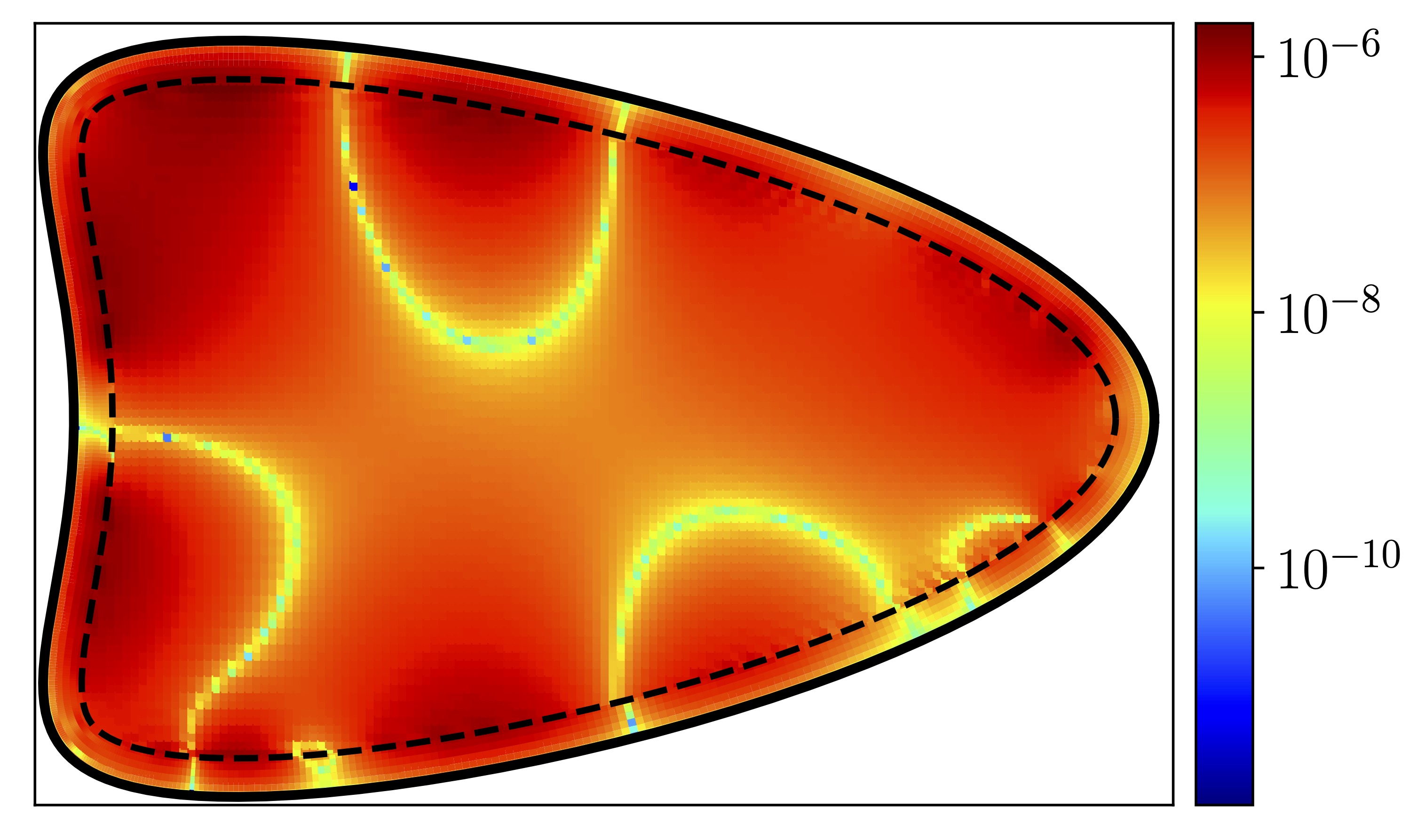}
    \caption{}
  \end{subfigure}
  \hspace{3em}
  \begin{subfigure}[c]{0.35\textwidth}
    \centering
    \includegraphics[width=\textwidth]{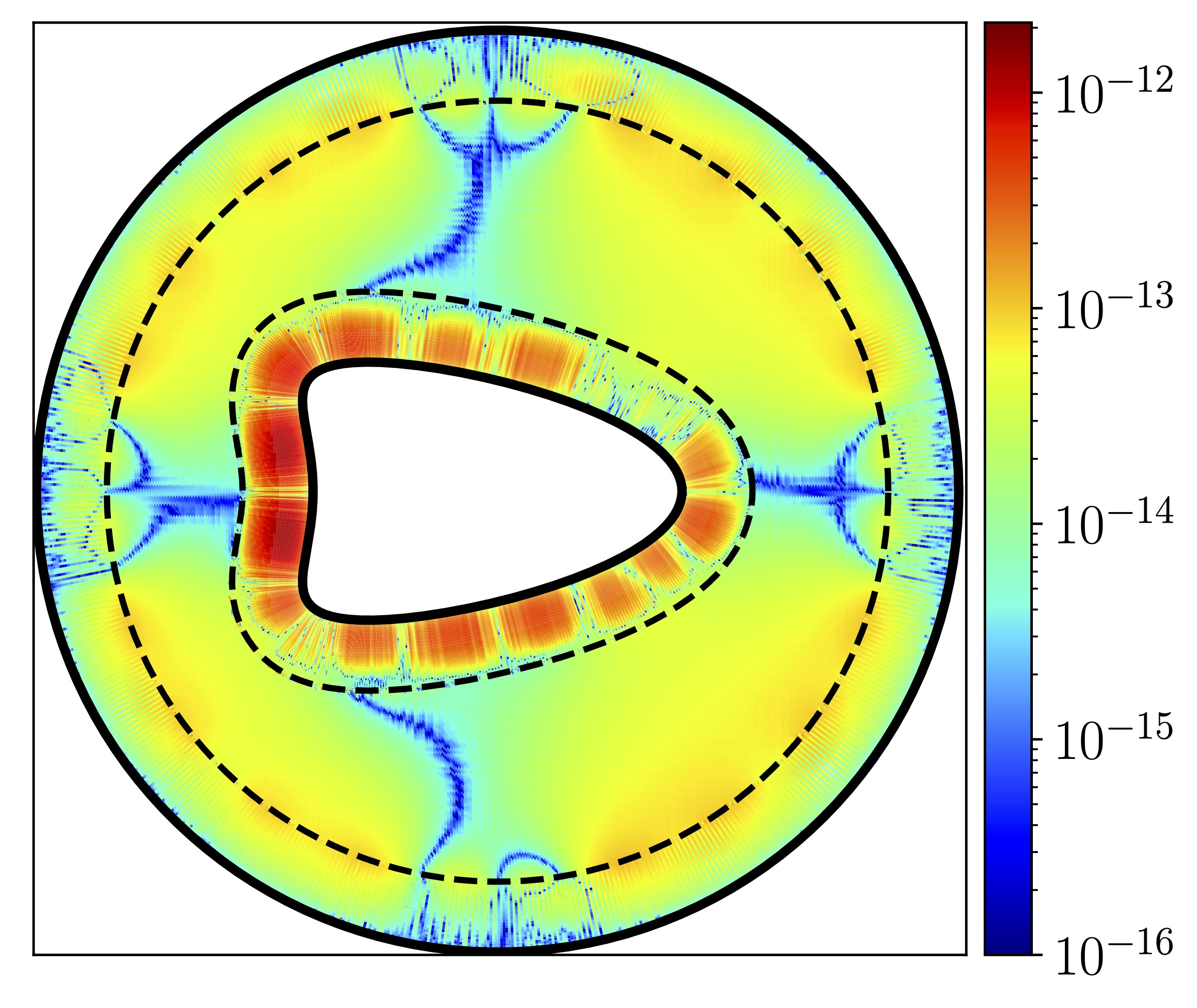}
    \caption{}
  \end{subfigure}\\
  \caption{Domains and errors for the interior (a) and exterior (b) problem set on the kite-shaped domain with smooth right-hand side from \cite{bruno2020two}. The kite-shaped domain here is the same physical size in each image, and $h=0.015$ is the same. Because the domain has high curvature regions where it is convex, $M$ is forced to be much smaller for the interior case ($M=6$) than the exterior case ($M=40$). When the right-hand side is much smoother than the geometry, as in this case, the difference in errors can be large.}
  \label{figure:fe_interior_exterior}
\end{figure}

We use this example as a way to further analyze domains on which function intension vs. extension will have an easier time. For interior problems, intension requires more resolution when there are high-curvature convex regions. The reverse is true for exterior problems: high-curvature convex regions present little problem, but high-curvature concave regions are challenging. To get a handle on this, we compare the solutions generated on this domain, for both interior and exterior problems, across a range of values of $k$. We first show results for both an interior and exterior problem, with $h=0.015$ and $k=2\pi$, in \Cref{figure:fe_interior_exterior}(a-b). For the exterior problem, an outer confining circle is added. As expected, for this smooth, low-frequency problem, errors are far lower --- by about 6 orders of magnitude at this value of $h$ --- for the exterior problem than the interior problem. In the exterior problem, $R_\textnormal{max}$ is large, and so $M$ can be big --- here $M=40$. For the interior problem, $M$ is instead $6$. Errors are dominated by resolving the geometry, rather than the function.

In \Cref{figure:fe_interior_exterior}(c), we show relative $L^2(\Omega)$ errors for both interior problems (solid lines) and exterior problems (dashed lines) across a range of values of $k$, with the most purple line corresponding to $k=2\pi$ and the most red line corresponding to $k=38\pi$, incrementing by $4\pi$. The black line shows the results from the FC method for $k=2\pi$ (reported only for the interior problem). As expected, errors produced by function intension on the exterior problem are smaller at all discretizations than those produced for the interior problem. To better analyze the error in these, we compute the ratio between these two errors, shown in Panel (d). These differences are most extreme at small $k$. In this case, achievable errors are severely limited in the interior problem by our ability to resolve the geometry, with the ratio maximized and very large (over 7 digits) at small $h$. As $k$ is increased, the length-scale of the RHS becomes more comparable to, and eventually smaller than, the length-scale associated with the boundary curvature, and the ratio is reduced to $\approx 3$ digits.

\begin{figure}[h!]
  \begin{subfigure}[c]{0.45\textwidth}
    \centering
    \includegraphics[width=\textwidth]{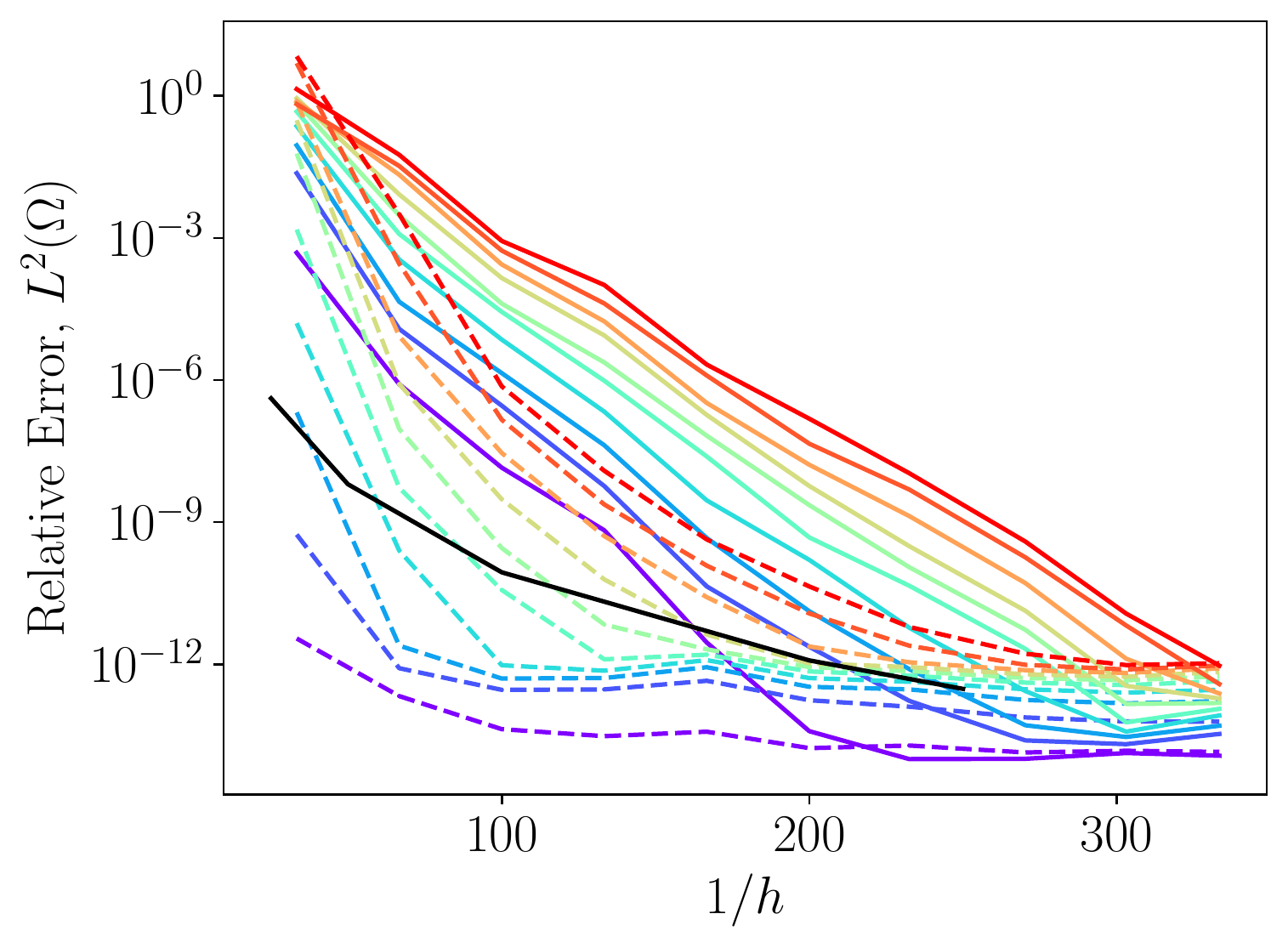}
  \end{subfigure}
  \hspace{3em}
  \begin{subfigure}[c]{0.45\textwidth}
    \centering
    \includegraphics[width=\textwidth]{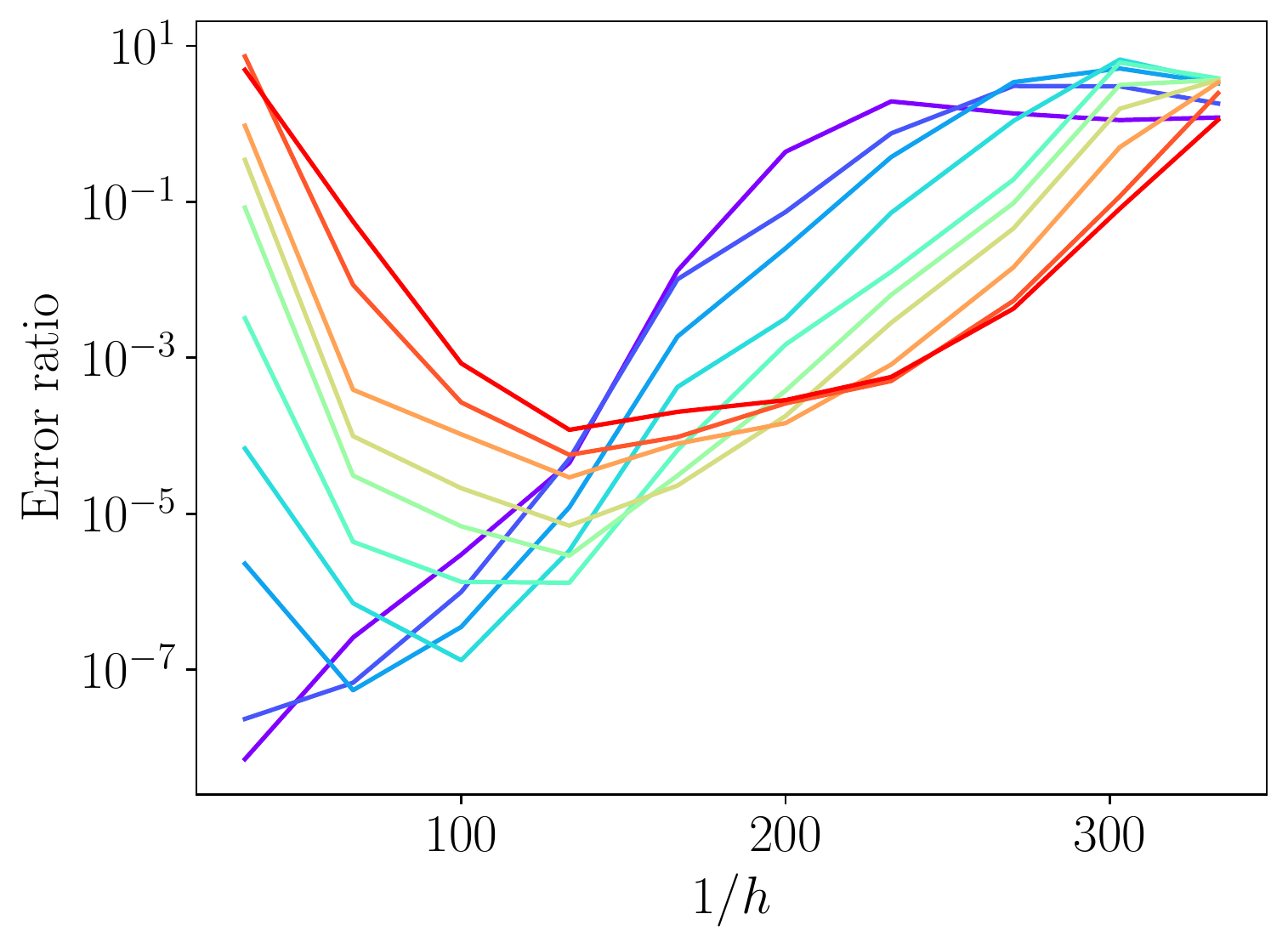}
  \end{subfigure}
  \label{figure:fe_interior_exterior_refinement}
  \caption{Panel (a) Shows relative $L^2$ errors on both the interior and exterior domains (see \Cref{figure:fe_interior_exterior}) over a range of $h$ values. Errors associated with the interior domain are shown as solid lines; errors associated with the exterior domain are shown as dashed lines. The results for the same interior problem produced by the Fourier continuation method, are shown as a black line \cite{bruno2020two}. Panel (b) shows the ratio between the errors from the exterior and interior solutions (exterior error / interior error).}
\end{figure}

\section{The modified-Helmholtz problem}
\label{section:modified_helmholtz}

\begin{figure}[h!]
  \centering
  \begin{subfigure}[c]{0.4\textwidth}
    \centering
    \includegraphics[width=\textwidth]{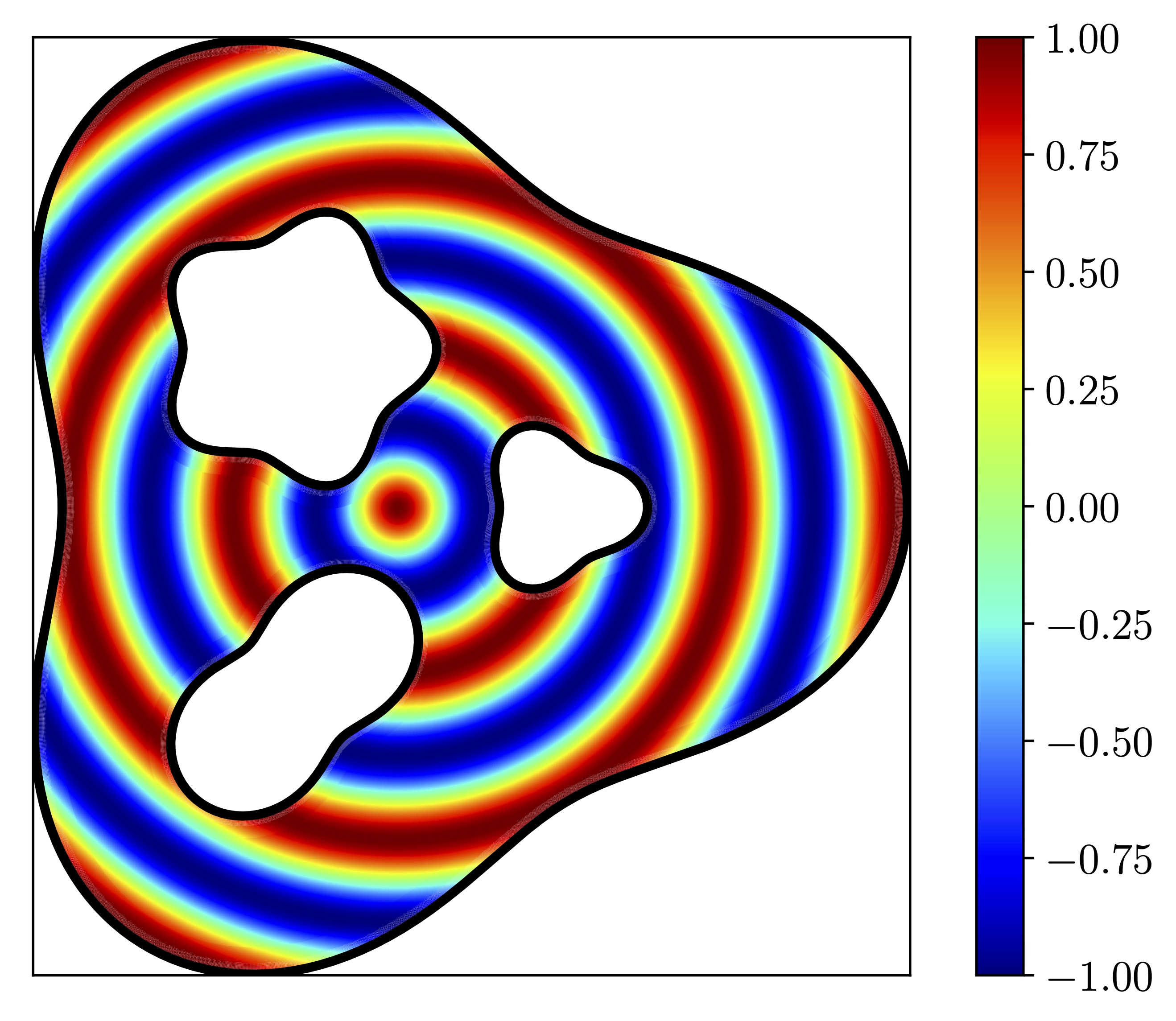}
  \end{subfigure}
  \hspace{5em}
  \begin{subfigure}[c]{0.4\textwidth}
    \centering
    \includegraphics[width=\textwidth]{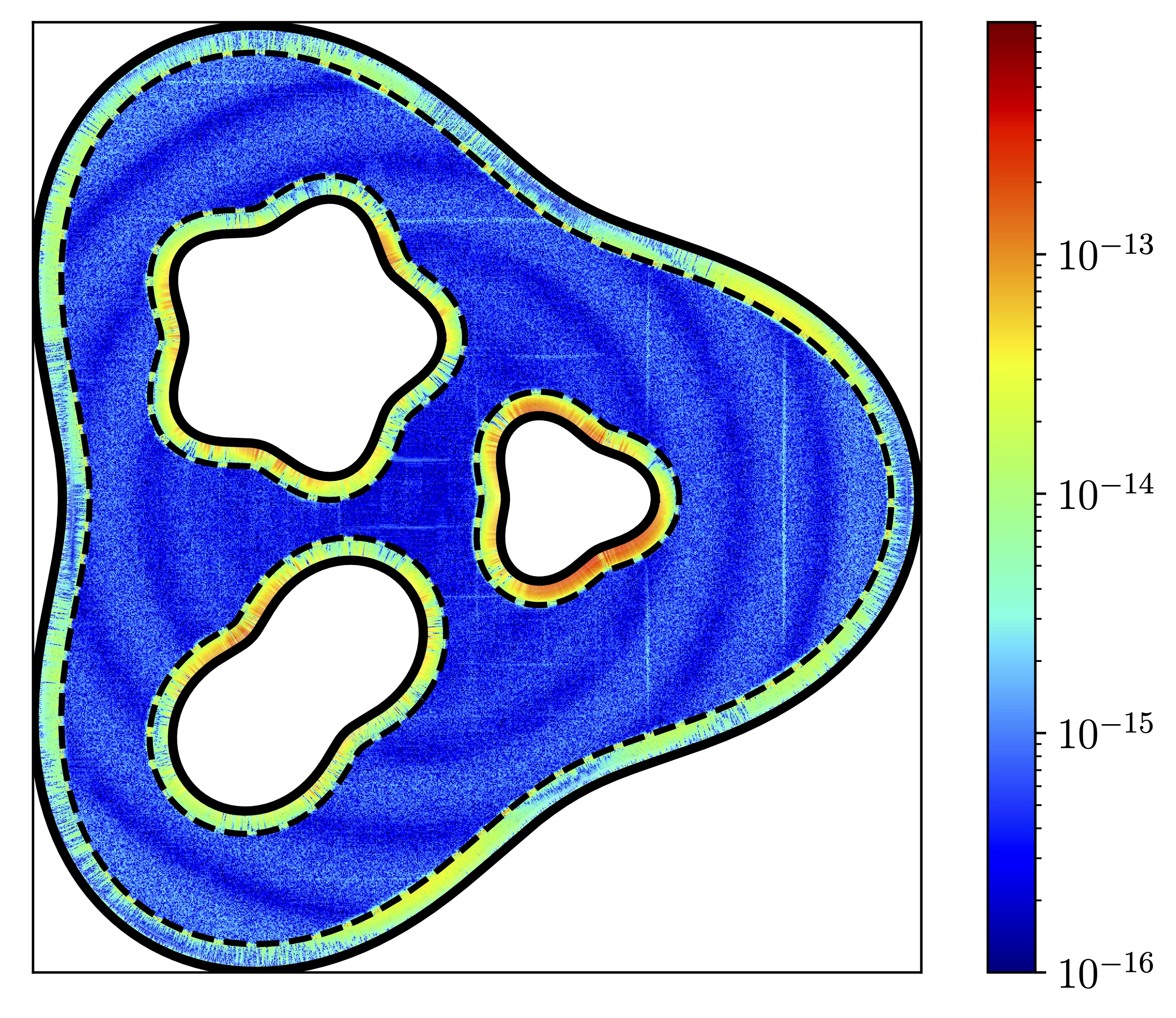}
  \end{subfigure}
  \label{figure:mh_solutions}
  \caption{Domain, solution, and errors to the modified-Helmholtz problem studied in \Cref{section:modified_helmholtz}. Panel (a) shows the solution $u$; panel (b) shows the error when $h=0.002$ and $\alpha^2=10^5$. The $L^\infty(\Omega)$ error in this case is $<10^{-12}$.}
\end{figure}

We now turn our attention to solving the inhomogeneous modified-Helmholtz problem. This problem arises when discrezing the heat equation in time: consider the simplest Forward-Euler/Backward-Euler IMEX scheme for the discretization of $u_t-\nu\Delta u=f$, which gives:
\begin{equation}
  (\mathbb{I}-\nu\Delta t\Delta)u(t+\Delta t) = u(t) + \Delta t f(t),
\end{equation}
subject to appropriate boundary conditions (if $f$ is analytically known, it can be taken at $t+\Delta t$, but typically $f$ is generated by non-linear terms, e.g. reactions or advection). Dividing through by $\nu\Delta t$ gives:
\begin{equation}
  (\alpha^2\mathbb{I}-\Delta)u(t+\Delta t) = \alpha^2 u(t) + \alpha^2\Delta t f(t),
\end{equation}
with $\alpha=1/\sqrt{\nu\Delta t}$. Clearly if the diffusion coefficient $\nu$, the timestep $\Delta t$, or both are small, $\alpha$ can be (very) large. We thus seek to solve the problem 
\begin{subequations}
  \begin{align}
    (\alpha^2\mathbb{I}-\Delta) u &= f &&\textnormal{in }\Omega,  \\
    u &= g &&\textnormal{on }\Gamma,
  \end{align}
  \label{equation:modified_helmholtz}
\end{subequations}
typically for $\alpha\gg1$. Only minor modifications need to be made to the algorithm as presented so far. We collect these here:
\begin{enumerate}
  \item Function intension is unchanged, and the regular problem is changed only in that the Fourier symbol of the differential operator is now $\alpha^2+\k^2$. This operator, unlike the Poisson operator, is invertible, with no solvability condition, and so the uniform grid can be taken tight to the boundary.
  \item The method for the annular problem is nearly unchanged. \Cref{equation:annular_laplacian_k} becomes instead:
  \begin{equation}
    \mathcal{F}[(\alpha^2\mathbb{I}-\Delta)\ua]_k = -\frac{1}{\psi}\left[\frac{\partial}{\partial r}\left(\psi\frac{\partial\widehat{\ua}_k}{\partial r}\right) + \left(\frac{1}{\psi}k^2 - \psi\alpha^2\right)\widehat{\ua}_k\right],
  \end{equation}
  which is again inverted using GMRES preconditioned by the (separable) inverse for a circular annulus.
  \item The radially symmetric Green's function is now $G_\alpha(r)=\alpha^2K_0(\alpha r)/(2\pi)$. Because $K_0(x)=-\log x + \xi(x)$ where $\xi$ is smooth, jump conditions for the single and double layer potentials are, up to constants, the same as those for $\log r$. The stitching step is thus the same, with jumps in the value corrected by double-layer potentials and jumps in the normal derivative corrected by single-layer potentials.
  \item The evaluation of these layer potentials is again done using the QFS-B aglorithm from \cite{stein2021quadrature}, with the singular on-surface evaluation done using 16th-order Alpert quadrature \cite{hao2014high}. When $\alpha$ is large, if the definition of the source curve used in \cite{stein2021quadrature} is used, $G_\alpha(r)$ decays too rapidly to convey information between the source and check curves. To avoid this, we upsample the source curve, and move it towards the boundary, by a factor of $\max(1, h\alpha)$. Note that this upsampling factor is only large when $h\gg1/\alpha$, in which case the problem is poorly resolved.
  \item The homogeneous correction problem is unchanged, and its solution is computed by solving a well-conditioned second-kind BIE; singular integral operators are computed using 16th-order Alpert quadrature \cite{hao2014high}.
\end{enumerate}
To analyze the convergence and behavior as a function of $\alpha$, we solve the Dirichlet problem generated by the solution $u=\cos\left(20\sqrt{x^2+y^2}\right)$, set on a multiply connected domain similar (though not identical) to the domain on which this same problem was analyzed in \cite{fryklund2020integral}, with solutions generated by a modified version of the PUX method, to which we compared our Poisson solver in \Cref{section:PUX_comparison}. The solution function $u$ and domain are shown in \Cref{figure:mh_solutions}, along with the pointwise error when $h=0.002$ and $\alpha^2=10^5$. Even for this relatively large $\alpha$, we are able to achieve a full 12 digits of accuracy, with errors in the annular solve evidently dominant.

\begin{figure}[h!]
  \centering
  \begin{subfigure}[c]{0.7\textwidth}
    \centering
    \includegraphics[width=\textwidth]{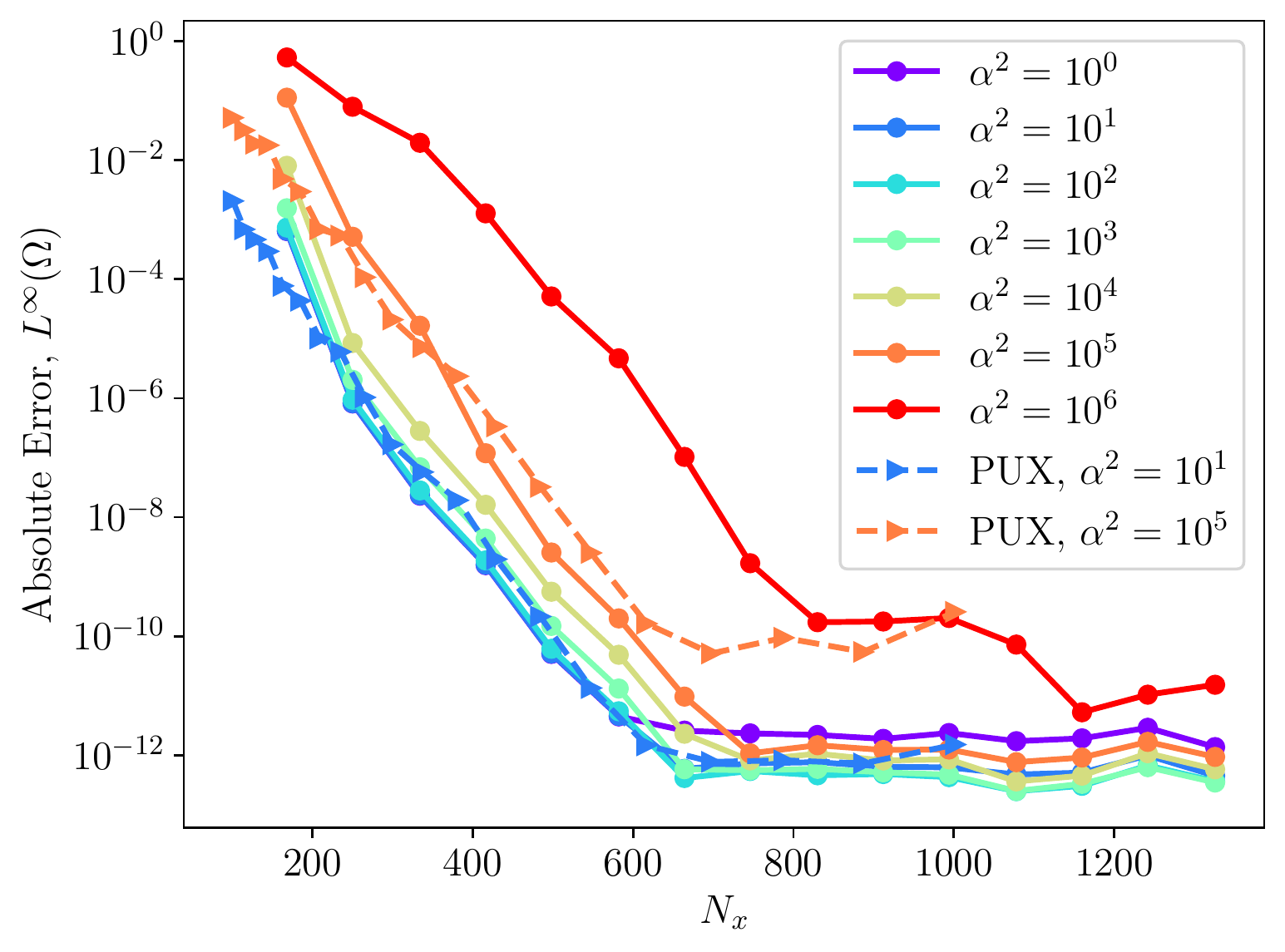}
  \end{subfigure}
  \label{figure:mh_refinement}
  \caption{Refinement study for the modified Helmholtz problem studied in \Cref{section:modified_helmholtz}. Solid lines are from the method introduced in this manuscript, with the indicated value of $\alpha$. Dashed lines are generated by the PUX method, extracted from Figure 8 in \cite{fryklund2020integral,Rohatgi2020}, again with the indicated value of $\alpha$. Both methods show exponential convergence with the same rate.}
\end{figure}

We now solve over a range of values of $h\in[0.00125, 0.01]$ and $\alpha^2=10^0,\ 10^1,\ \ldots,\ 10^6$, with results shown in \Cref{figure:mh_refinement}, although we have plotted the errors against the number of points discretizing our uniform grid in the $x$-direction, to enable comparison to the PUX method. For $\alpha^2$, up to $10^5$, we observe fast and stable convergence to $\approx10^{-12}$, with the convergence curve shifted up slightly for higher values of $\alpha$. For $\alpha^2=10^6$, our method requires markedly finer discretizations, and fails to reliably produce more than 10 digits. The slower convergence is perhaps not surprising: the length scale associated with the modified-Helmholtz equation in this case is $0.001$; the finest discretization we test is $h=0.00125$, and so it is perhaps remarkable that we achieve near-machine precision solutions when a length-scale in the problem is under-resolved; this may simply be due to the choice of a test problem whose solution does not have this length-scale present. The FMM library we were utilizing failed for higher values of $\alpha$ than this. Our results are qualitatively similar to those observed by the PUX method \cite{fryklund2020integral}. For comparison, we also plot the errors associated with their method, for $\alpha^2=10^1$ and $\alpha^2=10^5$. Because our domains are not exactly the same, we caution that comparing the errors directly is unwarranted\footnote{Indeed, it is likely that PUX is producing slightly better errors for a given $h$, as their domains necessarily include room for an extension.}; however, both methods, perhaps surprisingly, produce nearly exactly the same rate of exponential convergence. Just as with our method, errors are slightly worse when $\alpha$ is large. The only minor difference between the results is that our method appears to be slightly more stable for large $N$ at high $\alpha$, finding about two more good digits in the solution when $\alpha=10^5$.

\section{Stokes}
\label{section:stokes}
\begin{figure}[h!]
  \centering
  \begin{subfigure}[c]{0.4\textwidth}
    \centering
    \includegraphics[width=\textwidth]{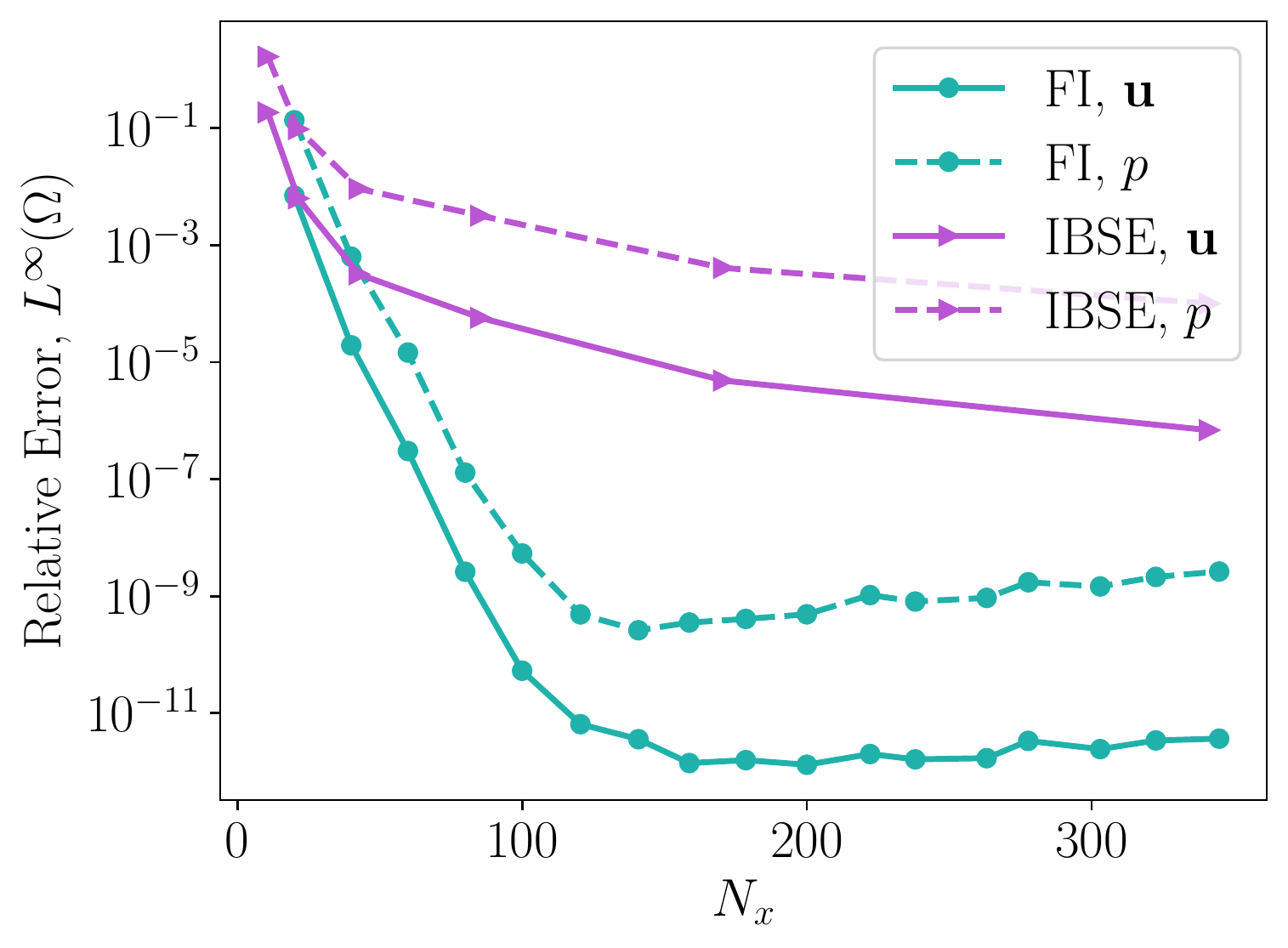}
    \caption{(a)}
  \end{subfigure}
  \hspace{2em}
  \begin{subfigure}[c]{0.4\textwidth}
    \centering
    \includegraphics[width=\textwidth]{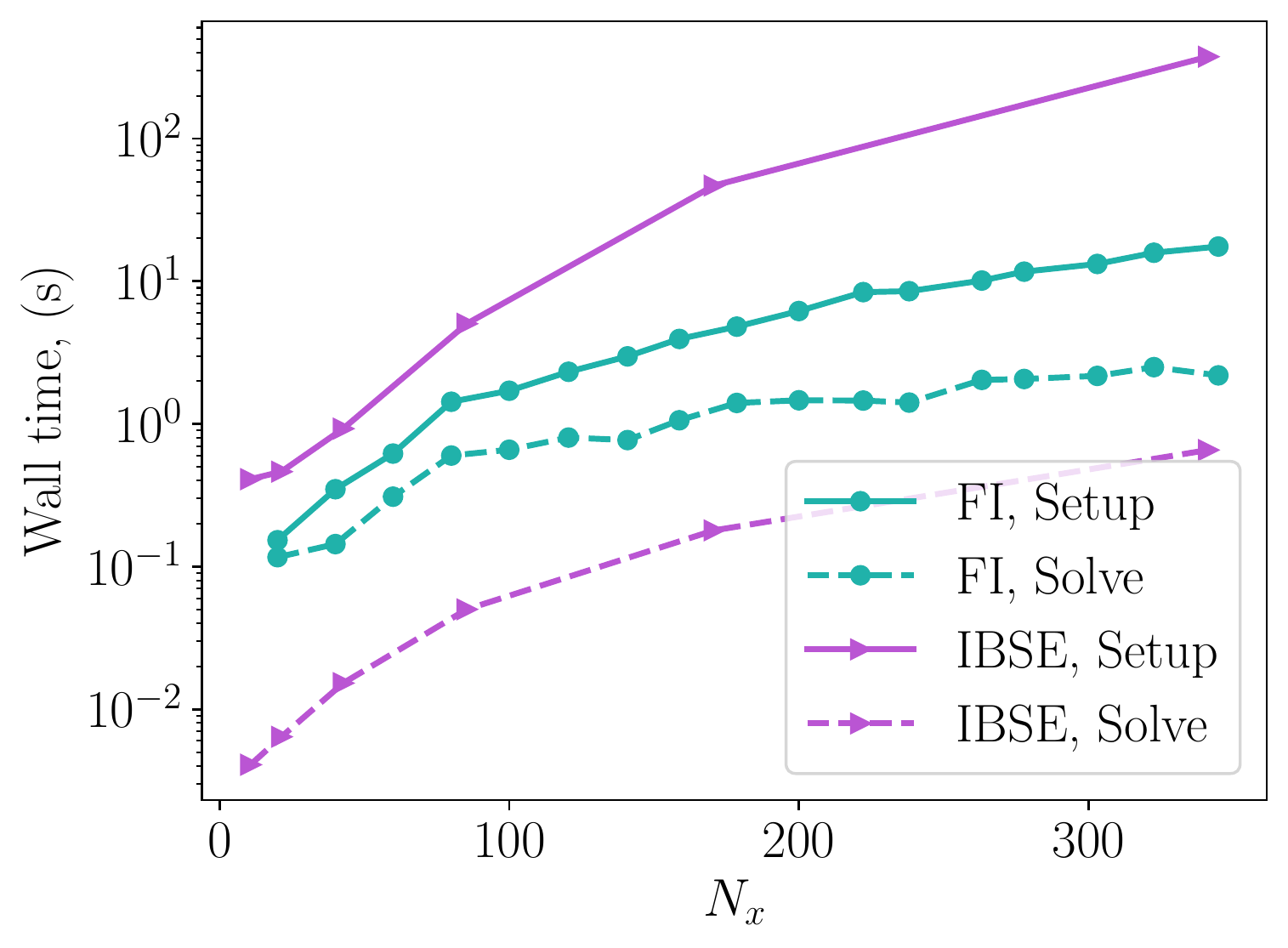}
    \caption{(b)}
  \end{subfigure}
  \caption{Comparison between function intension and the third order IBSE method, for the Stokes problem given in \Cref{section:stokes}. Panel (a) shows relative $L^\infty(\Omega)$ errors for $\u$ and $p$; panel (b) shows wall clock timings, broken down into both setup and solve portions for each respective solver.}
  \label{figure:ibse_comparison}
\end{figure}

Finally, we solve a Stokes problem with Dirichlet boundary conditions:
\begin{subequations}
  \begin{align}
    -\Delta\u + \grad p &= \f &&\textnormal{in }\Omega, \\
    \grad\cdot\u &= 0 &&\textnormal{in }\Omega, \\
    \u           &= \mathbf{g} &&\textnormal{on }\Gamma.
  \end{align}
\end{subequations}
As with the modified-Helmholtz problem shown in \Cref{section:modified_helmholtz}, only minor modifications to the method must be made:
\begin{enumerate}
	\item Both components of the force $\f$ must be independently adjusted to have $0$ mean on $\mathcal{C}$ (see \Cref{section:methods_preliminaries:peridoic_compatibility}).
	\item The method for the annular problem is somewhat more complicated. In our local coordinates system, the Laplacian of a vector field $\A$ is given by \cite{hirota1982vector}:
	\begin{subequations}
	    \begin{align}
	        (\Delta \A)_r &= \Delta A_r - \frac{2}{\psi^2}\frac{\partial\psi}{\partial r}\frac{\partial A_t}{\partial t} - \frac{A_r}{\psi^2}\left(\frac{\partial\psi}{\partial r}\right)^2 - \frac{A_t}{\psi}\frac{\partial}{\partial r}\left(\frac{1}{\psi}\frac{\partial\psi}{\partial t}\right),  \\
	        (\Delta\A)_t &= \Delta A_t + \frac{2}{\psi^2}\frac{\partial\psi}{\partial r}\frac{\partial A_r}{\partial t} + \frac{A_r}{\psi}\frac{\partial}{\partial t}\left(\frac{1}{\psi}\frac{\partial\psi}{\partial r}\right) - \frac{A_t}{\psi^2}\left(\frac{\partial\psi}{\partial r}\right)^2,
	    \end{align}
	\end{subequations}
	the gradient of a scalar field is:
	\begin{equation}
    \grad p = \frac{\partial p}{\partial r}\hat r + \frac{1}{\psi}\frac{\partial p}{\partial t}\hat t,
\end{equation}
and the divergence of a vector field is given by:
\begin{equation}
    \grad\cdot\u = \frac{1}{\psi}\left[\frac{\partial}{\partial r}\left(\psi u_r\textbf{}\right) + \frac{\partial u_t}{\partial t}\right],
\end{equation}
where in all cases $\Delta$ denotes the scalar Laplacian defined in \Cref{equation:annular_laplacian}. For the case of a circular annulus, the Stokes equations in these coordinates is separable, although the precise formulae are unwieldy, and writing them out offers no further insight. Our algorithm is the same as before: invert the Stokes operator for the real geometry, utilizing the inverse of the Stokes operator on the circular geometry as a preconditioner, with one significant caveat: discretizing all of the Fourier modes and Chebyshev modes leads to a checkerboard type instability in the pressure which prevents robust convergence of the iterative scheme. This is easily remedied by omitting the modes associated with the azimuthal Nyquist frequency.
\item The Green's function is now somewhat more complicated, and while the stitching step is the same at an abstract level, in this case jumps in the value of the velocity are corrected by double-layer potentials and jumps in the traction are corrected by single-layer potentials. Methodology for both the homogeneous correction problem and the close-evaluation of layer potentials can be found in \cite{stein2021quadrature}.
\item The pressure is defined only up to a constant.
\end{enumerate}
We now compare function intension against the 3rd-order Immersed Boundary Smooth Extension solver \cite{stein2017immersed}. The solution and force are manufactured from the solution:
\begin{subequations}
   \begin{align}
     u(x, y) & =e^{\sin{ax}}\cos{by},                    \\
     v(x, y) &= -\frac{a}{b}\cos{ax} e^{\sin{ax}}\sin{by}, \\
     p(x, y) &= \cos{cx} + e^{\sin{dy}}
   \end{align}
 \end{subequations}
 where $\u=u\hat{\mathbf{x}} + v\hat{\mathbf{y}}$.
\begin{figure}[t!]
  \centering
  \begin{subfigure}[c]{0.29\textwidth}
    \centering
    \includegraphics[width=\textwidth]{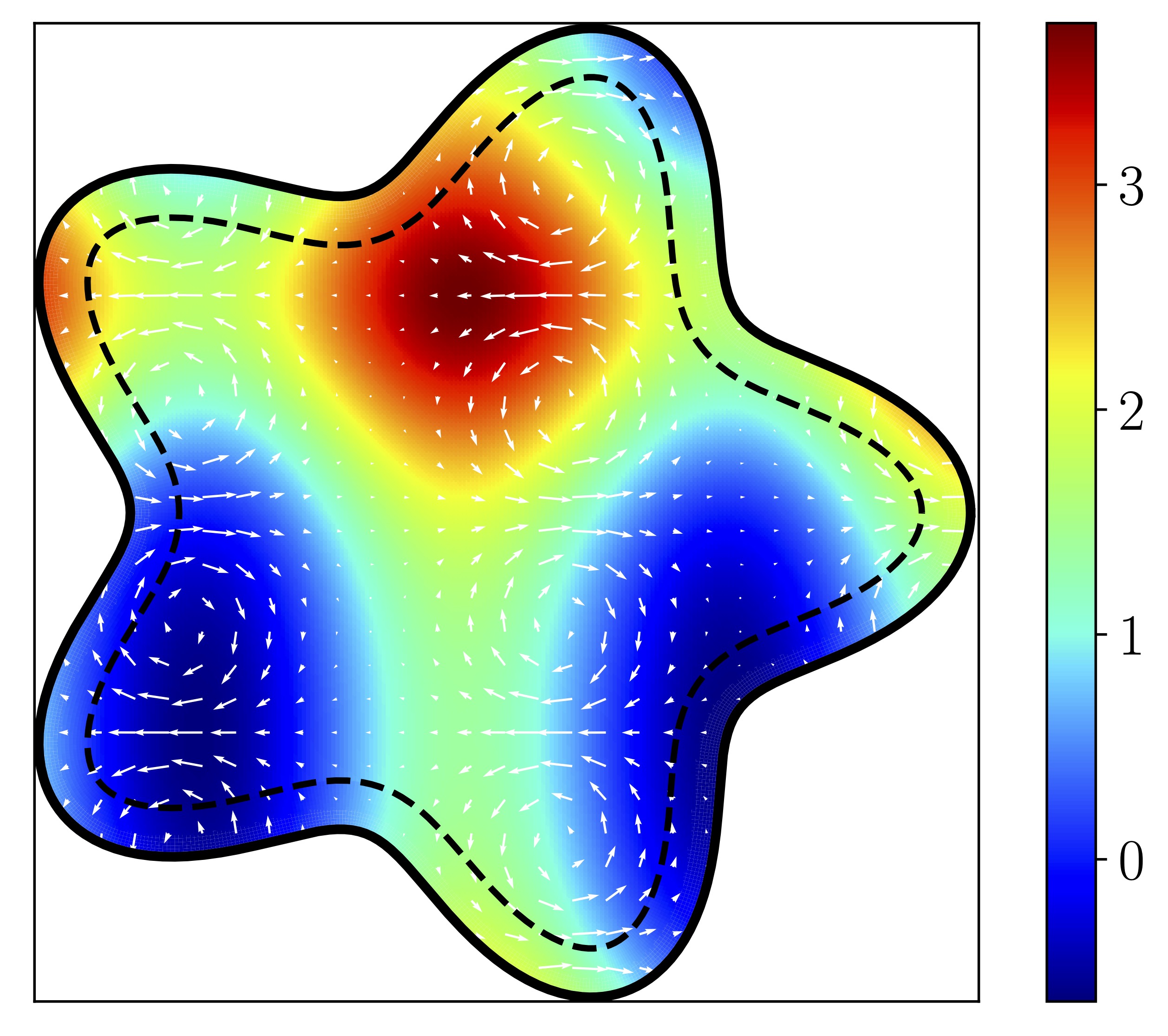}
    \caption{}
  \end{subfigure}
  \hfill
  \begin{subfigure}[c]{0.32\textwidth}
    \centering
    \includegraphics[width=\textwidth]{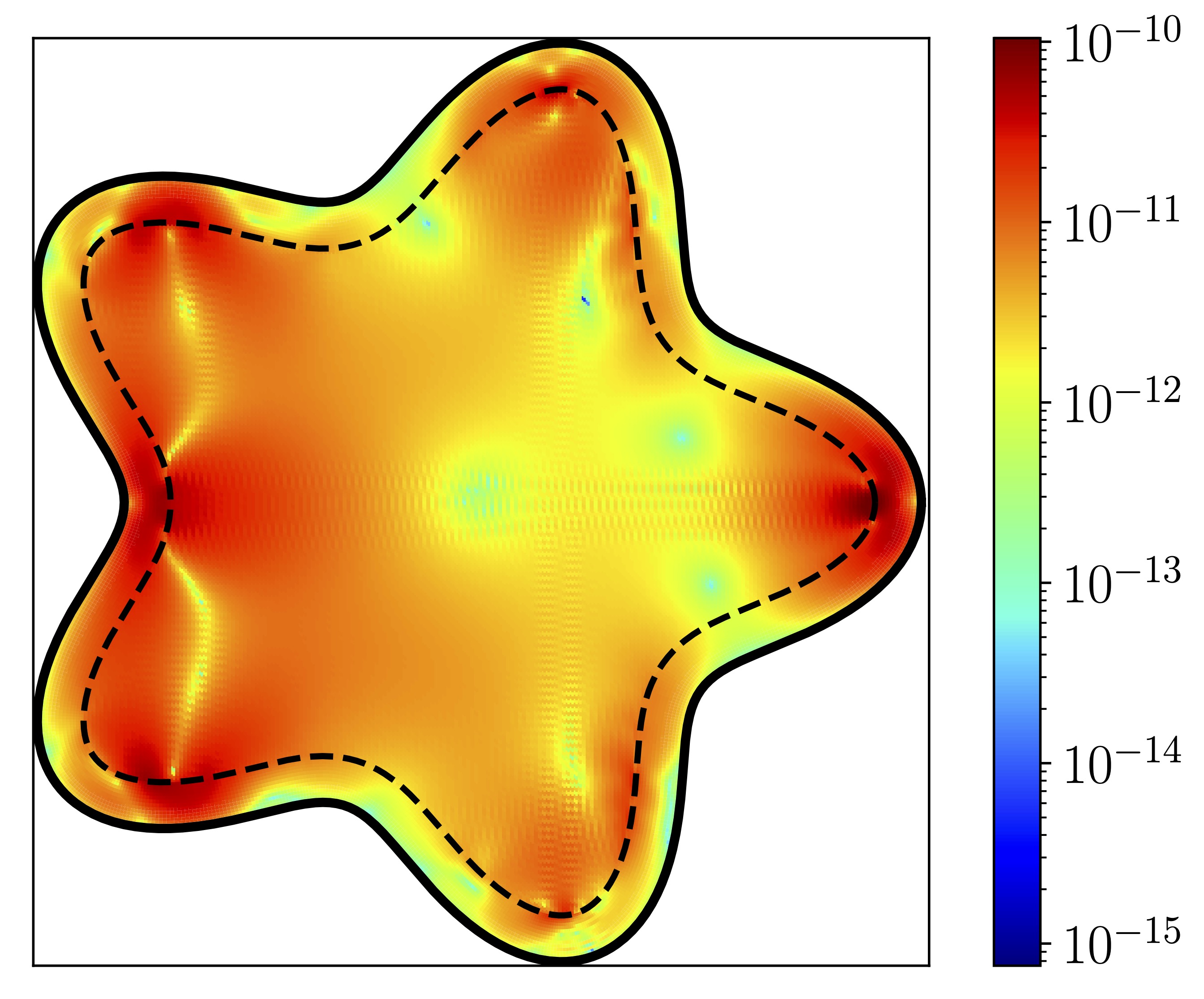}
    \caption{}
  \end{subfigure}
  \hfill
  \begin{subfigure}[c]{0.32\textwidth}
    \centering
    \includegraphics[width=\textwidth]{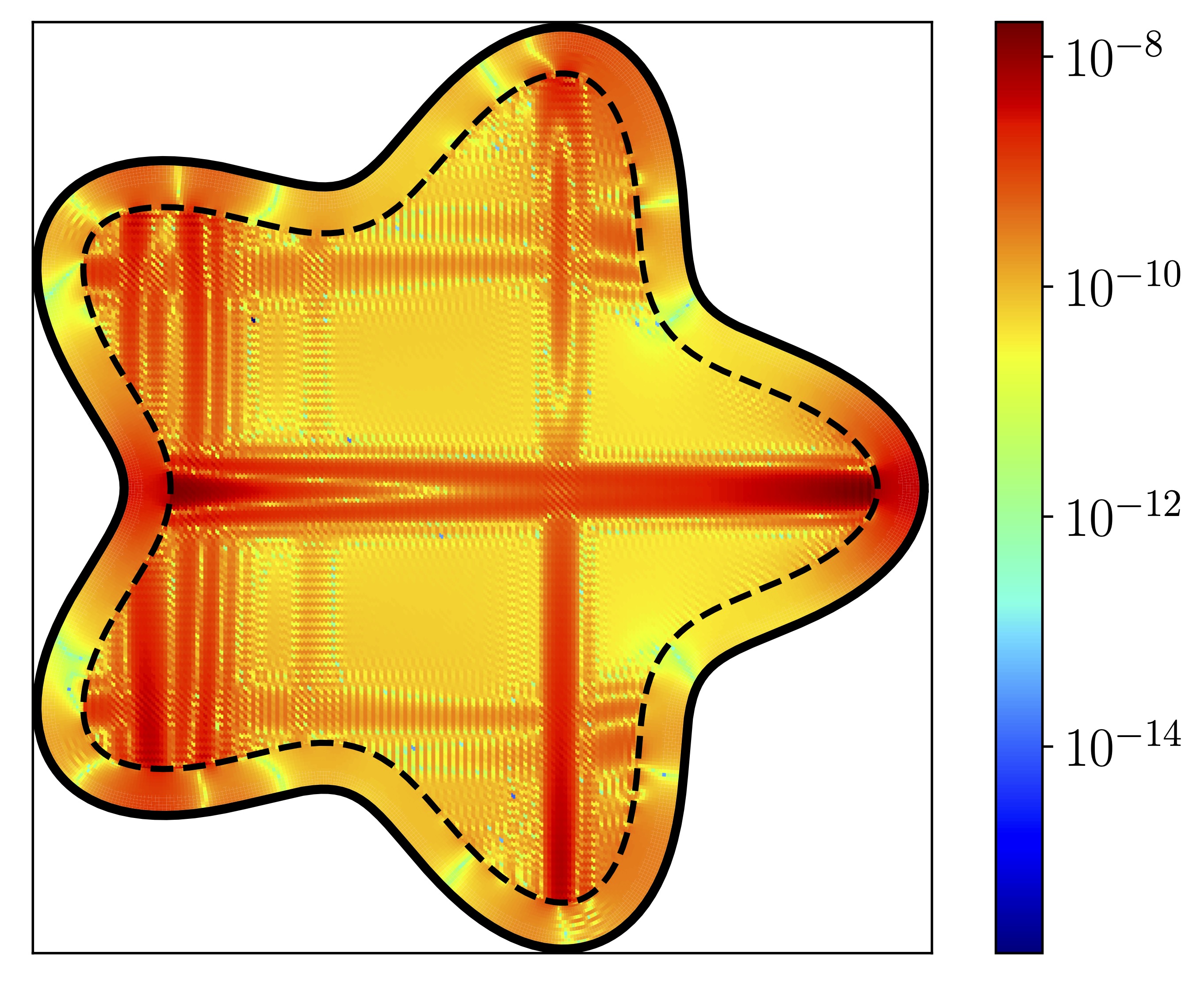}
    \caption{}
  \end{subfigure}
  \caption{Solution and errors for the Stokes problem analyzed in \Cref{section:stokes}. Panel (a) shows pressure (pseudocolor) and velocity field (white arrows) for the analytic solution. Panels (b) and (c) show relative $L^\infty(\Omega)$ errors for $\u$ and $p$, respectively, when $h=0.01$.}
  \label{figure:stokes_solutions}
\end{figure}
For this comparison we take $(a,b,c,d)=(7,6,5,3)$. The solution to this problem, with $\u$ plotted as a vector field overlaying a pseudocolor plot of the pressure field, is shown in \Cref{figure:stokes_solutions}(a), along with errors in $\u$ and $p$ in panels (b) and (c), respectively, when $h=0.01$. We compare solutions and wall-clock timings for this problem for both function intension and the third-order IBSE method\footnote{Results from IBSE generated by personal implementation.}, across a range of values of $h$, in \Cref{figure:ibse_comparison}. Panel (a) shows relative errors for both $\u$ and $p$. At large $h$, errors between the two methods are comparable, though unsurprisingly, convergence is far more rapid for function intension as $h$ is refined, with errors for $\u$ quickly reaching a small multiple of $\epsilon=10^{-12}$. Errors in the pressure function converge at a similar exponential rate, saturating about two digits worse than $\u$. Timings, broken down into ``setup'' and ``solve'' times (with ``setup'' being the re-usable portion of each solve for a fixed domain $\Omega$), are shown in Panel (b); both methods are personal implementations with similar amounts of effort expended on optimization. Across all values of $h$, function intension is considerably faster to setup. Once setup, IBSE produces solutions faster: solves reduce to just a few steps, dominated the FFT and LAPACK calls. Function intension is more complicated, with more computational work that cannot be directly farmed out to highly optimized external routines.  As the discretization is refined the computational load becomes dominated by calls to FMM and NUFFT routines that scale well, and so while still slower, the speed of the FI solves begins to approach that of the IBSE method.

\section{Discussion}
\label{section:discussion}

We have presented a spectrally accurate solver for a relatively wide range of constant-coefficient elliptic PDE, which utilizes the stable process of function \emph{intension} to convert a problem set on a general smooth domain to a problem set on a simple computational domain. In addition to solving a regular grid PDE with the smoothly truncated function acting as the right-hand side, we must additionally solve a PDE in an annulus localized along the boundary, and these solutions are then stitched together using techniques from boundary integral methods, which are also used to impose the physical boundary conditions. Although the analytic scaling to the implementation we use here is asymptotically optimal only after some boundary-dependent setup costs, this choice was made for convenience and can be remedied, as discussed, using established methods, reducing the asymptotic complexity to the same as the FFT used in solving the regular grid problem. Reasonably performant code implementing the method in Python is available in a repository maintained by the author [\footnote{For reviewers: currently located at https://github.com/dbstein/ipde; will be cleaned, better commented, tagged, and archived via Zenodo to accompany final manuscript.}].

To demonstrate the utility and versatility of the method, we solved Poisson, Modified-Helmholtz, and Stokes problems on a variety of domains. While we believe that the method given here has significant utility (and indeed has already found such use in \cite{young2021many}, without a detailed presentation of the numerical method), there are cases where it suffers, the most obvious being when the problem is highly multiscale in nature. There are two separate approaches to improve applicability of the method to multiscale problems. The first is to simply replace the Fourier method used in this manuscript with an adaptive regular grid method. This fits neatly within the paradigm presented herein, necessitating only changes to the regular grid solver and interpolation operators connecting the discretization of $\mathcal{C}$ to the annular grid $\A$, boundary $\Gamma$, and interface $\I$, although best parameter choices would need to be rethought. Full adaptivity requires further, and more fundamental changes, in particular to the discretization of $\Gamma$, the definition of the annular region $\A$, and the solver used to invert the PDE on $\A$.

\section{Acknowledgments}
I owe many thanks to Mike Shelley, Shravan Veerapaneni, and Alex Barnett for a great many useful discussions, as well as Manas Rachh, for help with some of the boundary integral components and the biharmonic FMM code that underlies the Stokes solver. I'd also like to thank Dan Fortunato, for coining the rather clever term \emph{function intension}, and for his continued work on an adaptive version of this methodology, which we hope to present soon.

\bibliographystyle{unsrturl}
\bibliography{refs}

\end{document}